\documentclass{article}

\usepackage[utf8]{inputenc}
\usepackage[T1]{fontenc}
\usepackage{adjustbox}
\usepackage{amsmath}
\usepackage{amsthm}
\usepackage{array}
\usepackage{relsize}
\usepackage{scalerel}[2016/12/29]
\usepackage{wasysym}
\usepackage{mathpartir}
\usepackage{mathtools}
\usepackage{caption}
\usepackage{rotating}
\usepackage{tcolorbox}
\tcbuselibrary{
  skins,
  breakable
}
\usepackage{quiver}
\usepackage{titling} 
\usepackage{titlesec}
\usepackage{float}
\usepackage{graphics, graphicx, tikz}
\usepackage{tikz-cd}
\usetikzlibrary{
  arrows,
  arrows.meta,
  decorations.markings,
  shapes.geometric,
  trees
}
\usepackage{xspace}
\usepackage{xcolor}
\usepackage{geometry}
\usepackage[colorlinks=true, linktocpage=true]{hyperref}

\usepackage{mathpazo}
\usepackage[scaled=0.95]{helvet}
\usepackage{courier}
\DeclareFontEncoding{LS1}{}{}
\DeclareFontSubstitution{LS1}{stix2}{m}{n}
\DeclareSymbolFont{arrows2}{LS1}{stix2sf}{m}{it}
\SetSymbolFont{arrows2}{bold}{LS1}{stix2sf}{b}{it}
\DeclareSymbolFontAlphabet{\mathsfit}{arrows2}
\DeclareMathSymbol{\lefttail}{\mathrel}{arrows2}{"B2}
\DeclareMathSymbol{\righttail}{\mathrel}{arrows2}{"B3}
\DeclareMathSymbol{\leftdbltail}{\mathrel}{arrows2}{"B4}
\DeclareMathSymbol{\rightdbltail}{\mathrel}{arrows2}{"B5}

\titlespacing*{\subsection} {0pt}{0ex plus 1ex minus .2ex}{2ex plus .2ex}

\definecolor{grey}{rgb}{0.5,0.5,0.5}
\definecolor{mypurple}{rgb}{0.7,0,0.7}
\definecolor{links_color}{rgb}{0.4,0.1,0.1}
\hypersetup{
  linkcolor = links_color,
  citecolor = links_color
}

\title{Two equivalent descriptions of opetopes:\\
  in terms of zoom complexes and of partial orders}
\author{Louise \textsc{Leclerc}}

\counterwithin{figure}{section}

\newcounter{paragraphCounter}[section]
\renewcommand{\theparagraphCounter}{\thesection.\arabic{paragraphCounter}}

\newtheoremstyle{break}
  {}{}
  {}{}
  {\bfseries}{}
  {\newline}{}
\theoremstyle{break}

\newtheorem{rmk}[paragraphCounter]{Remark}
\newtheorem{exm}[paragraphCounter]{Example}
\newtheorem*{exm*}{Example}

\newenvironment{thm}[1][]{
  \refstepcounter{paragraphCounter}
  \begin{tcolorbox}
    [ enhanced
    , breakable
    , title =
        \textcolor{black}{
          \ifstrempty{#1}
            {Theorem~\theparagraphCounter}
            {Theorem~\theparagraphCounter\,:~#1}
        }
    , attach boxed title to top left = 
        { xshift = 15pt
        , yshift* = -\tcboxedtitleheight/2
        }
    , boxed title style = 
        { size = small
        , colback = purple!20
        , top = 1pt
        , bottom = 1pt
        }
    , colback = white
    , after skip = 10pt
    , before skip = 6pt
    ]
  }{\end{tcolorbox}}

\newenvironment{thm*}[1][]{
  \begin{tcolorbox}
    [ enhanced
    , breakable
    , title =
        \textcolor{black}{
          \ifstrempty{#1}
            {Theorem}
            {Theorem\,:~#1}
        }
    , attach boxed title to top left = 
        { xshift = 15pt
        , yshift* = -\tcboxedtitleheight/2
        }
    , boxed title style = 
        { size = small
        , colback = purple!20
        , top = 1pt
        , bottom = 1pt
        }
    , colback = white
    , after skip = 10pt
    , before skip = 6pt
    ]
  }{\end{tcolorbox}}

\newenvironment{defin}[1][]{
  \refstepcounter{paragraphCounter}
  \begin{tcolorbox}
    [ enhanced
    , breakable
    , title =
        \textcolor{black}{
          \ifstrempty{#1}
            {Definition~\theparagraphCounter}
            {Definition~\theparagraphCounter\,:~#1}
        }
    , attach boxed title to top left = 
        { xshift = 15pt
        , yshift* = -\tcboxedtitleheight/2
        }
    , boxed title style = 
        { size = small
        , colback = blue!30
        , top = 1pt
        , bottom = 1pt
        }
    , colback = white
    , after skip = 10pt
    , before skip = 6pt
    ]
  }{\end{tcolorbox}}

\newenvironment{prop}[1][]{
  \refstepcounter{paragraphCounter}
  \begin{tcolorbox}
    [ enhanced
    , breakable
    , title =
        \textcolor{black}{
          \ifstrempty{#1}
            {Proposition~\theparagraphCounter}
            {Proposition~\theparagraphCounter\,:~#1}
        }
    , attach boxed title to top left = 
        { xshift = 15pt
        , yshift* = -\tcboxedtitleheight/2
        }
    , boxed title style = 
        { size = small
        , colback = green!50!blue!30!
        , top = 1pt
        , bottom = 1pt
        }
    , colback = white
    , after skip = 10pt
    , before skip = 6pt
    ]
  }{\end{tcolorbox}}

\newenvironment{lem}[1][]{
  \refstepcounter{paragraphCounter}
  \begin{tcolorbox}
    [ enhanced
    , breakable
    , title =
        \textcolor{black}{
          \ifstrempty{#1}
            {Lemma~\theparagraphCounter}
            {Lemma~\theparagraphCounter\,:~#1}
        }
    , attach boxed title to top left = 
        { xshift = 15pt
        , yshift* = -\tcboxedtitleheight/2
        }
    , boxed title style = 
        { size = small
        , colback = blue!30
        , top = 1pt
        , bottom = 1pt
        }
    , colback = white
    , after skip = 10pt
    , before skip = 6pt
    ]
  }{\end{tcolorbox}}

\newcommand{\para}{
  \refstepcounter{paragraphCounter}
  \noindent\bfseries{\theparagraphCounter.}
  \normalfont
}

\newcommand{\refsec}[1]{\hyperref[#1]{Section \ref*{#1}}}
\newcommand{\refrmk}[1]{\hyperref[#1]{Remark \ref*{#1}}}
\newcommand{\refexm}[1]{\hyperref[#1]{Example \ref*{#1}}}
\newcommand{\refthm}[1]{\hyperref[#1]{Theorem \ref*{#1}}}
\newcommand{\refdefin}[1]{\hyperref[#1]{Definition \ref*{#1}}}
\newcommand{\refprop}[1]{\hyperref[#1]{Proposition \ref*{#1}}}
\newcommand{\reflem}[1]{\hyperref[#1]{Lemma \ref*{#1}}}

\newcommand{\branch}[1]{~\raisebox{-0.2px}{$\scaleobj{1.2}{\lefttail}_{#1}$}~}
\newcommand{\branched}[1]{~\raisebox{-0.2px}{$\scaleobj{1.2}{\righttail}_{#1}$}~}

\newcommand{\edge}{\mathrel{\hspace{0.5px}\raisebox{-1.5px}{\rotatebox{90}{\scaleobj{0.8}{\lefttail}}}}}
\newcommand{\N}{\ensuremath{\mathbb{N}}\xspace}

\newcommand*{\intInter}[2]
  {\ensuremath{
    \left[\!\left[\, \ensuremath{{#1}},\, \ensuremath{{#2}} \,\right]\!\right]
  }}
\renewcommand{\dim}{\ensuremath{\mathsf{dim}}}
\newcommand{\core}{\ensuremath{\mathsf{core}}}

\newcommand{\trigup}[1]{\ensuremath{\triangleleft^{\scriptscriptstyle \,C_{#1},\,+}}}

\newcommand{\msfdots}{\mathsf{dots}}
\newcommand{\msfwdots}{\mathsf{whitedots}}
\newcommand{\msfbdots}{\mathsf{blackdots}}
\newcommand{\msfndots}{\mathsf{nulldots}}
\newcommand{\msfleaves}{\mathsf{leaves}}

\newcommand{\iedge}{\mathrel{\raisebox{-5px}{\rotatebox{90}{$\scriptstyle\multimap$}}}}

\newcommand{\DFC}{\ensuremath{\mathbf{DFC}}}

\newcommand{\KJBM}{\ensuremath{\mathbf{KJBM}}}
\newcommand{\MOP}{\ensuremath{\mathbf{MOP}}}

\renewcommand{\o}{\circ}

\newcommand{\rom}[1]{\textup{\uppercase\expandafter{\romannumeral#1}}}

\makeatletter
\newcommand{\hhat}[1]{% 
\begingroup%
  \let\macc@kerna\z@%
  \let\macc@kernb\z@%
  \let\macc@nucleus\@empty%
  \hat{\mathchoice%
    {\raisebox{.2ex}{\vphantom{\ensuremath{\displaystyle #1}}}}%
    {\raisebox{.2ex}{\vphantom{\ensuremath{\textstyle #1}}}}%
    {\raisebox{.16ex}{\vphantom{\ensuremath{\scriptstyle #1}}}}%
    {\raisebox{.14ex}{\vphantom{\ensuremath{\scriptscriptstyle #1}}}}%
    \smash{\hat{#1}}}%
\endgroup%
}
\makeatother

\begin{document}

  \maketitle

\begin{abstract}
  We introduce in this paper a definition of (non necessarily positive) opetopes where faces are organised in a poset.
  Then we show that this description is equivalent to that given in terms of constellations by \textsc{Kock}, \textsc{Joyal}, \textsc{Batanin} and \textsc{Mascari} in \cite{Kock2010}.
\end{abstract}

\section{Introduction}
\label{sec:intro}

Opetopes have been introduced first by \textsc{Baez} and \textsc{Dolan} in \cite{baez1997} as a kind of pasting diagram for higher category theory. They consist of formal gluings of operations, which may have arbritrary arities.
They enable an unbiased approach to higher algebra and may be combinatorially encoded in several ways.

\begin{exm}\label{exm:four_opetope}
  Here is an illustration of a (non-positive) $4$-opetope $\omega$:\\[10px]
  % \documentclass[landscape]{article}
% \usepackage[utf8]{inputenc}
% \usepackage{graphics, graphicx, tikz, geometry}
% \usepackage{amsfonts}
% \usepackage{amssymb,amsmath,amsthm,stmaryrd,mathrsfs,wasysym}

% %%% Set the fonts
% \usepackage{mathpazo}
% \usepackage[scaled=0.95]{helvet}
% \usepackage{courier}
% \usepackage{adjustbox}
% \geometry{hmargin=0.5cm,vmargin=1cm}
% \usetikzlibrary{arrows, arrows.meta}
% \linespread{1.05} % Palatino looks better with this
% \definecolor{mypurple}{rgb}{0.7,0,0.7}
% \usetikzlibrary{calc}
%
% \begin{document}
%
  \adjustbox{scale=0.8, center}
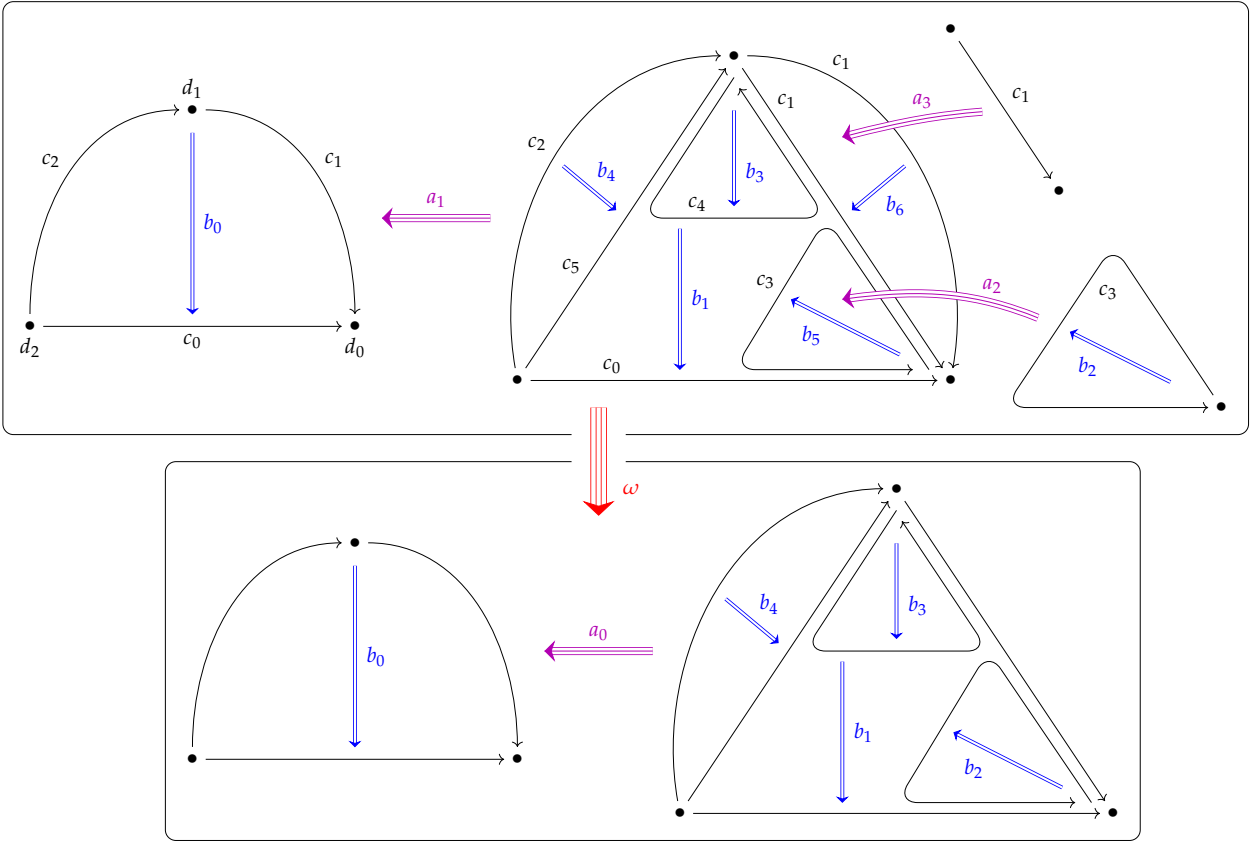
{
  \begin{tikzpicture}
    [
    , along/.style = {midway, inner sep = 0pt}
    , two/.style = {-{Stealth[length=3pt, width=6pt]}, double, double distance = 1pt, draw = blue, shorten >= 5pt, shorten <= 5pt}
    , three/.style={preaction={draw = mypurple, -{Stealth[length=5pt, width=10pt]}, double, double distance = 3pt}, draw = mypurple, -}
    , four/.style={preaction={draw = red, -{Stealth[length=7pt, width=15pt]}, double, double distance = 7pt, shorten >=0pt}, draw = red, -, double, double distance=2pt, shorten >=5pt}
    , scale = 0.9
    ]
    % \draw[help lines] (-12,-10) grid (14,8);

    % target de a1
    \node (d2_a0) at (-9, 1) {$\bullet$}
      node[below = 2] at (d2_a0) {$d_2$};
    \node (d1_a0) at (-6, 5) {$\bullet$}
      edge [<-, in = 90, out = 180]
        node[above left] {$c_2$}
        node[along] (c2_a0) {}
        (d2_a0)
        node[above = 2] at (d1_a0) {$d_1$};
    \node (d0_a0) at (-3, 1) {$\bullet$}
      edge [<-, in = 0, out = 90]
        node[above right] {$c_1$}
        node[along] (c1_a0) {}
        (d1_a0)
      edge [<-] 
        node[below] {$c_0$}
        node[along] (c0_a0) {}
        (d2_a0)
      node[below = 2] at (d0_a0) {$d_0$};
    \draw[two] (d1_a0) to
      node[right, color = blue, inner sep = 5pt] {$b_0$} 
      (c0_a0);

    % Source en a1
    \node (d2_a1)  at (0, 0)  {$\bullet$};
    \node (d1_a1)  at (4, 6)  {$\bullet$}
      edge [<-]
        node[above left = -1, pos = 0.7] {$c_5$}
        node[along] (c5_a1) {}
        (d2_a1)
      edge [<-, in = 100, out = 180]
        node[above left] {$c_2$}
        node[along] (c2_a1) {}
        (d2_a1);
    \node (d0_a1)  at (8, 0)   {$\bullet$}
      edge [<-]
        node[above right = -1, pos = 0.85] {$c_1$}
        node[along] (c1_int_a1) {}
        (d1_a1)
      edge [<-, in = 0, out = 80]
        node[above right = -1, pos = 0.8] {$c_1$}
        node[along] (c1_ext_a1) {}
        (d1_a1)
      edge [<-]
        node[above, pos = 0.8] {$c_0$}
        node[along] (c0_a1) {}
        (d2_a1);
    \pgfsetcornersarced{\pgfpoint{10}{10}}
    \draw 
          ($(d1_a1) - (0, 0.4)$) 
      --  ($(c5_a1) + (0.3, 0)$)
      --
        node[above, pos = 0.3] {$c_4$}
        node[along] (c4_a1) {}
          ($(c1_int_a1) - (0.3, 0)$) 
      --  ($(d1_a1) + (0.1, -0.6)$)
      [->];
    \draw 
          ($(d0_a1) + (-0.4, 0.2)$) 
      --  ($(c1_int_a1) - (0.3, 0)$)
      -- 
        node[above left = -1, pos = 0.5] {$c_3$}
        node[along] (c3_a1) {}
          ($(c0_a1) + (0, 0.2)$)
      --  ($(d0_a1) + (-0.7, 0.2)$)
      [->];
    \pgfsetcornersarced{\pgfpointorigin}
    \draw[two] (c2_a1) to
      node[color = blue, above right] {$b_4$}
      (c5_a1);
    \draw[two] (c1_ext_a1) to
      node[color = blue, below right] {$b_6$}
      (c1_int_a1);
    \draw[two, shorten <= 20] (d1_a1) to
      node[pos = 0.7, color = blue, right = 2] {$b_3$} 
      (c4_a1);
    \draw[two, shorten <= 20] (d0_a1) to
      node[pos = 0.7, color = blue, below left] {$b_5$} 
      (c3_a1);
    \draw[two] (3, 3) to
      node[color = blue, right = 2] {$b_1$} 
      (3,0);

    % Source en a2
    \node (d0_a2) at (13,-0.5) {$\bullet$};
    \pgfsetcornersarced{\pgfpoint{10}{10}}
    \draw 
          (d0_a2) 
      --  ($(d0_a2) + (-2, 3)$)
      --
        node[below right, pos = 0.2] {$c_3$}
        node[along] (c3_a2) {}
          ($(d0_a2) + (-4, 0)$)
      --  (d0_a2)
      [->];
    \pgfsetcornersarced{\pgfpointorigin}
    \draw[two, shorten <= 20] (d0_a2) to
      node[pos = 0.7, color = blue, below left] {$b_2$} 
      (c3_a2);

    % Source en a3
    \node (d1_a3) at (8, 6.5) {$\bullet$};
    \node (d0_a3) at (10, 3.5) {$\bullet$}
      edge [<-]
        node[above right] {$c_1$}
        node[along] (c1_a3) {}
        (d1_a3);

    % Target de omega
    \node (d2_t)  at (3, -8)  {$\bullet$};
    \node (d1_t)  at (7, -2)  {$\bullet$}
      edge [<-]
        node[along] (c5_t) {}
        (d2_t)
      edge [<-, in = 100, out = 180]
        node[along] (c2_t) {}
        (d2_t);
    \node (d0_t)  at (11, -8)   {$\bullet$}
      edge [<-] 
        node[along] (c1_t) {}
        (d1_t)
      edge [<-]
        node[along] (c0_t) {}
        (d2_t);
    \pgfsetcornersarced{\pgfpoint{10}{10}}
    \draw 
          ($(d1_t) - (0, 0.4)$) 
      --  ($(c5_t) + (0.3, 0)$)
      --
        node[along] (c4_t) {}
          ($(c1_t) - (0.3, 0)$) 
      --  ($(d1_t) + (0.1, -0.6)$)
      [->];
    \draw 
          ($(d0_t) + (-0.4, 0.2)$) 
      --  ($(c1_t) - (0.3, 0)$)
      --  
        node[along] (c3_t) {}
          ($(c0_t) + (0, 0.2)$)
      --  ($(d0_t) + (-0.7, 0.2)$)
      [->];
    \pgfsetcornersarced{\pgfpointorigin}
    \draw[two] (c2_t) to
      node[color = blue, above right] {$b_4$} 
      (c5_t);
    \draw[two, shorten <= 20] (d1_t) to
      node[pos = 0.7, color = blue, right = 2] {$b_3$} 
      (c4_t);
    \draw[two, shorten <= 20] (d0_t) to
      node[pos = 0.7, color = blue, below left] {$b_2$} 
      (c3_t);
    \draw[two] (6, -5) to
      node[color = blue, right = 2] {$b_1$} 
      (6,-8);

    % Target du target
    \node (d2_tt) at (-6, -7) {$\bullet$};
    \node (d1_tt) at (-3, -3) {$\bullet$}
      edge [<-, in = 90, out = 180] 
        node[along] (c2_tt) {}
        (d2_tt);
    \node (d0_tt) at (0, -7) {$\bullet$}
      edge [<-, in = 0, out = 90]
        node[along] (c1_tt) {}
        (d1_tt)
      edge [<-] 
        node[along] (c0_tt) {}
        (d2_tt);
    \draw[two] (d1_tt) to
      node[color = blue, right = 2] {$b_0$} 
      (c0_tt);

    % ------- Triple arrows -------

    % a1 vers a0
    \draw[three] (-0.5, 3) to 
      node[above, color = mypurple, inner sep = 5pt] {$a_1$}
      (-2.5, 3);
    % a2 vers a1
    \draw[three, shorten <= 10] (c3_a2) to[bend right = 15]
      node[above, color = mypurple, inner sep = 5pt, pos = 0.3] {$a_2$}
      (6, 1.5);
    % a3 vers a1
    \draw[three, shorten <= 10] (c1_a3) to[bend right = 5]
      node[above, color = mypurple, inner sep = 5pt] {$a_3$}
      (6, 4.5);
    % t vers tt
    \draw[three] (2.5, -5) to
      node[above, color = mypurple, inner sep = 5pt] {$a_0$}
      (0.5, -5);

    capsules
    \pgfsetcornersarced{\pgfpoint{5pt}{5pt}}
    \draw (1, -1) -- (-9.5, -1) -- (-9.5, 7) -- (13.5, 7) -- (13.5, -1) -- (2, -1);
    \draw (1, -1.5) -- (-6.5, -1.5) -- (-6.5, -8.5) -- (11.5, -8.5) -- (11.5, -1.5) -- (2, -1.5);
    \pgfsetcornersarced{\pgfpointorigin}

    \draw[four] (1.5, -0.5) to
      node[right, color = red, right = 8pt, pos = 0.75] {$\omega$}
      (1.5, -2.5);

  \end{tikzpicture}  
  }
%
% \end{document}
\end{exm}

A lengthier presentation of the concept may be found in \cite{ChengLauda2004Guidebook}, \cite{leinster2003} or (for the specific case of positive opetopes) in \cite{leclerc2023}. We assume the reader to be already familiar with the underlying idea of opetopes and present now the subject of this article.

First of all, a new definition of (non necessarily positive) opetopes will be introduced in section \ref{sec:DFC}, under the name dendritic face complexes (DFC). The idea is to store every cell of an opetope in a graded poset, and to characterise the posets which arise from opetopes. This approach is really akin to the one taken by \textsc{Zawadowski} in \cite{zawadowski2023positive} or \textsc{hadzihasanovic} in \cite{hadzihasanovic2019combinatorialtopological}.
It is an extension of the formalism already presented in \cite{leclerc2023} to the broader context of non-positive opetopes.
They form a category $\DFC$, which restricts to a groupoid $\DFC^\core$ when considering only the isomorphisms.

Then we give in section \ref{sec:zoom} a brief description of zoom complexes and the definition of opetopes introduced by \textsc{Kock}, \textsc{Joyal}, \textsc{Batanin} and \textsc{Mascari} in \cite{Kock2010}. We also denote their groupoid $\KJBM$.

Finally, we introduce two functors between those two categories (sections \ref{sec:DFC_to_Zoom} and \ref{sec:Zoom_to_DFC}) and show that the two differents point of view coincide modulo an equivalence of categories $\DFC^\core \simeq \KJBM$ established in section \ref{sec:cat_equiv}. This result is similar to the duality proven by \textsc{Zawadowski} in his work on positive opetopes \cite{zawadowski2023duality}.

\section{Dendritic face complexes}
\label{sec:DFC}

The poset structure that will be introduced below makes implicit use of a notion of rooted tree, which we start by defining below.
\begin{defin}[Rooted tree]\label{defin:rooted_tree}
  A \emph{rooted tree} $T$ consists of a finite set of \emph{nodes} $T^\bullet$, and a finite set of \emph{edges} $T^{\edge}$ such that:
  \begin{itemize}
    \item each node $a$ must be paired with exactly one \emph{target edge} $b$, which we denote $\branch{b} a$. We also say that $a$ is the \emph{source node} of $b$. Sourceless edges are called \emph{leaves} of $T$, and their set is written $T^\vert$.
    \item each edge $b$ may be paired with at most one \emph{target node} $a$, we denote $a \branch{b}$ such a relation. We also say that $b$ is a \emph{source} of $a$.
  \end{itemize}
  We moreover ask for a distinguished targetless edge $\rho(T) \in T^{\edge}$, called the \emph{root edge} of $T$, satisfying the following property:
  for each edge $b\in T^{\edge}$, there is \emph{descending path} in $T$
  $$
  \branched{b = b_1} a_1 \branched{b_2} a_2 \branched{b_2} \cdots \branched{b_p = \rho(T)}
  $$
  from $b$ to the root of $T$.
  More generally, an alternate sequence of nodes and edges $s = (x_1,\,x_2,\,\cdots,\,x_p)$ will be called
  \begin{itemize}
    \item an \emph{increasing path} if $\forall i<p$, $x_{i+1}$ is the target of $x_i$.
    \item a \emph{decreasing path} if $\forall i<p$, $x_i$ is the target of $x_{i+1}$.
    \item a \emph{zig-zag} if for any pair $(x_i,\, x_{i+1})$ where $x_i$ is an edge, we have $x_{i+1}$ source of $x_i$ and $x_{i-1}$ target of $x_i$ if $i>0$, or $x_{i+1}$ target of $x_i$ and $x_{i-1}$ source of $x_i$ if $i>0$.
      the zig-zag is said to be \emph{simple} when $\forall i\neq j,\, x_i \neq x_j$. A simple zig-zag corresponds to a succession of a decreasing and an increasing path, or vice-versa.
  \end{itemize}
  A \emph{subdivision} for $T$ consists of a finite totally ordered set $W(b)$ for each edge $b \in T^{\edge}$.
  A subdivided tree may always be seen as a rooted tree, by replacing each edge $b$ by a linear tree
  $$\branch{b_0} a_1 \branch{b_1} a_2 \branch{b_2} \cdots \branch{b_{p-1}} a_p \branch{b_p}$$
  where $a_1,\,\cdots,\,a_p$ are the elements of $W(b)$ following the increasing order.
  
  Following \cite{Kock2010}, we let 
  \begin{mathpar}
    \msfbdots(T) = T^\bullet
    \and
    \msfwdots(T) = \bigsqcup_{b \in T^{\edge}}{W(b)}
    \and
    \msfdots(T) = \msfbdots(T) \sqcup \msfwdots(T)
    \and
    \msfndots(T) = \{a \in T^\bullet \mid a \mbox{ sourceless}\}
    \and
    \msfleaves(T) = T^\vert = \{b \in T^{\edge} \mid b \mbox{ sourceless}\}
    \and
    T^{\iedge} = T^{\edge} \setminus T^\vert
    .
  \end{mathpar}
  A rooted tree $T$ is called a \emph{unit tree} iff it is has no node and only one leaf: its root. It is called a \emph{corolla} if it has exactly one node (and hence each edge is either a leaf or the root).
\end{defin}

\begin{rmk}\label{rmk:epi_unique_root}
  Notice that the root is uniquely determined as the unique edge without a target node, and there is exactly one path as above.
  Notice also that for any node $a$ there is a unique descending path from $a$ to $\rho(\omega)$.
\end{rmk}

\begin{exm}\label{exm:rooted_tree}
  In \ref{fig:rooted_tree_exm} is a representation of the rooted tree $T$ having
  \begin{itemize}
    \item as nodes $T^\bullet := \{a_1,\,a_2,\,a_3,\,a_4\}$, as edges $T^{\edge} := \{b_1,\,b_2,\,b_3,\,b_4,\,b_5\}$, and as root $b_1$.
    \item as relations
      \vspace{-5pt}
      \begin{mathpar}
        \branch{b_1} a_1
        \and
        a_1 \branch{b_2}
        \and
        a_1 \branch{b_3} a_2
        \and
        a_2 \branch{b_4} a_3
        \and
        a_2 \branch{b_5} a_4
      \end{mathpar}
  \end{itemize}
  \vspace{-15pt}
  \begin{figure}[H]
    \begin{minipage}{.5\textwidth}
      \[\begin{tikzpicture}
        [ level distance = 80px
        , bdot/.style = {circle, fill}
        , boundary/.style = {draw = none}
        , scale = 0.8
        ]
        \pgfmathsetmacro{\theta}{45}
        \node (root_3) {} [grow' = up, boundary, level distance = 50px]
          child {
            node (a1) [bdot, label=left:{\large$a_1$}] {} [clockwise from = 90+0.5*\theta, sibling angle = \theta]
              child[level distance = 80px]{
                node [boundary] {}
                edge from parent node [pos = 0.8, left] {\large$b_2$}
              }
              child[level distance = 60px]{
                node (a2) [bdot, label={left}:{\large$a_2$}] {} [clockwise from = 112.5, sibling angle = 45]
                  child[level distance = 50px]{
                    node (a3) [bdot, label={above}:{\large$a_3$}] {}
                    edge from parent node [left] {\large$b_4$}
                  }
                  child[level distance = 50px]{
                    node (a4) [bdot, label={above}:{\large$a_4$}] {}
                    edge from parent node [right] {\large$b_5$}
                  }
                edge from parent node [right] {\large$b_3$}
              }
            edge from parent node [right, pos = 0.2] {\large$b_1$}
          };
      \end{tikzpicture}\]
      \captionof{figure}{A rooted tree $T$}
      \label{fig:rooted_tree_exm}
    \end{minipage}
    \begin{minipage}{.5\textwidth}
      \[\begin{tikzpicture}
        [ level distance = 80px
        , bdot/.style = {circle, fill}
        , boundary/.style = {draw = none}
        , scale = 0.8
        ]
        \pgfmathsetmacro{\theta}{45}
        \node (root_3) {} [grow' = up, boundary, level distance = 50px]
          child {
            node (a1) [bdot, label=left:{\large$a_1$}] {} [clockwise from = 90+0.5*\theta, sibling angle = \theta]
              child[level distance = 120px]{
                node [boundary] {}
                edge from parent 
                  node [pos = 0.85, left] {\large$b_2$}
                  node (w1) [pos = 0.2] {}
                  node (w2) [pos = 0.4] {}
                  node (w3) [pos = 0.6] {}
              }
              child[level distance = 100px]{
                node (a2) [bdot, label={left}:{\large$a_2$}] {} [clockwise from = 112.5, sibling angle = 45]
                  child[level distance = 50px]{
                    node (a3) [bdot, label={above}:{\large$a_3$}] {}
                    edge from parent node [left] {\large$b_4$}
                  }
                  child[level distance = 50px]{
                    node (a4) [bdot, label={above}:{\large$a_4$}] {}
                    edge from parent node [right] {\large$b_5$}
                  }
                edge from parent 
                  node [pos = 0.8, right] {\large$b_3$}
                  node (w4) [pos = 0.25] {}
                  node (w5) [pos = 0.5] {}
              }
            edge from parent node [right, pos = 0.2] {\large$b_1$}
          };
        \draw (w1) circle (0.2)
          node [left = 4] {$w_1$};
        \draw (w2) circle (0.2)
          node [left = 4] {$w_2$};
        \draw (w3) circle (0.2)
          node [left = 4] {$w_3$};
        \draw (w4) circle (0.2)
          node [right = 4] {$w_4$};
        \draw (w5) circle (0.2)
          node [right = 4] {$w_5$};
      \end{tikzpicture}\]
      \captionof{figure}{A subdivison $T'$ of $T$}
      \label{fig:subdivided_rooted_tree_exm}
    \end{minipage}
  \end{figure}

  In \ref{fig:subdivided_rooted_tree_exm} is a representation of the subdivided tree $T' = (T,\, W)$, where
  \begin{mathpar}
    W(b_2) = (w_1 < w_2 < w_3)
    \and
    W(b_3) = (w_4 < w_5)
    \and
    W(b_1) = W(b_4) = W(b_5) = \emptyset 
  \end{mathpar}
\end{exm}

\begin{defin}[Many-to-one poset]\label{defin:MOP}
  A many-to-one poset (MOP) consists of:
  \begin{itemize}
    \item a finite set of elements $P$ (often called \emph{cells}),
    \item a gradation $\dim : P \to \N \cup \{-1\}$. We let $P_k := \dim^{-1}(\{k\})$ for any $k$,
    \item for each $k\ge -1$, two applications $\delta_k,\,\gamma_k : P_{k+1} \to \mathcal{P}(P_k)$ (we will omit the index $k$ for brievety),
  \end{itemize}
  such that the following properties hold:
  \begin{itemize}
    \item $\forall x \in P,\quad \dim(x)\ge0 \rightarrow \gamma(x) \mbox{ is a singleton}$.
    \item $\forall x \in P,\quad \delta(x) \cap \gamma(x) \neq \emptyset \rightarrow \delta(x) = \gamma(x)$.
    \item $\exists!(x=:\ast),\,\dim(x)=-1$,\, and \,$\forall x,\, \dim(x) = 0 \to \gamma(x)\setminus\delta(x) = \{\ast\}$.
  \end{itemize}

  We write \vspace{-4pt}
    \begin{mathpar}
      y \prec^{-} x \mbox{ \,iff\, } y \in \delta(x)\setminus\gamma(x) \and y \prec^{+} x \mbox{ \,iff\, } y \in \gamma(x)\setminus\delta(x) \\
      y \prec^{\o} x \mbox{ \,iff\, } y \in \delta(x)\cap\gamma(x) \and y \prec x \mbox{ \,iff\, } y \in \delta(x)\cup\gamma(x) \\
      \partial(x) = \delta(x) \cup \gamma(x) \and \delta^{-}(x) := \delta(x) \setminus \gamma(x) \and \gamma^{+}(x) = \gamma(x) \setminus \delta(x) \and \Lambda = P \setminus \gamma^{+}(P)
    \end{mathpar}

  In particular: $\prec,\,\prec^{-},\,\prec^{+}$ and $\prec^{\o}$ are asymmetric, and the reflexive transitive closure of $\prec$ equips $P$ with a structure of partially ordered set, such that $\dim$ is an increasing map.\\
  For $x$ with $\dim(x)\ge1$, we sometimes identify $\gamma(x)$ with its unique element, which we call the \emph{target} of $x$.
  The elements of $\delta(x)$ are called the \emph{sources} of $x$. If $y \prec^{\o} x$, $x$ is said to be \emph{a loop} on $y$.\\
  For $k\in\N \cup \{-1\}$, we also denote 
  \begin{mathpar}
    P_{\ge k} := \bigcup_{i\ge k}P_i \and P_{>k} := \bigcup_{i>k}P_i \and \dim(P) := \max\{\dim(x)\}_{x\in P} \mbox{ the \emph{dimension} of }P.
  \end{mathpar}

  \quad Finally for each $x \in \Lambda$ and $z \in \partial^2(x)$, we ask for a total order $\triangleleft$ on the set $\{y \mid z \prec^{\o} y \prec^{-} x\}$, called the \emph{local order}.
\end{defin}

\begin{rmk}\label{rmk:loop_facets}
  When $y \prec^{\o} x$ in a MOP, we have $\delta(x) = \gamma(x) = \{y\}$.
  Especially there is no $y' \prec^{\alpha} x$ for some $\alpha \neq \o$.
  In the following, we make implicit use of the monoid structure on $\{+,\,-,\,\o\}$, seen as $\{+1,\,-1,\,0\} \subseteq (\mathbb{Z},\,\times)$.
\end{rmk}

\para From now on, we will use a kind of \textsc{Hasse} diagrams to represent MOP.
The convention will be as follows:
% https://q.uiver.app/#q=WzAsOSxbMCwwLCJ4Il0sWzAsMSwieSJdLFswLDIsInhcXHByZWNeey19eSJdLFsyLDEsInkiXSxbMiwwLCJ4Il0sWzIsMiwieCBcXHByZWNeeyt9eSJdLFs0LDAsIngiXSxbNCwxLCJ5Il0sWzQsMiwieFxccHJlY157XFxhbHBoYX15Il0sWzEsMCwiIiwwLHsic3R5bGUiOnsiaGVhZCI6eyJuYW1lIjoibm9uZSJ9fX1dLFs0LDMsIiIsMCx7InN0eWxlIjp7ImJvZHkiOnsibmFtZSI6ImJhcnJlZCJ9LCJoZWFkIjp7Im5hbWUiOiJub25lIn19fV0sWzYsNywiXFxhbHBoYSIsMSx7InN0eWxlIjp7ImhlYWQiOnsibmFtZSI6Im5vbmUifX19XV0=
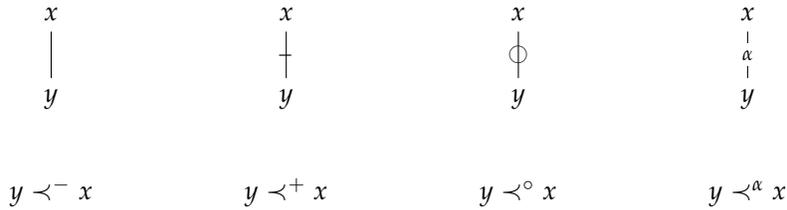
\[\begin{tikzcd}[]
  x && x && x && x \\
  y && y && y && y \\
  {y\prec^{-}x} && {y\prec^{+}x} && {y\prec^{\o}x} && {y\prec^{\alpha}x}
  \arrow[no head, from=2-1, to=1-1]
  \arrow["\shortmid"{marking}, no head, from=1-3, to=2-3]
  \arrow["\Circle"{marking}, no head, from=1-5, to=2-5]
  \arrow["\alpha"{description}, no head, from=1-7, to=2-7]
\end{tikzcd}\]
The idea for encoding opetopes (such as depicted in \refexm{exm:four_opetope}) as MOP is that they have cells, which will be the elements of the poset $P$, and the gradation $\dim$ really is the geometric dimension of those cells.
Then we will ask for a relation $y \prec x$ iff $y$ is a \emph{facet} of $x$. That is, a codimension one subface of $x$.
Since our opetopes has oriented faces (cells has sources and target), we moreover label the relations $y \prec x$ as
\begin{itemize}
  \item $y \prec^{+} x$ if $y$ is the target of $x$, but not among its source.
  \item $y \prec^{-} x$ if $y$ is the source of $x$, but not its target.
  \item $y \prec^{\o} x$ if $y$ is simultaneously the source and the target of $x$. That is, when $y$ is a \emph{loop} on $x$.
\end{itemize}
Of course, a MOP remember the data used to build an opetope, and how to reconstruct it, but fails to capture the properties of such objects.
This is why the refined version of \refdefin{defin:DFC} will be introduced later.

\begin{exm}\label{exm:opetope_loops}
  We depict below a $3$-dimensional opetope $\varrho$, who has five loops $b_3,\,b_4,\,b_5,\,b_6$ and $b_8$ on the same $0$-cell $c_1$.
  
  \vspace{15pt}
  % \documentclass[landscape]{article}
% \usepackage[utf8]{inputenc}
% \usepackage{graphics, graphicx, tikz, geometry}
% \usepackage{amsfonts}
% \usepackage{amssymb,amsmath,amsthm,stmaryrd,mathrsfs,wasysym}
%
% %%% Set the fonts
% \usepackage{mathpazo}
% \usepackage[scaled=0.95]{helvet}
% \usepackage{courier}
% \usepackage{adjustbox}
% \geometry{hmargin=0.5cm,vmargin=1cm}
% \usetikzlibrary{arrows, arrows.meta}
% \linespread{1.05} % Palatino looks better with this
% \definecolor{mypurple}{rgb}{0.7,0,0.7}
% \usetikzlibrary{calc}
%
% \begin{document}
%
  \adjustbox{scale=1, center}
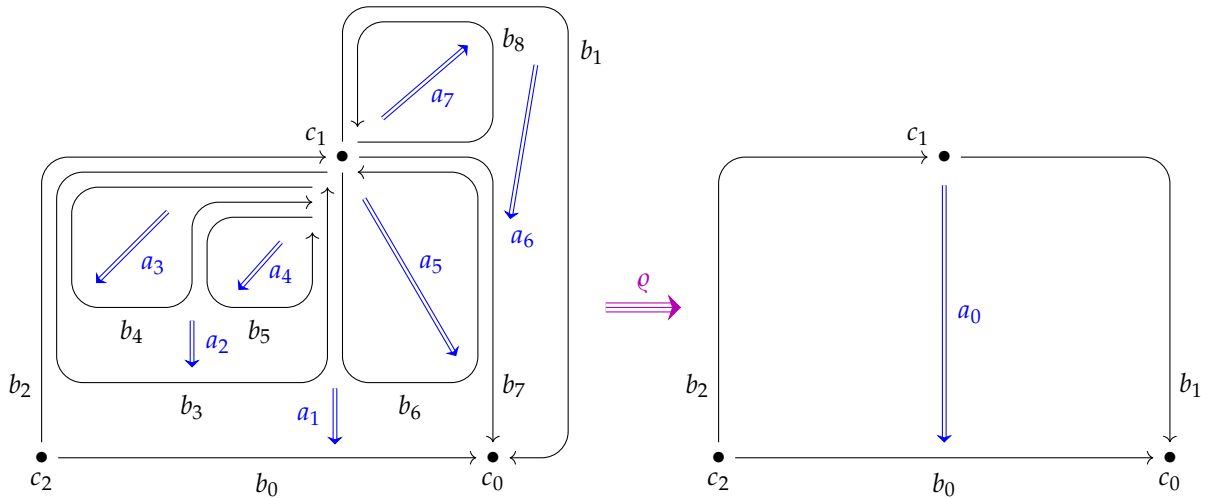
{
  \begin{tikzpicture}
    [
    , along/.style = {midway, inner sep = 0pt}
    , two/.style = {-{Stealth[length=3pt, width=6pt]}, double, double distance = 1pt, draw = blue, shorten >= 5pt, shorten <= 5pt}
    , three/.style={preaction={draw = mypurple, -{Stealth[length=5pt, width=10pt]}, double, double distance = 3pt}, draw = mypurple, -}
    , four/.style={preaction={draw = red, -{Stealth[length=7pt, width=15pt]}, double, double distance = 7pt, shorten >=0pt}, draw = red, -, double, double distance=2pt, shorten >=5pt}
    , scale = 1
    ]
    % \draw[help lines] (-12,-10) grid (14,8);

    % source
    \node (c2) at (-8, 0) {$\bullet$}
      node[below = 2] at (c2) {$c_2$};
    \node (c1) at (-4, 4) {$\bullet$}
      node[above = 8, left = 2] at (c1) {$c_1$};
    \node (c0) at (-2, 0) {$\bullet$}
      node[below = 2] at (c0) {$c_0$};
    \pgfsetcornersarced{\pgfpoint{10}{10}}
    \draw
          (c2)
      --  node[left, pos = 0.2] (b2) {$b_2$}
          (-8, 4)
      --  (c1)
      [->];
    \draw
          (c1)
      --  (-4, 6)
      --  (-1, 6)
      --  node[right = 1, pos = 0.1] (b1) {$b_1$} 
          (-1, 0)
      --  (c0)
      [->];
    \draw
          (c1)
      --  (-2, 4)
      --  node[right, pos = 0.8] (b7) {$b_7$} 
          (c0)
      [->];
    \draw (c2) -- node[along] (b0) {} (c0) [->];
    \node at (b0) [below = 2] {$b_0$};
    \draw
          ($(c1)-(0.2, 0.2)$)
      --  (-7.8, 3.8)
      --  (-7.8, 1)
      --  node[along] (b3) {}
          (-4.2, 1)
      --  (-4.2, 3.6)
      [->];
    \node at (b3) [below = 1] {$b_3$};
    \draw
          ($(c1)-(0.4, 0.4)$)
      --  (-7.6, 3.6)
      --  (-7.6, 2)
      --  node[along] (b4) {}
          (-6, 2)
      --  (-6, 3.4)
      --  (-4.4, 3.4)
      [->];
    \node at (b4) [below = 1] {$b_4$};
    \draw
          ($(c1)-(0.4, 0.8)$)
      --  (-5.8, 3.2)
      --  (-5.8, 2)
      --  node[along] (b5) {}
          (-4.4, 2)
      --  (-4.4, 3)
      [->];
    \node at (b5) [below = 1] {$b_5$};
    \draw
          (c1)
      --  (-4, 1)
      --  node[along] (b6) {}
          (-2.2, 1)
      --  (-2.2, 3.8)
      --  (-3.8, 3.8)
      [->];
    \node at (b6) [below = 1] {$b_6$};
    \draw
          ($(c1)+(0.2, 0.2)$)
      --  (-2, 4.2)
      --  node[along] (b8) {} 
          (-2, 5.8)
      --  (-3.8, 5.8)
      --  (-3.8, 4.4)
      [->];
    \node at (b8) [right = 8, above = 8] {$b_8$};

    \draw[two] ($(b4)+(0.6,1.4)$) to
      node[below right, color = blue] {$a_3$} 
      ($(b4)-(0.6,-0.2)$);
    \draw[two] ($(b5)+(0.4,1)$) to
      node[below = 3, right, color = blue] {$a_4$} 
      ($(b5)-(0.4,-0.1)$);
    \draw[two] ($(b3)+(0,1)$) to
      node[right, color = blue, inner sep = 5pt] {$a_2$} 
      (b3);
    \draw[two] (-3.8, 3.6) to
      node[above = 5, right = 0, color = blue] {$a_5$} 
      (-2.4, 1.2);
    \draw[two] (-3.6, 4.4) to
      node[below = 6, right = -3, color = blue, inner sep = 5pt] {$a_7$} 
      (-2.2, 5.6);
    \draw[two] (-4.1, 1.1) to
      node[left, color = blue, inner sep = 5pt] {$a_1$} 
      (-4.1, 0);
    \draw[two] (-1.4, 5.4) to
      node[below = 28, color = blue, inner sep = 5pt] {$a_6$} 
      (-1.8, 3);

    % target
    \node (c2_t) at (1, 0) {$\bullet$}
      node[below = 2] at (c2_t) {$c_2$};
    \node (c1_t) at (4, 4) {$\bullet$}
      node[above = 8, left = 2] at (c1_t) {$c_1$};
    \node (c0_t) at (7, 0) {$\bullet$}
      node[below = 2] at (c0_t) {$c_0$};
    \pgfsetcornersarced{\pgfpoint{10}{10}}
    \draw
          (c2_t)
      --  node[left, pos = 0.2] (b2_t) {$b_2$}
          (1, 4)
      --  (c1_t)
      [->];
    \draw
          (c1_t)
      --  (7, 4)
      --  node[right, pos = 0.8] (b1_t) {$b_1$} 
          (c0_t)
      [->];
    \draw (c2_t) -- node[along] (b0_t) {} (c0_t) [->];
    \node at (b0_t) [below = 2] {$b_0$};
    \draw[two] ($(b0_t) + (0, 3.8)$) to
      node[right, color = blue, inner sep = 5pt] {$a_0$}
      (b0_t);
    
    % ------- Triple arrows -------

    \draw[three] (-0.5, 2) to
      node[above, color = mypurple, inner sep = 5pt] {$\varrho$}
      (0.5, 2);

  \end{tikzpicture}
  }
%
% \end{document}

  Which we may encode as the following MOP.
  \vspace{15pt}
  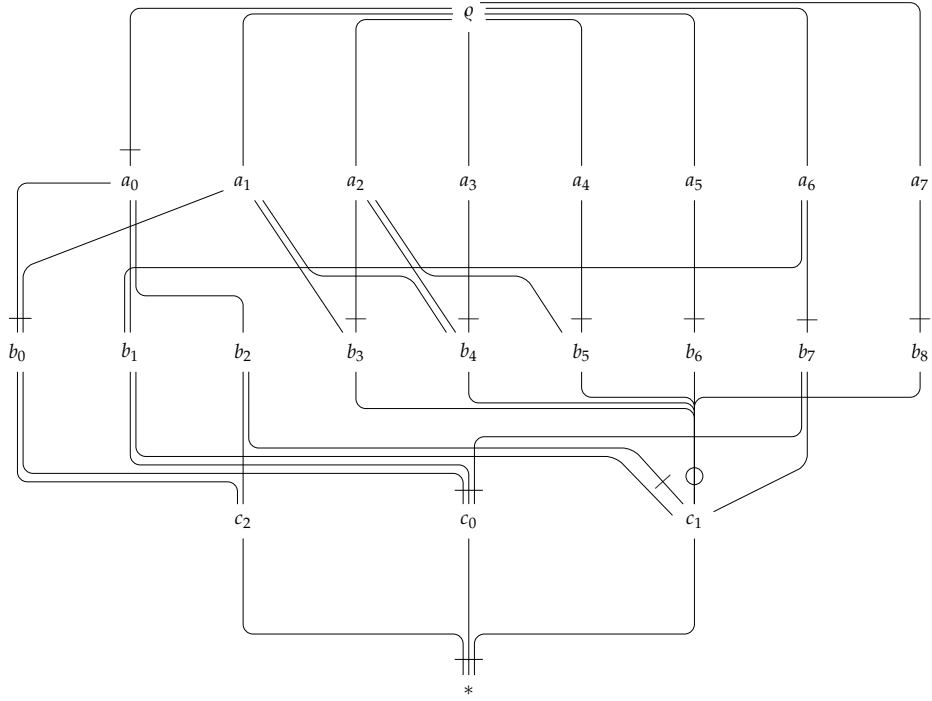
\begin{figure}[H]
    % \documentclass[landscape]{article}
% \usepackage[utf8]{inputenc}
% \usepackage{graphics, graphicx, tikz, geometry}
% \usepackage{amsfonts}
% \usepackage{amssymb,amsmath,amsthm,stmaryrd,mathrsfs,wasysym}
%
% %%% Set the fonts
% \usepackage{mathpazo}
% \usepackage[scaled=0.95]{helvet}
% \usepackage{courier}
% \usepackage{adjustbox}
% \geometry{hmargin=0.5cm,vmargin=1cm}
% \usetikzlibrary{arrows, arrows.meta}
% \linespread{1.05} % Palatino looks better with this
% \definecolor{mypurple}{rgb}{0.7,0,0.7}
% \usetikzlibrary{calc}
%
% \begin{document}
%
  \adjustbox{scale=0.75, center}{

    \pgfmathsetmacro{\eps}{0.2}
    \begin{tikzpicture}
      [ scale = 0.5
      , every node/.style = {outer sep = 2}]
      \node[draw = none] (xeps) at (\eps, 0) {};
      \node[draw = none] (yeps) at (0, \eps) {};

      % --- Noeuds ---
      
      % dim 3
      \node (rho) at (16, 24) {$\varrho$};

      % dim 2
      \node (a0) at (4, 18) {$a_0$};
      \node (a1) at (8, 18) {$a_1$};
      \node (a2) at (12, 18) {$a_2$};
      \node (a3) at (16, 18) {$a_3$};
      \node (a4) at (20, 18) {$a_4$};
      \node (a5) at (24, 18) {$a_5$};
      \node (a6) at (28, 18) {$a_6$};
      \node (a7) at (32, 18) {$a_7$};

      % dim 1
      \node (b0) at (0, 12) {$b_0$};
      \node (b1) at (4, 12) {$b_1$};
      \node (b2) at (8, 12) {$b_2$};
      \node (b3) at (12, 12) {$b_3$};
      \node (b4) at (16, 12) {$b_4$};
      \node (b5) at (20, 12) {$b_5$};
      \node (b6) at (24, 12) {$b_6$};
      \node (b7) at (28, 12) {$b_7$};
      \node (b8) at (32, 12) {$b_8$};

      % dim 0
      \node (c2) at (8, 6) {$c_2$};
      \node (c0) at (16, 6) {$c_0$};
      \node (c1) at (24, 6) {$c_1$};

      % dim -1
      \node (ast) at (16, 0) {$\ast$};

      % --- relations ---
      \pgfsetcornersarced{\pgfpoint{5}{5}}

      % dim 2 & 3
      \draw 
        ($(rho)+(-0.6, 0.2)$)
          -- (4, 24.2)
          -- node[sloped, pos = 0.9] {$\mid$} (a0)
        (rho)
          -- (8, 24)
          -- (a1)
        ($(rho)+(-0.6, -0.2)$)
          -- (12, 23.8)
          -- (a2)
        (rho)
          -- (a3)
        ($(rho)+(0.6, -0.2)$)
          -- (20, 23.8)
          -- (a4)
        (rho)
          -- (24, 24)
          -- (a5)
        ($(rho)+(0.6, 0.2)$)
          -- (28, 24.2)
          -- (a6)
        ($(rho)+(0.4, 0.4)$)
          -- (32, 24.4)
          -- (a7)
        ;

      % dim 1 & 2
      \draw % enfants de a0
        (a0)
          -- (0, 18)
          -- (b0)
        (a0)
          -- (b1)
        ($(a0.south) + (xeps)$)
          -- (4.2, 14)
          -- (8, 14)
          -- (b2)
        ;
      \draw % enfants de a1
        (a1)
          -- (0.2, 15)
          -- ($(b0.north)+(xeps)$)
        (a1)
          -- (b3)
        ($(a1.south east)$)
          -- (10.5, 14.7)
          -- (13.9, 14.7)
          -- ($(b4.north west) - (0.1, 0)$)
        ;
      \draw % b0 target de a0 et a1
        (-0.3, 13.2) -- (0.5, 13.2)
        ;
      \draw % enfants de a2
          (a2)
            -- node[sloped, pos = 0.9] {$\mid$} (b3)
          (a2)
            -- (b4)
          ($(a2.south east)$)
            -- (14.5, 14.7)
            -- (18, 14.7)
            -- ($(b5.north west)$)
        ;
      \draw % enfants de a3
        (a3)
          -- node[sloped, pos = 0.9] {$\mid$} (b4)
        ;
      \draw % enfants de a4
        (a4)
          -- node[sloped, pos = 0.9] {$\mid$} (b5)
        ;
      \draw % enfants de a5
        (a5)
          -- node[sloped, pos = 0.9] {$\mid$} (b6)
        ;
      \draw % enfants de a6
        (a6)
          -- node[sloped, pos = 0.9] {$\mid$} (b7)
        ($(a6.south)-(xeps)$)
          -- (27.8, 15)
          -- (3.8, 15)
          -- ($(b1.north)-(xeps)$)
        ;
      \draw % enfants de a7
        (a7)
          -- node[sloped, pos = 0.9] {$\mid$} (b8)
        ;

      % dim 0 & 1
      \draw % enfants de b0
        (b0)
          -- (0, 7.4)
          -- (7.8, 7.4)
          -- ($(c2.north)-(xeps)$)
        ($(b0.south)+(xeps)$)
          -- (0.2, 7.7)
          -- (15.8, 7.7)
          -- ($(c0.north)-(xeps)$)
        ;
      \draw % enfants de b1
        (b1)
          -- (4, 8)
          -- (16, 8)
          -- (c0)
        ($(b1.south)+(xeps)$)
          -- (4.2, 8.3)
          -- (21.2, 8.3)
          -- ($(c1.north west) - (0.1, 0.4)$)
        ;
      \draw % enfants de b2
        (b2)
          -- (c2)
        ($(b2.south)+(xeps)$)
          -- (8.2, 8.6)
          -- (21.8, 8.6)
          -- node[sloped, pos = 0.6] {$\mid$} ($(c1.north)-2*(xeps)$)
        ;
      \draw % enfants de b3
        (b3)
          -- (12, 10)
          -- (24, 10)
          -- (c1)
        ;
      \draw % enfants de b4
        (b4)
          -- (16, 10.2)
          -- (24, 10.2)
          -- (c1)
        ;
      \draw % enfants de b5
        (b5)
          -- (20, 10.4)
          -- (24, 10.4)
          -- (c1)
        ;
      \draw % enfants de b6
        (b6)
          -- (c1)
        ;
      \draw % enfants de b7
        (b7)
          -- (28, 8)
          -- (c1)
        ($(b7.south)-(xeps)$)
          -- (27.8, 9)
          -- (16.2, 9)
          -- ($(c0.north)+(xeps)$)
        ;
      \draw % enfants de b6
        (b8)
          -- (32, 10.4)
          -- (24, 10.4)
          -- (c1)
        ;
      \draw (24, 7.6) circle (0.3);
      \draw % c0 target de b0, b1 et b7
        (15.5, 7.1) -- (16.5, 7.1)
        ;
      
      % dim 4 & 5
      \draw % enfant de c2
        (c2)
          -- (8, 2)
          -- (15.8, 2)
          -- ($(ast.north)-(xeps)$)
        ;
      \draw % enfant de c0
        (c0)
          -- (ast)
        ;
      \draw % enfant de c1
        (c1)
          -- (24, 2)
          -- (16.2, 2)
          -- ($(ast.north)+(xeps)$)
        ;
      \draw % \ast target de c0, c1 et c7
        (15.5, 1.1) -- (16.5, 1.1)
        ;

      \pgfsetcornersarced{\pgfpointorigin}

    \end{tikzpicture}

  }
% 
% \end{document}
    \captionof{figure}{A $3$-dimensional DFC}
    \label{fig:opetope_loops_poset}
  \end{figure}
  Where we have the local orders $b_5 \triangleleft b_4$ (as chains $c_1 \prec^{\o} \bullet \prec^{-} a_2$) and $b_6 \triangleleft b_3$ (as chains $c_1 \prec^{\o} \bullet \prec^{-} a_1$).
\end{exm}

\para Non-positive opetopes may have loops (cells whose only source is also their target), which make them behave more complicatedly than \emph{positive} opetopes (that is, opetopes where sources and target are always disjoint).
  When reduced to this simpler case, the formalism of MOP boils down to \emph{positive-to-one posets}, their positive counterpart defined in \cite{leclerc2023}.
  One of the complications of the non-positive case is the following: there is no canonical order for composing two loops $x,\,x' \succ^{\o} y$ on a same cell $y$.
  For instance, one may consider \refexm{exm:opetope_loops} of the non-positive opetope $\varrho$.
  There is no canonical order for comparing $b_4,\, b_5$ and $b_6$. Hence we need to add such a choice of order to the structure of many-to-one posets, so we know how we should compose the cells of the opetope.
  This is the order $\triangleleft$ of \refdefin{defin:MOP}.

\begin{defin}[Dendritic face complex]\label{defin:DFC}
  A dendritic face complex (DFC) is a many-to-one poset $C$ satisfying the following extra axioms:
  \begin{itemize}
    \item (\emph{greatest element})\\
      There is a greatest element of positive dimension in $C$, for the partial order induced by $\prec$.
    \item (\emph{loops})\\
      When $z\prec^{\o}y$ in $C$, there is some $x \in C$ with $y \prec^{+} x$.
    \item (\emph{oriented thinness})
      \begin{itemize}
        \item 
          For $z \prec^{\beta} y \prec^{\alpha} x$ with $\alpha,\,\beta \neq \o$ in $P$, there is a unique $y'\neq y$ in $P$ such that $z \prec y' \prec x$.\\
          Hence there is a lozenge as in \ref{fig:dfc_def_lozenge} below.\\ Moreover, we have the \emph{sign rule} $\alpha \beta = - \alpha' \beta'$.
        \item 
          For $z \prec^{\o} y \prec^{-} x$, there is a $y' (\neq y)$ such that $z \prec^{\beta} y' \prec^{\alpha} x$ with $(\alpha,\,\beta) \neq (-,\,\o)$.
      \end{itemize}
      In both cases, when finding such a $y'$ we say that we \emph{complete the half lozenge} (or \emph{chain}) $z \prec y \prec x$.
    \item (\emph{acyclicity})\\
      For $x\in P_{\ge1}$, there is no cycle as in \ref{fig:dfc_def_cycle} below.
  \end{itemize}
\end{defin}
\vspace{-5pt}
\begin{figure}[H]
  \begin{minipage}{.5\textwidth}
    % https://q.uiver.app/#q=WzAsNCxbMSwyLCJ6Il0sWzAsMSwieSJdLFsyLDEsInknIl0sWzEsMCwieCJdLFszLDEsIlxcYWxwaGEiLDEseyJzdHlsZSI6eyJoZWFkIjp7Im5hbWUiOiJub25lIn19fV0sWzEsMCwiXFxiZXRhIiwxLHsic3R5bGUiOnsiaGVhZCI6eyJuYW1lIjoibm9uZSJ9fX1dLFszLDIsIlxcYWxwaGEnIiwxLHsic3R5bGUiOnsiaGVhZCI6eyJuYW1lIjoibm9uZSJ9fX1dLFsyLDAsIlxcYmV0YSciLDEseyJzdHlsZSI6eyJoZWFkIjp7Im5hbWUiOiJub25lIn19fV1d
    \[\begin{tikzcd}[row sep = 15pt]
      & x \\
      y && {y'} \\
      & z
      \arrow["\alpha"{description}, no head, from=1-2, to=2-1]
      \arrow["\beta"{description}, no head, from=2-1, to=3-2]
      \arrow["{\alpha'}"{description}, no head, from=1-2, to=2-3]
      \arrow["{\beta'}"{description}, no head, from=2-3, to=3-2]
    \end{tikzcd}\]
    \captionof{figure}{Lozenge}
    \label{fig:dfc_def_lozenge}
  \end{minipage}
  \begin{minipage}{.5\textwidth}
    % https://q.uiver.app/#q=WzAsMTAsWzIsMCwieCJdLFswLDEsInkgPSB5XzEiXSxbMCwyLCJcXGdhbW1hKHlfMikiXSxbMSwxLCJ5XzIiXSxbMSwyLCJcXGdhbW1hKHlfMykiXSxbMiwxLCJ5XzMiXSxbMiwyLCJcXGdhbW1hKHlfNCkiXSxbMywxLCJ5X3AiXSxbMywyLCJcXGdhbW1hKHlfMSkiXSxbNCwxLCJ5XzEiXSxbMSwyLCIiLDAseyJzdHlsZSI6eyJoZWFkIjp7Im5hbWUiOiJub25lIn19fV0sWzMsMiwiIiwyLHsic3R5bGUiOnsiYm9keSI6eyJuYW1lIjoiYmFycmVkIn0sImhlYWQiOnsibmFtZSI6Im5vbmUifX19XSxbMyw0LCIiLDAseyJzdHlsZSI6eyJoZWFkIjp7Im5hbWUiOiJub25lIn19fV0sWzQsNSwiIiwwLHsic3R5bGUiOnsiYm9keSI6eyJuYW1lIjoiYmFycmVkIn0sImhlYWQiOnsibmFtZSI6Im5vbmUifX19XSxbNSw2LCIiLDAseyJzdHlsZSI6eyJoZWFkIjp7Im5hbWUiOiJub25lIn19fV0sWzcsNiwiIiwyLHsic3R5bGUiOnsiYm9keSI6eyJuYW1lIjoiZGFzaGVkIn0sImhlYWQiOnsibmFtZSI6Im5vbmUifX19XSxbNyw4LCIiLDAseyJzdHlsZSI6eyJoZWFkIjp7Im5hbWUiOiJub25lIn19fV0sWzgsOSwiIiwwLHsic3R5bGUiOnsiYm9keSI6eyJuYW1lIjoiYmFycmVkIn0sImhlYWQiOnsibmFtZSI6Im5vbmUifX19XSxbMCwxLCIiLDAseyJjdXJ2ZSI6MSwic3R5bGUiOnsiaGVhZCI6eyJuYW1lIjoibm9uZSJ9fX1dLFswLDMsIiIsMCx7InN0eWxlIjp7ImhlYWQiOnsibmFtZSI6Im5vbmUifX19XSxbMCw1LCIiLDAseyJzdHlsZSI6eyJoZWFkIjp7Im5hbWUiOiJub25lIn19fV0sWzAsNywiIiwwLHsic3R5bGUiOnsiaGVhZCI6eyJuYW1lIjoibm9uZSJ9fX1dLFswLDksIiIsMix7ImN1cnZlIjotMSwic3R5bGUiOnsiaGVhZCI6eyJuYW1lIjoibm9uZSJ9fX1dXQ==
    \[\begin{tikzcd}[row sep = 15pt, column sep = small]
      && x \\
      {y = y_1} & {y_2} & {y_3} & {y_p} & {y_1} \\
      {\gamma(y_2)} & {\gamma(y_3)} & {\gamma(y_4)} & {\gamma(y_1)}
      \arrow[no head, from=2-1, to=3-1]
      \arrow["\shortmid"{marking}, no head, from=2-2, to=3-1]
      \arrow[no head, from=2-2, to=3-2]
      \arrow["\shortmid"{marking}, no head, from=3-2, to=2-3]
      \arrow[no head, from=2-3, to=3-3]
      \arrow["\shortmid"{marking}, dotted, no head, from=2-4, to=3-3]
      \arrow[no head, from=2-4, to=3-4]
      \arrow["\shortmid"{marking}, no head, from=3-4, to=2-5]
      \arrow[curve={height=6pt}, no head, from=1-3, to=2-1]
      \arrow[no head, from=1-3, to=2-2]
      \arrow[no head, from=1-3, to=2-3]
      \arrow[no head, from=1-3, to=2-4]
      \arrow[curve={height=-6pt}, no head, from=1-3, to=2-5]
    \end{tikzcd}\]
    \captionof{figure}{Cycle}
    \label{fig:dfc_def_cycle}
  \end{minipage}
\end{figure}

\begin{exm}\label{exm:opetope_loops_poset}
  The MOP of \refexm{exm:opetope_loops} is a DFC.
\end{exm}

\begin{exm}\label{exm:four_opetope_poset}
  We may encode the opetope $\omega$ drawn in \refexm{exm:four_opetope}.
  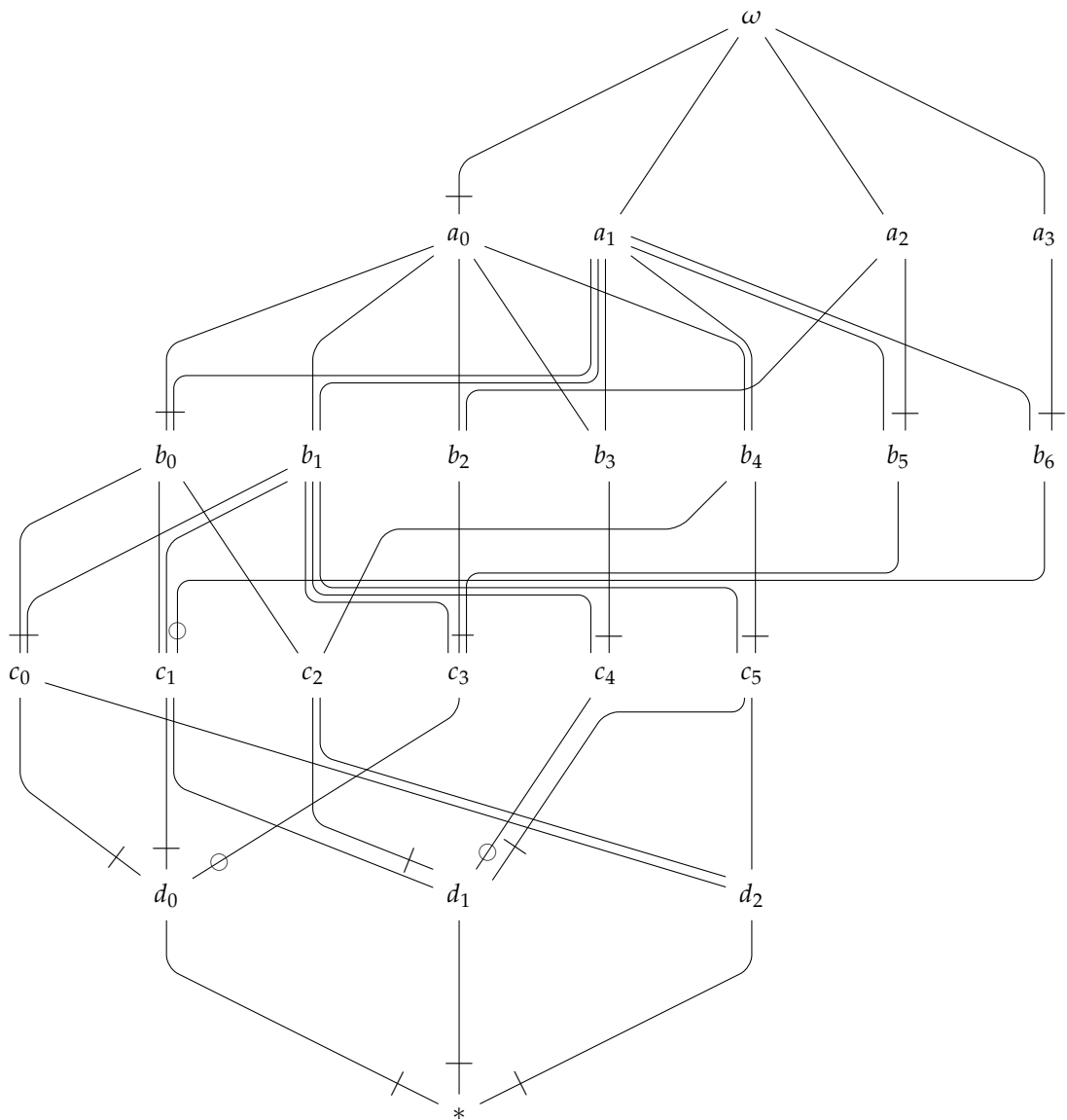
\begin{figure}[H]
    % \documentclass[landscape]{article}
% \usepackage[utf8]{inputenc}
% \usepackage{graphics, graphicx, tikz, geometry}
% \usepackage{amsfonts}
% \usepackage{amssymb,amsmath,amsthm,stmaryrd,mathrsfs,wasysym}
%
% %%% Set the fonts
% \usepackage{mathpazo}
% \usepackage[scaled=0.95]{helvet}
% \usepackage{courier}
% \usepackage{adjustbox}
% \geometry{hmargin=0.5cm,vmargin=1cm}
% \usetikzlibrary{arrows, arrows.meta}
% \linespread{1.05} % Palatino looks better with this
% \definecolor{mypurple}{rgb}{0.7,0,0.7}
% \usetikzlibrary{calc}
%
% \begin{document}
%
  \adjustbox{scale=1, center}{

    \pgfmathsetmacro{\eps}{0.2}
    \begin{tikzpicture}
      [ scale = 0.5
      , every node/.style = {outer sep = 2}]
      \node[draw = none] (xeps) at (\eps, 0) {};

      % --- Noeuds ---
      
      % dim 4
      \node (omega) at (20, 30) {$\omega$};

      % dim 3
      \node (a0) at (12, 24) {$a_0$};
      \node (a1) at (16, 24) {$a_1$};
      \node (a2) at (24, 24) {$a_2$};
      \node (a3) at (28, 24) {$a_3$};

      % dim 2
      \node (b0) at (4, 18) {$b_0$};
      \node (b1) at (8, 18) {$b_1$};
      \node (b2) at (12, 18) {$b_2$};
      \node (b3) at (16, 18) {$b_3$};
      \node (b4) at (20, 18) {$b_4$};
      \node (b5) at (24, 18) {$b_5$};
      \node (b6) at (28, 18) {$b_6$};

      % dim 1
      \node (c0) at (0, 12) {$c_0$};
      \node (c1) at (4, 12) {$c_1$};
      \node (c2) at (8, 12) {$c_2$};
      \node (c3) at (12, 12) {$c_3$};
      \node (c4) at (16, 12) {$c_4$};
      \node (c5) at (20, 12) {$c_5$};

      % dim 0
      \node (d0) at (4, 6) {$d_0$};
      \node (d1) at (12, 6) {$d_1$};
      \node (d2) at (20, 6) {$d_2$};

      % dim -1
      \node (ast) at (12, 0) {$\ast$};

      % --- relations ---
      \pgfsetcornersarced{\pgfpoint{5}{5}}

      % dim 0 & 1
      \draw 
        (omega)
          -- (12, 26)
          -- node[sloped, pos = 0.65] {$\mid$} (a0)
        (omega)
          -- (a1)
        (omega)
          -- (a2)
        (omega)
          -- (28, 26)
          -- (a3)
        ;

      % dim 1 & 2
      \draw % enfants de a0
        (a0)
          -- (4, 21)
          -- (b0)
        (a0)
          -- (8, 21)
          -- (b1)
        (a0)
          -- (b2)
        (a0)
          -- (b3)
        (a0)
          -- (20 - \eps,21)
          -- ($(b4.north) - (xeps)$)
        ;
      \draw % enfants de a1
        ($(a1.south) - 2*(xeps)$)
          -- (16 - 2*\eps, 20 + \eps)
          -- (4 + \eps, 20 + \eps)
          -- ($(b0.north) + (xeps)$)
        ($(a1.south) - (xeps)$)
          -- (16 - \eps, 20)
          -- (8 + \eps, 20)
          -- ($(b1.north) + (xeps)$)
        (a1)
          -- (b3)
        (a1)
          -- (20, 21)
          -- (b4)
        (a1)
          -- (24 - 2*\eps, 21)
          -- ($(b5.north) - 2*(xeps)$)
        (a1.east)
          -- (28 - 2*\eps, 19.5 + \eps)
          -- ($(b6.north) - 2*(xeps)$)
        ;
      \draw % enfants de a2
        (a2)
          -- (20, 20 - \eps)
          -- (12 + \eps, 20 - \eps)
          -- ($(b2.north)+(xeps)$)
        ($(a2.south) + (xeps)$)
          -- node[sloped, pos = 0.9] {$\mid$} ($(b5.north) + (xeps)$)
        ;
      \draw % enfants de a3
        ($(a3.south) + (xeps)$)
          -- node[sloped, pos = 0.9] {$\mid$} ($(b6.north) + (xeps)$)
        ;
      \draw % b0 target de a0 et a1
        (3.7, 19.2) -- (4.5, 19.2)
        ;

      % dim 2 & 3
      \draw % enfants de b0
        (b0)
          -- (0, 16)
          -- (c0)
        ($(b0.south) - (xeps)$)
          -- ($(c1.north) - (xeps)$)
        (b0)
          -- (c2)
        ;
      \draw % enfants de b1
        (b1)
          -- (\eps, 14)
          -- ($(c0.north) + (xeps)$)
        ($(b1.south west)$)
          -- (4, 16 - 2*\eps)
          -- (c1)
        ($(b1.south) - (xeps)$)
          -- (8 - \eps, 14)
          -- (12 - 1.5*\eps, 14)
          -- ($(c3.north)-1.5*(xeps)$)
        (b1)
          -- (8, 14 + \eps)
          -- (16 - 2*\eps, 14 + \eps)
          -- ($(c4.north)-2*(xeps)$)
        ($(b1.south) + (xeps)$)
          -- (8 + \eps, 14 + 2*\eps)
          -- (20 - 2*\eps, 14 + 2*\eps)
          -- ($(c5.north)-2*(xeps)$)
        ;
      \draw % enfants de b2
        (b2)
          -- (c3)
        ;
      \draw % enfants de b3
        ($(b3.south)+.5*(xeps)$)
          -- node[sloped, pos = 0.9] {$\mid$} ($(c4.north)+.5*(xeps)$)
        ;
      \draw % enfants de b4
        (b4)
          -- (18, 16)
          -- (10, 16)
          -- (c2)
        ($(b4.south)+.5*(xeps)$)
          -- node[sloped, pos = 0.9] {$\mid$} ($(c5.north)+.5*(xeps)$)
        ;
      \draw % enfants de b5
        (b5)
          -- (24, 14 + 4*\eps)
          -- (12 + \eps, 14 + 4*\eps)
          -- ($(c3.north)+(xeps)$)
        ;
      \draw % enfants de b6
        (b6)
          -- (28, 14 + 3*\eps)
          -- (4 + 1.5*\eps, 14 + 3*\eps)
          -- node[sloped, pos = 0.7] {$\Circle$} ($(c1.north)+1.5*(xeps)$)
        ;
      \draw % c0 target de b0 et b1
        (-0.3, 13.1) -- (0.5, 13.1)
        ;
      \draw % c3 target de b2 et b5
        (11.8, 13.1) -- (12.4, 13.1)
        ;

      % dim 3 & 4
      \draw % enfants de c0
        (c0)
          -- (0, 9)
          -- node[sloped, pos = 0.8] {$\mid$} (d0)
        (c0)
          -- (d2)
        ;
      \draw % enfants de c1
        (c1)
          -- node[sloped, pos = 0.88] {$\mid$} (d0)
        ($(c1.south)+(xeps)$)
          -- (4 + \eps, 9)
          -- ($(d1.west)+(0,\eps)$)
        ;
      \draw % enfants de c2
        (c2)
          -- (8, 8)
          -- node[sloped, pos = 0.8] {$\mid$} (d1.north west)
        ($(c2.south)+(xeps)$)
          -- (8 + \eps, 9.8)
          -- ($(d2.north west)-(0,\eps)$)
        ;
      \draw % enfants de c3
        (c3)
          -- (12, 11)
          -- node[sloped, pos = 0.9] {$\Circle$} (d0)
        ;
      \draw % enfants de c4
        (c4)
          -- node[sloped, pos = 0.9] {$\Circle$} (d1)
        ;
      \draw % enfants de c5
        ($(c5.south) - (xeps)$)
          -- (20 - \eps, 11)
          -- (16, 11)
          -- node[sloped, pos = 0.8] {$\mid$} ($(d1.east)+(\eps,2*\eps)$)
        (c5)
          -- (d2)
        ;
      
      % dim 4 & 5
      \draw % enfant de d0
        (d0)
          -- (4, 4)
          -- node [sloped, pos = 0.85] {$\mid$} (ast)
        ;
      \draw % enfant de d1
        (d1)
          -- node [sloped, pos = 0.82] {$\mid$} (ast)
        ;
      \draw % enfant de d2
        (d2)
          -- (20, 4)
          -- node [sloped, pos = 0.85] {$\mid$} (ast)
        ;

      \pgfsetcornersarced{\pgfpointorigin}

    \end{tikzpicture}

  }
%
% \end{document}
    \captionof{figure}{A $4$-dimensional DFC}
    \label{fig:four_opetope_poset}
  \end{figure}
\end{exm}

\begin{defin}[Morphism of many-to-one posets]\label{defin:DFC_hom}
  Let $C$ and $D$ be two many-to-one posets.
  A morphism $f : C \to D$ of MOP is a map $f$ between the underlying sets of $C$ and $D$, such that:
  \begin{itemize}
    \item $f$ preserves the dimension.
    \item For all $x\in P_{\ge0}$, $f_x := f\vert_{\delta(x)} : \delta(x) \to \delta(f(x))$ is a bijection and $f(\gamma(x)) = \gamma(f(x))$.
    \item $f$ preserves the local orders.
  \end{itemize}
  We let $\MOP$ be the category of many-to-one posets and their morphisms.
\end{defin}

\begin{defin}[\DFC]\label{defin:dfc_cat}
  The category $\DFC$ is the full subcategory of $\MOP$ whose objects are the DFC.
\end{defin}

\section{Zoom complexes and opetopes}
\label{sec:zoom}

In this section, we present the notion of opetopes as introduced by \textsc{Kock}, \textsc{Joyal}, \textsc{Batanin} and \textsc{Mascari} in \cite{Kock2010}.
We give a slightly different formalisation of their concept by using the notion of rooted tree introduced in \refdefin{defin:rooted_tree}.

\begin{defin}[Constellation]\label{defin:constellation}
  A \emph{constellation} $C : T \to U$ between two rooted trees $T$ and $U$ consists of:
  \begin{itemize}
    \item a subdivision $T'$ of $T$,
    \item a bijection $\sigma_\bullet : \msfbdots(T') \to \msfleaves(U)$,
    \item a bijection $\sigma_\circ : \msfwdots(T') \to \msfndots(U)$,
  \end{itemize}
  such that the sum map $\sigma = \sigma_\bullet + \sigma_\circ$ satisfies the \emph{kernel rule}:
  $$
  \mbox{For each } x \in \msfdots(U), \mbox{the nodes } \{ t \in \msfdots(T') \mid \sigma(t) \sqsubseteq x \} \mbox{ induces a connected full subgraph of } T'.
  $$
  where $x' \sqsubseteq x$ means that there is a descending path from $x'$ to $x$ in $U$.\\
  A constellation will be said to be \emph{exact} when the maps $\sigma_\bullet$ and $\sigma_\circ$ are identities.
\end{defin}

\begin{rmk}\label{rmk:constellation picture}
  Following \cite{Kock2010}, a constellation $C : T \to U$ may be pictured as a nesting of the tree $T'$, such as\\[5px]
  \adjustbox{scale = 0.7, center}
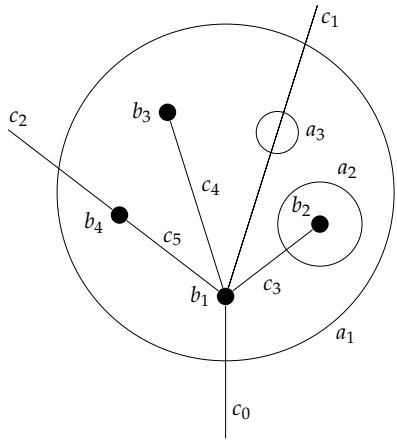
{\begin{tikzpicture}
    [ level distance = 80px
    , bdot/.style = {circle, fill}
    , boundary/.style = {draw = none}
    , scale = 0.8
    ]
    \pgfmathsetmacro{\theta}{35}
    \node (root_3) {} [grow' = up, boundary, level distance = 100px]
      child {
        node (b1) [bdot, label=left:{\large$b_1$}] {} [clockwise from = 90+1.5*\theta, sibling angle = \theta]
          child[level distance = 90px]{
            node (b4) [bdot, label={left, yshift = -4}:{\large$b_4$}] {} [level distance = 100px]
              child{
                node [boundary] {}
                edge from parent node [pos = 0.9, above = 3] {\large$c_2$}
              }
            edge from parent node [above = 2] {\large$c_5$}
          }
          child[level distance = 130px]{
            node (b3) [bdot, label=left:{\large$b_3$}] {}
            edge from parent node [right = 2, pos = 0.6] {\large$c_4$}
          }
          child[level distance = 210px]{
            node [boundary] {}
            edge from parent node [pos = 0.55] (a3) {}
            edge from parent node [pos = 0.95, right = 1] {\large$c_1$}
          }
          child[level distance = 80px]{
            node (b2) [bdot] {}
            edge from parent node [below right, pos = 0.3] {\large$c_3$}
          }
        edge from parent node [right, pos = 0.2] {\large$c_0$}
      };
    \draw ($1.7*(b1)$) circle (4)
      node [right = 65, below = 70] {\large$a_1$};
    \draw (a3) circle (1/2)
      node [right = 12] {\large$a_3$};
    \draw (b2) circle (1)
      node [above = 29, right = 6] {\large$a_2$};
    \node at (b2) [xshift = -10, yshift = 10] {\large$b_2$};
  \end{tikzpicture}}
  where the $c_i$'s are the edges of $T$, the $b_i$'s are the nodes of $T$ and the $a_i$'s are the nodes of $U$.
\end{rmk}

\begin{defin}[Zoom complex]\label{defin:zoom_complex}
  A \emph{zoom complex} of \emph{degree} (or \emph{dimension}) $n$ is a sequence
  $$
  T_0 \overset{C_1}{\longrightarrow} T_1 \overset{C_2}{\longrightarrow} \cdots \overset{C_{n}}{\longrightarrow} T_n
  $$
  of rooted trees and constellations. It will be said to be \emph{exact} when all the $C_i$'s are.
\end{defin}

\begin{defin}[Opetope]\label{defin:opetope}
  For $n\in\N$, a $n$-\emph{opetope} is an exact zoom complex $T_0 \to \cdots \to T_n$ of degree $n$ such that:
  \begin{itemize}
    \item $T_0$ is a tree with one root, one leaf and one node.
    \item If $n \ge 1$, $T_1$ has the same form as $T_0$.
    \item If $n \ge 2$, $T_2$ is a linear tree.
  \end{itemize}
\end{defin}

\begin{figure}[H]
  \adjustbox{scale=0.7, center}{

  \begin{tikzpicture}
    [ blackdot/.style = {circle, fill}
    , whitedot/.style = {circle, draw}
    ]

    % constellation 2
    \node (X2) at (0,0) {
      \begin{tikzpicture}
        [ level distance = 80px
        , bdot/.style = {circle, fill}
        , boundary/.style = {draw = none}
        , scale = 0.8
        ]
        \node {} [grow' = up, boundary, level distance = 100px]
          child {
            node (d2) [bdot, label=left:{\large$d_2$\hspace{-1px}~}] {}
              child{
                node {} [boundary]
              }
          };
        \draw (d2) circle (1)
          node [xshift = 25, yshift = -20] {\large$c_2$};
        \draw (d2) circle (2)
          node [xshift = 38, yshift = -40] {\large$c_1$};
      \end{tikzpicture}
    };

    % zoom
    \node at (3,0) [whitedot] {}
      child[level distance = 20px, grow = right]{
        node [blackdot] {}
      };

    % constellation 3
    \node (X3) at (7,0) {
      \begin{tikzpicture}
        [ level distance = 80px
        , bdot/.style = {circle, fill}
        , boundary/.style = {draw = none}
        , scale = 0.8
        ] 
        \node (root_3) {} [grow' = up, boundary, level distance = 100px]
          child {
            node (c1) [bdot, label=left:{\large$c_1$}] {}
              child[level distance = 110px]{
                node (c2) [bdot, label=left:{\large$c_2$}] {} [level distance = 80px]
                  child{
                    node {} [boundary]
                    edge from parent node [right = 1, pos = 0.9] {\large$d_2$}
                  }
              }
            edge from parent node [right = 1, pos = 0.1] {\large$d_0$}
          };
        \node (mid_3) at ($(c1)!0.4!(c2)$) {};
        \draw (mid_3) ellipse (60px and 110px) 
          node [xshift = 58] {\large$b_1$};
        \draw (mid_3) [yshift = 10] circle (1/2)
          node [xshift = 20] {\large$b_3$};
        \node at (mid_3) [yshift = -17, xshift = 10] {\large$d_1$};
        \draw (c2) circle (1)
          node [xshift = 28, yshift = -15] {\large$b_4$};
        \draw ($(root_3)!0.6!(c1)$) circle (1/2)
          node [xshift = 19, yshift = 12] {\large$b_2$};
      \end{tikzpicture}
    };

    % zoom
    \node at (10,0) [whitedot] {}
      child[level distance = 20px, grow = right]{
        node [blackdot] {}
      };

    % constellation 4
    \node (X4) at (15,0) {
      \begin{tikzpicture}
        [ level distance = 80px
        , bdot/.style = {circle, fill}
        , boundary/.style = {draw = none}
        , scale = 0.8
        ]
        \pgfmathsetmacro{\theta}{35}
        \node (root_3) {} [grow' = up, boundary, level distance = 100px]
          child {
            node (b1) [bdot, label=left:{\large$b_1$}] {} [clockwise from = 90+1.5*\theta, sibling angle = \theta]
              child[level distance = 90px]{
                node (b4) [bdot, label={left, yshift = -4}:{\large$b_4$}] {} [level distance = 100px]
                  child{
                    node [boundary] {}
                    edge from parent node [pos = 0.9, above = 3] {\large$c_2$}
                  }
                edge from parent node [above = 2] {\large$c_5$}
              }
              child[level distance = 130px]{
                node (b3) [bdot, label=left:{\large$b_3$}] {}
                edge from parent node [right = 2, pos = 0.6] {\large$c_4$}
              }
              child[level distance = 210px]{
                node [boundary] {}
                edge from parent node [pos = 0.55] (a3) {}
                edge from parent node [pos = 0.95, right = 1] {\large$c_1$}
              }
              child[level distance = 80px]{
                node (b2) [bdot] {}
                edge from parent node [below right, pos = 0.3] {\large$c_3$}
              }
            edge from parent node [right, pos = 0.2] {\large$c_0$}
          };
        \draw ($1.7*(b1)$) circle (4)
          node [right = 65, below = 70] {\large$a_1$};
        \draw (a3) circle (1/2)
          node [right = 12] {\large$a_3$};
        \draw (b2) circle (1)
          node [above = 29, right = 6] {\large$a_2$};
        \node at (b2) [xshift = -10, yshift = 10] {\large$b_2$};
      \end{tikzpicture}
    };

    \end{tikzpicture}

}
  \vspace{-15pt}
  \captionof{figure}{A $4$-opetope}
  \label{fig:opetope_zoom_drawing}
\end{figure}
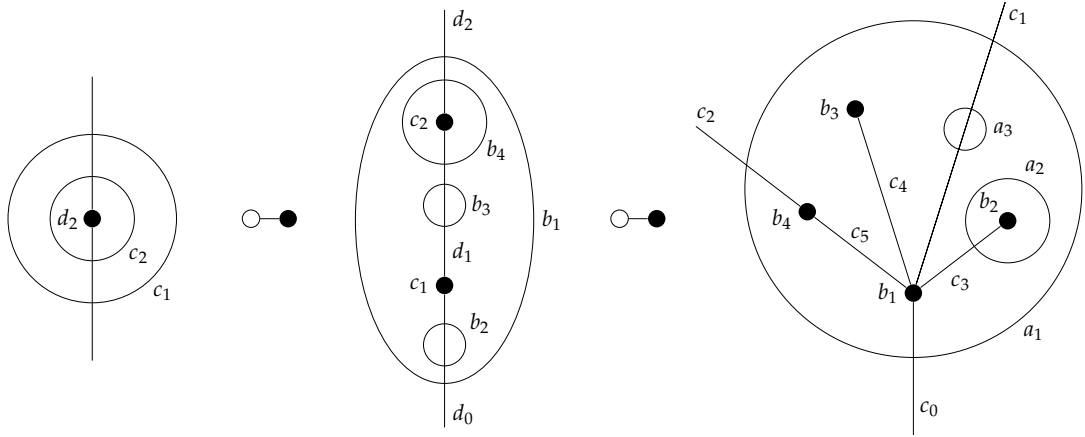

\begin{exm}\label{exm:four_opetope_zoom}
  The opetope associated to \refexm{exm:four_opetope} may be drawn as in \ref{fig:opetope_zoom_drawing} (the first constellation is ommited).
\end{exm}
\vspace{-20pt}

\begin{exm}\label{exm:opetope_zoom_loops}
  The opetope associated to \refexm{exm:four_opetope} may be drawn as in \ref{fig:opetope_zoom_drawing}.
  Notice that in the rightmost constellations, the smallest circles (the minimal ones for the inclusion) really are the whitedots of \refdefin{defin:rooted_tree}.
\end{exm}
\vspace{-25pt}
\begin{figure}[H]
  % \documentclass[landscape]{article}
% \usepackage[utf8]{inputenc}
% \usepackage{graphics, graphicx, tikz, geometry}
% \usepackage{amsfonts}
% \usepackage{amssymb,amsmath,amsthm,stmaryrd,mathrsfs,wasysym}
% 
% %%% Set the fonts
% \usepackage{mathpazo}
% \usepackage[scaled=0.95]{helvet}
% \usepackage{courier}
% \usepackage{adjustbox}
% \geometry{hmargin=0.5cm,vmargin=1cm}
% \usetikzlibrary{arrows, arrows.meta}
% \linespread{1.05} % Palatino looks better with this
% \definecolor{mypurple}{rgb}{0.7,0,0.7}
% \usetikzlibrary{calc, trees}
% 
% \begin{document}
% 
\adjustbox{scale=0.7, center}{

  \begin{tikzpicture}
    [ blackdot/.style = {circle, fill}
    , whitedot/.style = {circle, draw}
    ]

    % constellation 1
    \node (X1) at (-6,0) {
      \begin{tikzpicture}
        [ level distance = 80px
        , bdot/.style = {circle, fill}
        , boundary/.style = {draw = none}
        , scale = 0.8
        ] 
        \node {} [grow' = up, boundary]
          child {
            node (mid_1) [bdot] {}
              child{
                node {} [boundary]
              }
          };
        \draw (mid_1) circle (1)
          node [xshift = 26, yshift = -20] {\large$c_2$};
      \end{tikzpicture}
    };

    % zoom
    \node at (-3.8,0) [whitedot] {}
      child[level distance = 20px, grow = right]{
        node [blackdot] {}
      };

    % constellation 2
    \node (X2) at (0,0) {
      \begin{tikzpicture}
        [ level distance = 80px
        , bdot/.style = {circle, fill}
        , boundary/.style = {draw = none}
        , scale = 0.8
        ]
        \node {} [grow' = up, boundary, level distance = 100px]
          child {
            node (c2) [bdot, label=left:{\large$c_2$\hspace{-1px}~}] {}
              child{
                node {} [boundary]
              }
          };
        \draw (c2) circle (1)
          node [xshift = 25, yshift = -20] {\large$b_2$};
        \draw (c2) circle (2)
          node [xshift = 38, yshift = -40] {\large$b_1$};
      \end{tikzpicture}
    };

    % zoom
    \node at (3,0) [whitedot] {}
      child[level distance = 20px, grow = right]{
        node [blackdot] {}
      };

    % constellation 3
    \node (X3) at (7,0) {
      \begin{tikzpicture}
        [ level distance = 80px
        , bdot/.style = {circle, fill}
        , boundary/.style = {draw = none}
        , scale = 0.8
        ] 
        \node (root_3) {} [grow' = up, boundary, level distance = 60px]
          child {
            node (b1) [bdot] {}
              child[level distance = 190px]{
                node (b2) [bdot, label=left:{\large$b_2$}] {} [level distance = 50px]
                  child{
                    node {} [boundary]
                    edge from parent node [right = 1, pos = 0.9] {\large$c_2$}
                  }
              }
            edge from parent node [right = 1, pos = 0.1] {\large$c_0$}
          };
        \draw ($(b1)!0.48!(b2)$) ellipse (60px and 125px) 
          node [xshift = 58] {\large$a_1$};
        \draw ($(b1)!0.725!(b2)$) ellipse (1 and 1.3)
          node [xshift = -31] {\large$a_2$};
        \draw ($(b1)!0.8!(b2)$) circle (0.4) 
          node [xshift = 20, fill=white] {\large$a_3$};
        \draw ($(b1)!0.65!(b2)$) circle (0.4) 
          node [xshift = 20, fill=white] {\large$a_4$};
        \draw ($(b1)!0.41!(b2)$) circle (0.4) 
          node [xshift = 20, fill=white] {\large$a_5$};
        \draw ($(b1)!0.1!(b2)$) ellipse (1 and 1.3)
          node [xshift = -29, yshift = 18] {\large$a_6$};
        \draw ($(b1)!0.15!(b2)$) circle (0.4) 
          node [xshift = 20, fill=white] {\large$a_7$};
        \node at (b1) [xshift = -15, yshift = 5, fill=white] {\large$b_1$};
      \end{tikzpicture}
    };

    \end{tikzpicture}

}
% 
% \end{document}
  \vspace{-30pt}
  \captionof{figure}{A $3$-opetope}
  \label{fig:opetope_zoom_drawing_loops}
\end{figure}
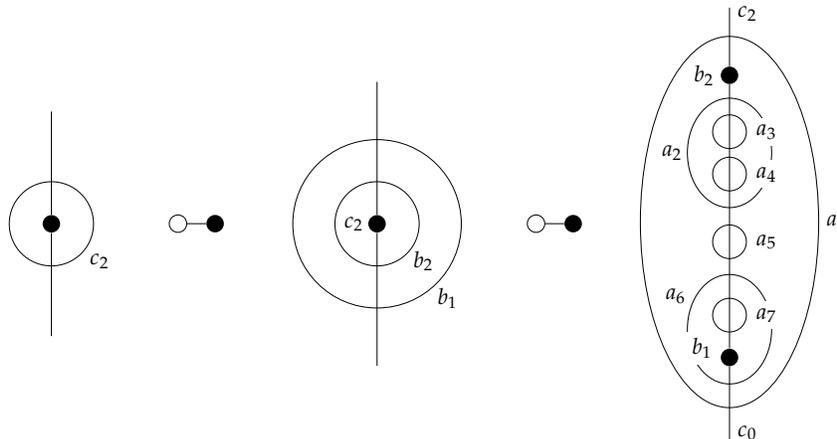

\begin{defin}[Isomorphism of opetopes]
  An isomorphism between two opetopes $(T_i)$ and $(U_i)$ of the same dimsension is a sequence of (subdivided) tree isomorphisms compatible with the structural maps of the constellations.
\end{defin}

\begin{defin}[\KJBM]\label{defin:op_cat}
  We let $\KJBM$ be the category of opetopes and their isomorphisms.
\end{defin}

\section{From dendritic face complexes to opetopes}
\label{sec:DFC_to_Zoom}

From now on, $C$ is a fixed DFC. We let $\omega$ be its greatest element and $n$ its dimension.
Our aim in this section is to build a functor $Z : \DFC^\core \to \KJBM$, where $\DFC^\core$ denotes the \emph{core} subcategory of $\DFC$ whose objects are the DFC and the morphisms are the isomorphisms of DFC.

\begin{prop}[Confinement]
  Let $z \prec^{\o} y \prec^{+} x$ in $C$ and $y' \in \delta(x)$, then $z \prec^{\o} y'$.
  % https://q.uiver.app/#q=WzAsNCxbMSwwLCJ4Il0sWzAsMSwieSJdLFsyLDEsInknIl0sWzEsMiwieiJdLFswLDEsIiIsMCx7InN0eWxlIjp7ImJvZHkiOnsibmFtZSI6ImJhcnJlZCJ9LCJoZWFkIjp7Im5hbWUiOiJub25lIn19fV0sWzEsMywiIiwwLHsic3R5bGUiOnsiaGVhZCI6eyJuYW1lIjoibm9uZSJ9fX1dLFswLDIsIiIsMix7InN0eWxlIjp7ImhlYWQiOnsibmFtZSI6Im5vbmUifX19XSxbMiwzLCIiLDIseyJzdHlsZSI6eyJib2R5Ijp7Im5hbWUiOiJkb3R0ZWQifSwiaGVhZCI6eyJuYW1lIjoibm9uZSJ9fX1dXQ==
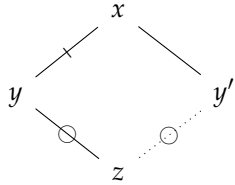
  \[\begin{tikzcd}[row sep = 15pt]
    & x \\
    y && {y'} \\
    & z
    \arrow["\shortmid"{marking}, no head, from=1-2, to=2-1]
    \arrow["\Circle"{marking}, no head, from=2-1, to=3-2]
    \arrow[no head, from=1-2, to=2-3]
    \arrow["\Circle"{marking}, dotted, no head, from=2-3, to=3-2]
  \end{tikzcd}\]
  \label{prop:source_confinement}
\end{prop}
\begin{proof}
  In the configuration $z \prec^{\o} y \prec^{+} x \succ^{-} y'$, let $z' = \gamma(y')$.\\
  We first show that $z' \prec^{\o} y'$. Suppose the opposite, \textit{i.e.} $z' \prec^{+} y'$. Then by completing lozenges from left to right, we may find a pattern as below (the process ends by acyclicity).
  % https://q.uiver.app/#q=WzAsMTAsWzIsMCwieCJdLFsxLDEsIlxcYnVsbGV0Il0sWzAsMSwieSciXSxbMCwyLCJ6Il0sWzEsMiwiXFxidWxsZXQiXSxbMiwxLCJcXGJ1bGxldCJdLFsyLDIsIlxcYnVsbGV0Il0sWzQsMSwieSJdLFszLDEsIlxcYnVsbGV0Il0sWzMsMiwiXFxnYW1tYSh5KSJdLFsyLDAsIiIsMCx7ImN1cnZlIjotMiwic3R5bGUiOnsiaGVhZCI6eyJuYW1lIjoibm9uZSJ9fX1dLFszLDIsIiIsMCx7InN0eWxlIjp7ImJvZHkiOnsibmFtZSI6ImJhcnJlZCJ9LCJoZWFkIjp7Im5hbWUiOiJub25lIn19fV0sWzAsMSwiIiwwLHsic3R5bGUiOnsiaGVhZCI6eyJuYW1lIjoibm9uZSJ9fX1dLFsxLDMsIiIsMCx7InN0eWxlIjp7ImhlYWQiOnsibmFtZSI6Im5vbmUifX19XSxbMSw0LCIiLDAseyJzdHlsZSI6eyJib2R5Ijp7Im5hbWUiOiJiYXJyZWQifSwiaGVhZCI6eyJuYW1lIjoibm9uZSJ9fX1dLFswLDUsIiIsMCx7InN0eWxlIjp7ImhlYWQiOnsibmFtZSI6Im5vbmUifX19XSxbNSw0LCIiLDEseyJzdHlsZSI6eyJoZWFkIjp7Im5hbWUiOiJub25lIn19fV0sWzUsNiwiIiwwLHsic3R5bGUiOnsiaGVhZCI6eyJuYW1lIjoibm9uZSJ9fX1dLFs2LDgsIiIsMCx7InN0eWxlIjp7ImJvZHkiOnsibmFtZSI6ImRvdHRlZCJ9LCJoZWFkIjp7Im5hbWUiOiJub25lIn19fV0sWzAsOCwiIiwwLHsic3R5bGUiOnsiaGVhZCI6eyJuYW1lIjoibm9uZSJ9fX1dLFswLDcsIiIsMCx7ImN1cnZlIjotMiwic3R5bGUiOnsiYm9keSI6eyJuYW1lIjoiYmFycmVkIn0sImhlYWQiOnsibmFtZSI6Im5vbmUifX19XSxbOCw5LCIiLDAseyJzdHlsZSI6eyJib2R5Ijp7Im5hbWUiOiJiYXJyZWQifSwiaGVhZCI6eyJuYW1lIjoibm9uZSJ9fX1dLFs3LDksIiIsMSx7InN0eWxlIjp7ImJvZHkiOnsibmFtZSI6ImJhcnJlZCJ9LCJoZWFkIjp7Im5hbWUiOiJub25lIn19fV1d
  
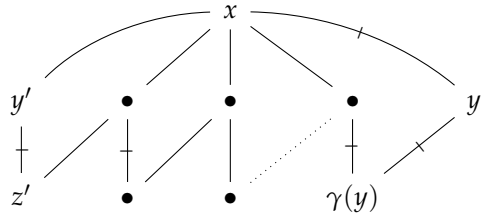
\[\begin{tikzcd}
    && x \\
    {y'} & \bullet & \bullet & \bullet & y \\
    {z'} & \bullet & \bullet & {\gamma(y)}
    \arrow[curve={height=-12pt}, no head, from=2-1, to=1-3]
    \arrow["\shortmid"{marking}, no head, from=3-1, to=2-1]
    \arrow[no head, from=1-3, to=2-2]
    \arrow[no head, from=2-2, to=3-1]
    \arrow["\shortmid"{marking}, no head, from=2-2, to=3-2]
    \arrow[no head, from=1-3, to=2-3]
    \arrow[no head, from=2-3, to=3-2]
    \arrow[no head, from=2-3, to=3-3]
    \arrow[dotted, no head, from=3-3, to=2-4]
    \arrow[no head, from=1-3, to=2-4]
    \arrow["\shortmid"{marking}, curve={height=-12pt}, no head, from=1-3, to=2-5]
    \arrow["\shortmid"{marking}, no head, from=2-4, to=3-4]
    \arrow["\shortmid"{marking}, no head, from=2-5, to=3-4]
  \end{tikzcd}\]
  Hence, $z = \gamma(y) \prec^{+} y$, but we assumed $z \prec^{\o} y$. Whence $z' \prec^{\o} y'$.\\
  We then prove $z = z'$. Using oriented thinness, we may complete $z' \prec^{\o} y' \prec^{-} x$ as $z' \prec^{\beta} y'' \prec^{\alpha} x$ with $(\alpha,\,\beta) \neq (-,\,\o)$. % $\alpha \beta \neq 0$ or $(\alpha,\,\beta) = (+,\,\o)$.
  \begin{itemize}
    \item Suppose $(\alpha,\,\beta) = (-,\,-)$, then oriented thinness applied to the chain $(z' \prec y'' \prec x)$ yields a lozenge completion $(z' \prec^{-} y''' \prec^{+} x)$ or $(z' \prec^{-} y''' \prec^{+} x)$.
      In the first case, we have $y''' = \gamma(x) = y$, but $y$ is a loop, hence $\delta(y) = \emptyset$ and $z' \prec^{-} y$ is impossible. In the second case we are reduced to the following case (replacing $y''$ by $y'''$).
    \item Suppose $(\alpha,\,\beta) = (-,\,+)$, then we have $z \prec^{\o} y \prec^{+} x \succ^{-} y''$, with $z' = \gamma(y'')$.
      Using the same reasoning as at the beginning of the proof, we see that $z' \prec^{\o} y'''$, which contradicts $z' \prec^{+} y'''$.
    \item Suppose $\alpha = \o$, then $y'' = y = \gamma(x)$, but $y \prec^{+} x$ and $y'' \prec^{\o} x$. This is impossible.
  \end{itemize}
  Hence necessarily $\alpha = +$ and $y'' = y$. Since $z' \in \partial(y'') = \partial(y) = \{z\}$, we have $z' = z$.
\end{proof}

\begin{prop}[Lozenge completion for chains $z \prec^{\o} y \prec^{-} x$]\label{prop:loop_lozenge}
  For $z \prec^{\o} y \prec^{-} x$ in $C$, exactly one of the two properties below holds.
  \begin{itemize}
    \item[(i)] There are two distinct completions $z \prec^{\beta} y' \prec^{\alpha} x$ with $\alpha\beta \in \{-,\,+\}$.
    \item[(ii)] There is a completion $z \prec^{\o} y' \prec^{+} x$.
  \end{itemize}
\end{prop}
\begin{proof}
  That either (i) or (ii) hold is a consequence of oriented thinness.\\
  Suppose (i) and (ii). Then choosing by (i) a completion $z \prec^{\beta} y' \prec^{\alpha} x$ with $\alpha = -,\, \beta \neq \o$ and by (ii) a completion $z \prec^{\o} y'' \prec^{-} x$ yields a lozenge
  % https://q.uiver.app/#q=WzAsNCxbMSwwLCJ4Il0sWzAsMSwieSciXSxbMiwxLCJ5JyciXSxbMSwyLCJ6Il0sWzAsMSwiIiwwLHsic3R5bGUiOnsiaGVhZCI6eyJuYW1lIjoibm9uZSJ9fX1dLFsxLDMsIlxcYmV0YSIsMSx7InN0eWxlIjp7ImhlYWQiOnsibmFtZSI6Im5vbmUifX19XSxbMCwyLCIiLDIseyJzdHlsZSI6eyJib2R5Ijp7Im5hbWUiOiJiYXJyZWQifSwiaGVhZCI6eyJuYW1lIjoibm9uZSJ9fX1dLFsyLDMsIiIsMix7InN0eWxlIjp7ImhlYWQiOnsibmFtZSI6Im5vbmUifX19XV0=
  \[\begin{tikzcd}
    & x \\
    {y'} && {y''} \\
    & z
    \arrow[no head, from=1-2, to=2-1]
    \arrow["\beta"{description}, no head, from=2-1, to=3-2]
    \arrow["\shortmid"{marking}, no head, from=1-2, to=2-3]
    \arrow["\Circle"{marking}, no head, from=2-3, to=3-2]
  \end{tikzcd}\]
  Contradicting \refprop{prop:source_confinement}.
\end{proof}

\begin{defin}[$\Gamma,\,\Lambda,\,\Omega,\,N$]
  \vspace{-12pt}
  \begin{mathpar}
    \mbox{We let} \and
    \Gamma := \delta^{+}(C) 
    \and \Lambda := C \setminus \Gamma 
    \and \Omega := \{x \in C \mid x \mbox{ is a loop}\}
    \and N := \{x \in C \mid \delta(x) = \emptyset\}
    \and \Lambda_k := \Lambda \cap C_k,\,\Gamma_k := \Gamma \cap C_k,\, \Omega_k := \Omega \cap C_k,\, N_k := N \cap C_k \quad(k\in\N).
  \end{mathpar}
\end{defin}

\begin{defin}[$\triangleleft^{C_k,\,+},\,\triangleleft^{C_k,\,-},\,<^{C_k,\,+},\,<^{C_k,\,-}$]\label{defin:path_order}
  We define the following relations:
  \begin{itemize}
    \item $<^{C_0,\,-}$ is the empty relation.\\
      For $k > 0$, $<^{C_k,\,-}$ is the transitive closure of
      the relation $\triangleleft^{C_k,\,-}$ on $C_k$, such that $x \triangleleft^{C_k,-} x'$ iff $x \succ^{+} \gamma(x) \prec^{-} x'$. 
      We write $x \bot^{-} x'$ iff either $x <^{-} x'$ or $x' <^{-} x$, and we write $x \le^{-} x'$ iff either $x = x'$
      or $x <^{-} x'$.\\
      A \emph{lower path} is a sequence $x_0 \succ^{+} y_0 \prec^{-} \cdots \prec^{-} x_{p-1} \succ^{+} y_{p} \prec^{-} x_{p}$.
    \item For $k>0$, $<^{C_k,\,+}$ is the transitive closure of the relation $\triangleleft^{C_k,+}$ on $C_k$, such that
      $x \triangleleft^{C_k,+} x'$ iff there is $w \in S_{k+1}$, such that $x \prec^{-} w \succ^{+} x'$.
      We write $x \bot^{+} x'$ iff either $x <^{+} x'$ or $x' <^{+} x$, and we write $x \le^{+} x'$ iff either $x = x'$
      or $x <^{+} x'$.\\
      An \emph{upper path} is a sequence $y_0 \prec^{-} x_1 \succ^{+} y_1 \prec^{-} \cdots \prec^{-} x_p \succ^{+} y_p$.
  \end{itemize}
\end{defin}

\begin{prop}[strictness]\label{prop:DFC_strictness_1}
  The partial order $<^{+}$ is strict.
\end{prop}
\begin{proof}
  We shall prove the following three properties:
  $$
  \begin{array}{ccl}
      \left(k\in\intInter{0}{n}\right) & \mathcal{P}_k : & <^{\scriptscriptstyle C_k,\,+} \mbox{ is a strict partial order}. \\
      \left(k\in\intInter{0}{n-1}\right) & \mathcal{Q}_k : & \delta^{-}(C_{k+1}) = \delta^{-}(\Lambda_{k+1}) \wedge \gamma^{+}(C_{k+1}) = \gamma^{+}(\Lambda_{k+1}). \\
      \left(k\in\intInter{0}{n-1}\right) & \mathcal{R}_k : & \Lambda_k \subseteq \delta(\gamma^{(k+1)}\omega).
  \end{array}
  $$
  by induction on the codimension $n-k$.
  \begin{itemize}
    \item $\boxed{k = n}$\:
      Since $\omega$ is a greatest element, we have $C_n = \{\omega\}$. Hence $\mathcal{P}_n$ is clear.
    \item $\boxed{k = n-1}$\:
      Because $\omega$ is a greatest element, when $a\trigup{n-1}a'$, the only possible situation is:
      % https://q.uiver.app/?q=WzAsMyxbMSwwLCJcXG9tZWdhIl0sWzAsMSwiYSJdLFsyLDEsImEnIl0sWzAsMSwiIiwyLHsic3R5bGUiOnsiaGVhZCI6eyJuYW1lIjoibm9uZSJ9fX1dLFswLDIsIiIsMCx7InN0eWxlIjp7ImJvZHkiOnsibmFtZSI6ImJhcnJlZCJ9LCJoZWFkIjp7Im5hbWUiOiJub25lIn19fV1d
      
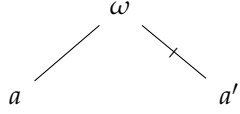
\[\begin{tikzcd}
        & \omega \\
        a && {a'}
        \arrow[no head, from=1-2, to=2-1]
        \arrow["\shortmid"{marking}, no head, from=1-2, to=2-3]
      \end{tikzcd}\]
      Hence a cycle $a_0\trigup{n}a_1\trigup{n}\cdots\trigup{n}a_0$
      must have the form $a_0\trigup{n}a_0$ with $n=0$ and $a_0 = \gamma(\omega)$.
      But $\gamma(\omega)\nprec^{-}\omega$. So there is no such cycle.
      Whence $\mathcal{P}_{n-1}$.

      Because $\omega\in\Lambda_n$ is a greatest element, $\mathcal{Q}_{n-1}$ and $\mathcal{R}_{n-1}$ are clear.

    \item \boxed{\mbox{Induction}}
      We suppose now $\mathcal{P}_{k+1},\,\mathcal{Q}_{k+1},\,\mathcal{R}_{k+1}$, and $\mathcal{R}_{k+2}$ if $k<n-2$.\\
      First, we prove $\mathcal{Q}_k$.
      Let $e \in C_k$ with $e \prec^{\alpha} d$ for some $d =: d_1$ and $\alpha \in \{-,\,+\}$.
      If $d_1$ is in $\Lambda_{k+1}$ then we are done.
      Else, let $c_1$ be such that $d_1 \prec^{+} c_1$. Then we may complete $e\prec^{\alpha}d_1\prec^{+}c_1$ into a lozenge as in \ref{fig:strictness_lozenge} below.
      If $d_2\in\Lambda_{k+1}$, we have finished. Else, we continue with $d_2$, taking $c_2$ such that $d_2\prec^{+}c_2$ \textit{etc.}.
      While iterating this construction, we cannot produce a loop as in \ref{fig:strictness_loop} below:
      \begin{figure}[H]
        \begin{minipage}{.32\textwidth}
          % https://q.uiver.app/#q=WzAsNCxbMCwxLCJkXzEiXSxbMSwyLCJlIl0sWzEsMCwiY18xIl0sWzIsMSwiZF8yIl0sWzAsMSwiXFxhbHBoYSIsMSx7InN0eWxlIjp7ImhlYWQiOnsibmFtZSI6Im5vbmUifX19XSxbMCwyLCIiLDIseyJzdHlsZSI6eyJib2R5Ijp7Im5hbWUiOiJiYXJyZWQifSwiaGVhZCI6eyJuYW1lIjoibm9uZSJ9fX1dLFsyLDMsIiIsMix7InN0eWxlIjp7ImhlYWQiOnsibmFtZSI6Im5vbmUifX19XSxbMywxLCJcXGFscGhhIiwxLHsic3R5bGUiOnsiaGVhZCI6eyJuYW1lIjoibm9uZSJ9fX1dXQ==
          \[\begin{tikzcd}[column sep = tiny]
            & {c_1} \\
            {d_1} && {d_2} \\
            & e
            \arrow["\alpha"{description}, no head, from=2-1, to=3-2]
            \arrow["\shortmid"{marking}, no head, from=2-1, to=1-2]
            \arrow[no head, from=1-2, to=2-3]
            \arrow["\alpha"{description}, no head, from=2-3, to=3-2]
          \end{tikzcd}\]
          \captionof{figure}{Lozenge}
          \label{fig:strictness_lozenge}
        \end{minipage}
        \begin{minipage}{.32\textwidth}
          % https://q.uiver.app/#q=WzAsOSxbMCwxLCJkXzEiXSxbMSwxLCJkXzIiXSxbMiwyLCJlIl0sWzQsMSwiZF8xIl0sWzAsMCwiY18xIl0sWzEsMCwiY18yIl0sWzQsMCwiY18xIl0sWzMsMSwiZF9xIl0sWzMsMCwiY19xIl0sWzQsMCwiIiwyLHsic3R5bGUiOnsiYm9keSI6eyJuYW1lIjoiYmFycmVkIn0sImhlYWQiOnsibmFtZSI6Im5vbmUifX19XSxbNCwxLCIiLDAseyJzdHlsZSI6eyJoZWFkIjp7Im5hbWUiOiJub25lIn19fV0sWzEsNSwiIiwwLHsic3R5bGUiOnsiYm9keSI6eyJuYW1lIjoiYmFycmVkIn0sImhlYWQiOnsibmFtZSI6Im5vbmUifX19XSxbMyw2LCIiLDEseyJzdHlsZSI6eyJib2R5Ijp7Im5hbWUiOiJiYXJyZWQifSwiaGVhZCI6eyJuYW1lIjoibm9uZSJ9fX1dLFs3LDgsIiIsMCx7InN0eWxlIjp7ImJvZHkiOnsibmFtZSI6ImJhcnJlZCJ9LCJoZWFkIjp7Im5hbWUiOiJub25lIn19fV0sWzgsMywiIiwxLHsic3R5bGUiOnsiaGVhZCI6eyJuYW1lIjoibm9uZSJ9fX1dLFswLDIsIlxcYWxwaGEiLDEseyJjdXJ2ZSI6Miwic3R5bGUiOnsiaGVhZCI6eyJuYW1lIjoibm9uZSJ9fX1dLFsxLDIsIlxcYWxwaGEiLDEseyJjdXJ2ZSI6MSwic3R5bGUiOnsiaGVhZCI6eyJuYW1lIjoibm9uZSJ9fX1dLFs3LDIsIlxcYWxwaGEiLDEseyJjdXJ2ZSI6LTEsInN0eWxlIjp7ImhlYWQiOnsibmFtZSI6Im5vbmUifX19XSxbMywyLCJcXGFscGhhIiwxLHsiY3VydmUiOi0yLCJzdHlsZSI6eyJoZWFkIjp7Im5hbWUiOiJub25lIn19fV0sWzUsNywiIiwwLHsic3R5bGUiOnsiYm9keSI6eyJuYW1lIjoiZG90dGVkIn0sImhlYWQiOnsibmFtZSI6Im5vbmUifX19XV0=
          \[\begin{tikzcd}[column sep = small]
            {c_1} & {c_2} && {c_q} & {c_1} \\
            {d_1} & {d_2} && {d_q} & {d_1} \\
            && e
            \arrow["\shortmid"{marking}, no head, from=1-1, to=2-1]
            \arrow[no head, from=1-1, to=2-2]
            \arrow["\shortmid"{marking}, no head, from=2-2, to=1-2]
            \arrow["\shortmid"{marking}, no head, from=2-5, to=1-5]
            \arrow["\shortmid"{marking}, no head, from=2-4, to=1-4]
            \arrow[no head, from=1-4, to=2-5]
            \arrow["\alpha"{description}, curve={height=12pt}, no head, from=2-1, to=3-3]
            \arrow["\alpha"{description}, curve={height=6pt}, no head, from=2-2, to=3-3]
            \arrow["\alpha"{description}, curve={height=-6pt}, no head, from=2-4, to=3-3]
            \arrow["\alpha"{description}, curve={height=-12pt}, no head, from=2-5, to=3-3]
            \arrow[dotted, no head, from=1-2, to=2-4]
          \end{tikzcd}\]
          \captionof{figure}{Impossible loop}
          \label{fig:strictness_loop}
        \end{minipage}
        \begin{minipage}{.32\textwidth}
          % https://q.uiver.app/?q=WzAsOSxbMCwyLCJlXzAiXSxbMSwyLCJlXzEiXSxbMSwxLCJjXzEiXSxbNCwxLCJjX3EiXSxbMiwxLCJjXzIiXSxbMywyLCJcXGNkb3RzIl0sWzIsMiwiZV8yIl0sWzQsMiwiZV9xIl0sWzIsMCwidF57ay0yfVxcb21lZ2EiXSxbMSwyLCIiLDAseyJzdHlsZSI6eyJib2R5Ijp7Im5hbWUiOiJiYXJyZWQifSwiaGVhZCI6eyJuYW1lIjoibm9uZSJ9fX1dLFswLDIsIiIsMCx7InN0eWxlIjp7ImhlYWQiOnsibmFtZSI6Im5vbmUifX19XSxbMSw0LCIiLDIseyJzdHlsZSI6eyJoZWFkIjp7Im5hbWUiOiJub25lIn19fV0sWzYsNCwiIiwwLHsic3R5bGUiOnsiYm9keSI6eyJuYW1lIjoiYmFycmVkIn0sImhlYWQiOnsibmFtZSI6Im5vbmUifX19XSxbNiwzLCIiLDIseyJzdHlsZSI6eyJib2R5Ijp7Im5hbWUiOiJkb3R0ZWQifSwiaGVhZCI6eyJuYW1lIjoibm9uZSJ9fX1dLFszLDcsIiIsMix7InN0eWxlIjp7ImhlYWQiOnsibmFtZSI6Im5vbmUifX19XSxbMiw4LCIiLDAseyJzdHlsZSI6eyJoZWFkIjp7Im5hbWUiOiJub25lIn19fV0sWzgsNCwiIiwwLHsic3R5bGUiOnsiaGVhZCI6eyJuYW1lIjoibm9uZSJ9fX1dLFs4LDMsIiIsMSx7InN0eWxlIjp7ImJvZHkiOnsibmFtZSI6ImRvdHRlZCJ9LCJoZWFkIjp7Im5hbWUiOiJub25lIn19fV1d
          \[\begin{tikzcd}[column sep = small]
            && {\gamma^{(k+2)}(\omega)} \\
            & {d_1} & {d_2} & {d_q} \\
            {e_0} & {e_1} & {e_2} & {e_q}
            \arrow["\shortmid"{marking}, no head, from=3-2, to=2-2]
            \arrow[no head, from=3-1, to=2-2]
            \arrow[no head, from=3-2, to=2-3]
            \arrow["\shortmid"{marking}, no head, from=3-3, to=2-3]
            \arrow[dotted, no head, from=3-3, to=2-4]
            \arrow["\shortmid"{marking}, no head, from=2-4, to=3-4]
            \arrow[no head, from=2-2, to=1-3]
            \arrow[no head, from=1-3, to=2-3]
            \arrow[dotted, no head, from=1-3, to=2-4]
          \end{tikzcd}\]
          \captionof{figure}{Path}
          \label{fig:strictness_path}
        \end{minipage}
      \end{figure}
      because it would contradict $\mathcal{P}_{k+1}$.
      Hence, each time $d_i\notin\Lambda_{k+1}$, the next element $d_{i+1}$ is a new one, and this construction must finish because the poset is finite.
      So when the construction ends, it produces an element $e\prec^{\alpha}d_q\in\Lambda_{k+1}$. Whence the existence.
      
    \vspace{10pt}
    We then prove $\mathcal{R}_k$.
    Let $e\in\Lambda_{k}$. Then $e \in \delta(d)$ for some $d$.\\
    Suppose first that $e \in \delta^{-}(d)$. Then, using $\mathcal{Q}_k$ we may assume $d \in \Lambda_{k+1}$, and by $\mathcal{R}_{k+1}$, $d \in \delta(\gamma^{(k+2)}\omega)$.\\
    Either $d \prec^{\o} \gamma^{(k+2)}\omega$, in which case $d = \gamma^{(k+1)}\omega \succ^{-} e$ and we are done.
    Or $d \prec^{-} \gamma^{(k+2)}\omega$ and, because $e$ is not a source, the only possible lozenge completion for $e \prec^{-} d \prec^{-} \gamma^{(k+2)}\omega$ is $e \prec^{-} \gamma^{(k+1)}\omega \prec^{+} \gamma^{(k+2)}\omega$, and we are done.\\
    Suppose now that $e \prec^{\o} d$. Then there is some $c$ with $d \prec^{+} c$. Using $\mathcal{Q}_{k+1}$, we may assume $c \in \Lambda_{k+2}$.\\
    If $k = n-2$, then $c = \omega$, $d = \gamma(\omega)$ and we are done.
    If $k < n-2$, then using $\mathcal{R}_{k+2}$ we have $c \in \delta(\gamma^{(k+3)}\omega)$.
    Now either $c = \gamma^{(k+2)}\omega \prec^{\o} \gamma^{(k+3)}\omega$, hence $d = \gamma^{(k+1)}\omega$ and we are done, or $c \prec^{-} \gamma^{(k+3)}\omega$ and we are in the situation $e \prec^{\o} d \prec^{+} c \prec^{-} \gamma^{(k+3)}\omega$.
    By successively completing lozenges from left to right, we may construct one of the two diagrams below.
    
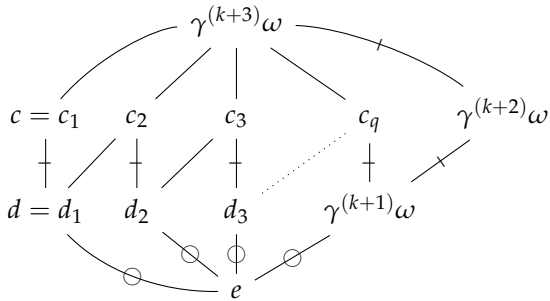
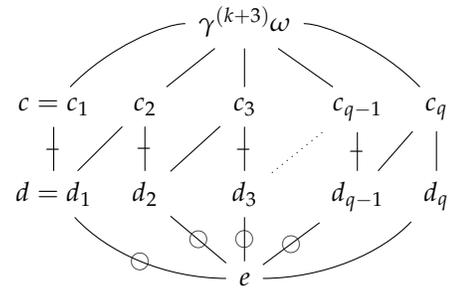
\begin{figure}[H]
      \hspace{10pt}
      \begin{minipage}{.5\textwidth}
        % https://q.uiver.app/#q=WzAsMTEsWzIsMCwidF57KGsrMyl9XFxvbWVnYSJdLFswLDEsImM9Y18xIl0sWzAsMiwiZD1kXzEiXSxbMiwzLCJlIl0sWzEsMSwiY18yIl0sWzEsMiwiZF8yIl0sWzIsMSwiY18zIl0sWzIsMiwiZF8zIl0sWzMsMSwiY19xIl0sWzMsMiwidF57KGsrMSl9XFxvbWVnYSJdLFs0LDEsInReeyhrKzIpfVxcb21lZ2EiXSxbMiwzLCIiLDEseyJjdXJ2ZSI6Miwic3R5bGUiOnsiaGVhZCI6eyJuYW1lIjoibm9uZSJ9fX1dLFswLDEsIiIsMSx7ImN1cnZlIjoyLCJzdHlsZSI6eyJoZWFkIjp7Im5hbWUiOiJub25lIn19fV0sWzEsMiwiIiwxLHsic3R5bGUiOnsiYm9keSI6eyJuYW1lIjoiYmFycmVkIn0sImhlYWQiOnsibmFtZSI6Im5vbmUifX19XSxbMCw0LCIiLDEseyJzdHlsZSI6eyJoZWFkIjp7Im5hbWUiOiJub25lIn19fV0sWzQsMiwiIiwxLHsic3R5bGUiOnsiaGVhZCI6eyJuYW1lIjoibm9uZSJ9fX1dLFs0LDUsIiIsMSx7InN0eWxlIjp7ImJvZHkiOnsibmFtZSI6ImJhcnJlZCJ9LCJoZWFkIjp7Im5hbWUiOiJub25lIn19fV0sWzUsMywiIiwxLHsic3R5bGUiOnsiaGVhZCI6eyJuYW1lIjoibm9uZSJ9fX1dLFswLDYsIiIsMSx7InN0eWxlIjp7ImhlYWQiOnsibmFtZSI6Im5vbmUifX19XSxbNiw1LCIiLDEseyJzdHlsZSI6eyJoZWFkIjp7Im5hbWUiOiJub25lIn19fV0sWzYsNywiIiwxLHsic3R5bGUiOnsiYm9keSI6eyJuYW1lIjoiYmFycmVkIn0sImhlYWQiOnsibmFtZSI6Im5vbmUifX19XSxbNywzLCIiLDEseyJzdHlsZSI6eyJoZWFkIjp7Im5hbWUiOiJub25lIn19fV0sWzAsOCwiIiwxLHsic3R5bGUiOnsiaGVhZCI6eyJuYW1lIjoibm9uZSJ9fX1dLFs4LDcsIiIsMSx7InN0eWxlIjp7ImJvZHkiOnsibmFtZSI6ImRvdHRlZCJ9LCJoZWFkIjp7Im5hbWUiOiJub25lIn19fV0sWzgsOSwiIiwxLHsic3R5bGUiOnsiYm9keSI6eyJuYW1lIjoiYmFycmVkIn0sImhlYWQiOnsibmFtZSI6Im5vbmUifX19XSxbOSwzLCIiLDEseyJzdHlsZSI6eyJoZWFkIjp7Im5hbWUiOiJub25lIn19fV0sWzAsMTAsIiIsMSx7InN0eWxlIjp7ImJvZHkiOnsibmFtZSI6ImJhcnJlZCJ9LCJoZWFkIjp7Im5hbWUiOiJub25lIn19fV0sWzEwLDksIiIsMSx7InN0eWxlIjp7ImJvZHkiOnsibmFtZSI6ImJhcnJlZCJ9LCJoZWFkIjp7Im5hbWUiOiJub25lIn19fV1d
        \[\begin{tikzcd}[row sep = 15pt, column sep = tiny]
          && {\gamma^{(k+3)}\omega} \\
          {c=c_1} & {c_2} & {c_3} & {c_q} & {\gamma^{(k+2)}\omega} \\
          {d=d_1} & {d_2} & {d_3} & {\gamma^{(k+1)}\omega} \\
          && e
          \arrow["\Circle"{marking}, curve={height=12pt}, no head, from=3-1, to=4-3]
          \arrow[curve={height=12pt}, no head, from=1-3, to=2-1]
          \arrow["\shortmid"{marking}, no head, from=2-1, to=3-1]
          \arrow[no head, from=1-3, to=2-2]
          \arrow[no head, from=2-2, to=3-1]
          \arrow["\shortmid"{marking}, no head, from=2-2, to=3-2]
          \arrow["\Circle"{marking}, no head, from=3-2, to=4-3]
          \arrow[no head, from=1-3, to=2-3]
          \arrow[no head, from=2-3, to=3-2]
          \arrow["\shortmid"{marking}, no head, from=2-3, to=3-3]
          \arrow["\Circle"{marking}, no head, from=3-3, to=4-3]
          \arrow[no head, from=1-3, to=2-4]
          \arrow[dotted, no head, from=2-4, to=3-3]
          \arrow["\shortmid"{marking}, no head, from=2-4, to=3-4]
          \arrow["\Circle"{marking}, no head, from=3-4, to=4-3]
          \arrow[curve={height=-12pt}, "\shortmid"{marking}, no head, from=1-3, to=2-5]
          \arrow["\shortmid"{marking}, no head, from=2-5, to=3-4]
        \end{tikzcd}\]
      \end{minipage}
      \begin{minipage}{.5\textwidth}
        \[\begin{tikzcd}[row sep = 15pt, column sep = tiny]
          && {\gamma^{(k+3)}\omega} \\
          {c=c_1} & {c_2} & {c_3} & {c_{q-1}} & {c_q} \\
          {d=d_1} & {d_2} & {d_3} & {d_{q-1}} & {d_q} \\
          && e
          \arrow["\Circle"{marking}, curve={height=12pt}, no head, from=3-1, to=4-3]
          \arrow[curve={height=12pt}, no head, from=1-3, to=2-1]
          \arrow["\shortmid"{marking}, no head, from=2-1, to=3-1]
          \arrow[no head, from=1-3, to=2-2]
          \arrow[no head, from=2-2, to=3-1]
          \arrow["\shortmid"{marking}, no head, from=2-2, to=3-2]
          \arrow["\Circle"{marking}, no head, from=3-2, to=4-3]
          \arrow[no head, from=1-3, to=2-3]
          \arrow[no head, from=2-3, to=3-2]
          \arrow["\shortmid"{marking}, no head, from=2-3, to=3-3]
          \arrow["\Circle"{marking}, no head, from=3-3, to=4-3]
          \arrow[no head, from=1-3, to=2-4]
          \arrow[dotted, no head, from=2-4, to=3-3]
          \arrow["\shortmid"{marking}, no head, from=2-4, to=3-4]
          \arrow["\Circle"{marking}, no head, from=3-4, to=4-3]
          \arrow[curve={height=-12pt}, no head, from=1-3, to=2-5]
          \arrow[no head, from=2-5, to=3-4]
          \arrow[no head, from=2-5, to=3-5]
          \arrow[curve={height=-12pt}, no head, from=3-5, to=4-3]
        \end{tikzcd}\]
      \end{minipage}
    \end{figure}
    In the first case (left one), we are done, and in the second (right one) we are reduced to the case $e \in \delta^{-}(d')$ for some $d'$. This completes the proof of $\mathcal{R}_k$.

    \vspace{10pt}
    Now we prove $\mathcal{P}_k$.
    The strategy is to replace a related pair $e\trigup{k}e'$ by a path\\
    $e = e_0 \trigup{k} e_1 \trigup{k} \cdots \trigup{k} e_q = e'$ as in \ref{fig:strictness_path} above.
    Indeed, if we are able to do this, every cyclic path $e_0\trigup{k}e_1\trigup{k}\cdots\trigup{k}e_q = e_0$ will induce a longer path
    $e_0 = e'_0\trigup{k}e'_1\trigup{k}\cdots\trigup{k}e'_m = e_0$ as above, which cannot exist because it would contradict the acyclicity axiom.
    Notice that either $q=1$ and we are done, or if we know that each $d_i$ is in $\Lambda_{k+1}$ then we know that $d_i \prec^{-} \gamma^{(k+2)}\omega$ (by $\mathcal{R}_{k+1}$, and because $\gamma^{(k+2)}$ has at least two sources, it cannot be a loop).
    So we will be able to conclude if we prove \reflem{lem:DFC_path_in_Lambda} below.
      
    \begin{lem}\label{lem:DFC_path_in_Lambda}
      if $e\trigup{k}e'$, then we have a sequence as follows, with $d_1,\,d_2,\,\cdots,\,d_q\in\Lambda_{k+1}$.
      % https://q.uiver.app/?q=WzAsOCxbMCwxLCJlID0gZV8wIl0sWzEsMSwiZV8xIl0sWzEsMCwiY18xIl0sWzQsMCwiY19xIl0sWzIsMCwiY18yIl0sWzMsMSwiXFxjZG90cyJdLFsyLDEsImVfMiJdLFs0LDEsImVfcSA9IGUnIl0sWzEsMiwiIiwwLHsic3R5bGUiOnsiYm9keSI6eyJuYW1lIjoiYmFycmVkIn0sImhlYWQiOnsibmFtZSI6Im5vbmUifX19XSxbMCwyLCIiLDAseyJzdHlsZSI6eyJoZWFkIjp7Im5hbWUiOiJub25lIn19fV0sWzEsNCwiIiwyLHsic3R5bGUiOnsiaGVhZCI6eyJuYW1lIjoibm9uZSJ9fX1dLFs2LDQsIiIsMCx7InN0eWxlIjp7ImJvZHkiOnsibmFtZSI6ImJhcnJlZCJ9LCJoZWFkIjp7Im5hbWUiOiJub25lIn19fV0sWzYsMywiIiwyLHsic3R5bGUiOnsiYm9keSI6eyJuYW1lIjoiZG90dGVkIn0sImhlYWQiOnsibmFtZSI6Im5vbmUifX19XSxbMyw3LCIiLDIseyJzdHlsZSI6eyJoZWFkIjp7Im5hbWUiOiJub25lIn19fV1d
      
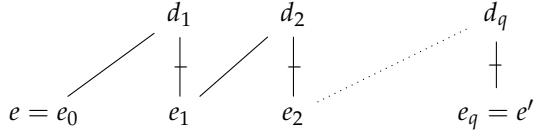
\[\begin{tikzcd}
        & {d_1} & {d_2} && {d_q} \\
        {e = e_0} & {e_1} & {e_2} && {e_q = e'}
        \arrow["\shortmid"{marking}, no head, from=2-2, to=1-2]
        \arrow[no head, from=2-1, to=1-2]
        \arrow[no head, from=2-2, to=1-3]
        \arrow["\shortmid"{marking}, no head, from=2-3, to=1-3]
        \arrow[dotted, no head, from=2-3, to=1-5]
        \arrow["\shortmid"{marking}, no head, from=1-5, to=2-5]
      \end{tikzcd}\]
    \end{lem}
    \textit{Proof.}
      If $e\prec^{-}d\succ^{+}e'$ with $d\in\Lambda_{k+1}$, then the result is proven.
      So we may suppose that $d$ is the target of some element: $d = \gamma(c)$, as in \ref{fig:dfc_to_zpo_hasse_1} below.
      \begin{figure}[H]
        \begin{minipage}{.5\textwidth}
          % https://q.uiver.app/?q=WzAsNCxbMCwyLCJlIl0sWzEsMSwiZCJdLFsyLDIsImUnIl0sWzEsMCwiYyJdLFswLDEsIiIsMCx7InN0eWxlIjp7ImhlYWQiOnsibmFtZSI6Im5vbmUifX19XSxbMSwyLCIiLDAseyJzdHlsZSI6eyJib2R5Ijp7Im5hbWUiOiJiYXJyZWQifSwiaGVhZCI6eyJuYW1lIjoibm9uZSJ9fX1dLFsxLDMsIiIsMCx7InN0eWxlIjp7ImhlYWQiOnsibmFtZSI6Im5vbmUifX19XV0=
          \[\begin{tikzcd}[row sep = small]
            & c \\
            & d \\
            e && {e'}
            \arrow[no head, from=3-1, to=2-2]
            \arrow["\shortmid"{marking}, no head, from=2-2, to=3-3]
            \arrow["\shortmid"{marking}, no head, from=2-2, to=1-2]
          \end{tikzcd}\]
          \captionof{figure}{}
          \label{fig:dfc_to_zpo_hasse_1}
        \end{minipage}
        \begin{minipage}{.5\textwidth}
          % https://q.uiver.app/?q=WzAsNCxbMCwyLCJlIl0sWzEsMSwiZF8xIl0sWzEsMCwiYyJdLFswLDEsImQiXSxbMCwxLCIiLDAseyJzdHlsZSI6eyJoZWFkIjp7Im5hbWUiOiJub25lIn19fV0sWzEsMiwiIiwwLHsic3R5bGUiOnsiaGVhZCI6eyJuYW1lIjoibm9uZSJ9fX1dLFsyLDMsIiIsMCx7InN0eWxlIjp7ImJvZHkiOnsibmFtZSI6ImJhcnJlZCJ9LCJoZWFkIjp7Im5hbWUiOiJub25lIn19fV0sWzMsMCwiIiwwLHsic3R5bGUiOnsiaGVhZCI6eyJuYW1lIjoibm9uZSJ9fX1dXQ==
          \[\begin{tikzcd}[row sep = small]
            & c \\
            d & {d_1} \\
            e
            \arrow[no head, from=3-1, to=2-2]
            \arrow[no head, from=2-2, to=1-2]
            \arrow["\shortmid"{marking}, no head, from=1-2, to=2-1]
            \arrow[no head, from=2-1, to=3-1]
          \end{tikzcd}\]
          \captionof{figure}{}
          \label{fig:dfc_to_zpo_hasse_2}
        \end{minipage}
      \end{figure}
      First, complete the left chain with some $d_1$ as in \ref{fig:dfc_to_zpo_hasse_2} above.
      
      Then we consider its target $\gamma(d_1)$, and complete $\gamma(d_1)\prec^{+}d_1\prec^{-}c$ as a lozenge. We have two possibilities.
      Either the completion is of type \rom{1} as in \ref{fig:dfc_to_zpo_hasse_3} below, or it is of type \rom{2} as in \ref{fig:dfc_to_zpo_hasse_4} below:
      \begin{figure}[H]
        \begin{minipage}{.5\textwidth}
          % https://q.uiver.app/?q=WzAsNixbMCwyLCJlIl0sWzEsMSwiZF8xIl0sWzEsMCwiYyJdLFswLDEsImQiXSxbMSwyLCJ0LmRfMSA9IGUnIl0sWzIsMSwiZF8yID0gZCJdLFswLDEsIiIsMCx7InN0eWxlIjp7ImhlYWQiOnsibmFtZSI6Im5vbmUifX19XSxbMSwyLCIiLDAseyJzdHlsZSI6eyJoZWFkIjp7Im5hbWUiOiJub25lIn19fV0sWzIsMywiIiwwLHsic3R5bGUiOnsiYm9keSI6eyJuYW1lIjoiYmFycmVkIn0sImhlYWQiOnsibmFtZSI6Im5vbmUifX19XSxbMywwLCIiLDAseyJzdHlsZSI6eyJoZWFkIjp7Im5hbWUiOiJub25lIn19fV0sWzEsNCwiIiwwLHsic3R5bGUiOnsiYm9keSI6eyJuYW1lIjoiYmFycmVkIn0sImhlYWQiOnsibmFtZSI6Im5vbmUifX19XSxbMiw1LCIiLDAseyJzdHlsZSI6eyJib2R5Ijp7Im5hbWUiOiJiYXJyZWQifSwiaGVhZCI6eyJuYW1lIjoibm9uZSJ9fX1dLFs1LDQsIiIsMSx7InN0eWxlIjp7ImJvZHkiOnsibmFtZSI6ImJhcnJlZCJ9LCJoZWFkIjp7Im5hbWUiOiJub25lIn19fV1d
          \[\begin{tikzcd}[row sep = small]
            & c \\
            d & {d_1} & {d_2 = d} \\
            e & {\gamma(d_1) = e'}
            \arrow[no head, from=3-1, to=2-2]
            \arrow[no head, from=2-2, to=1-2]
            \arrow["\shortmid"{marking}, no head, from=1-2, to=2-1]
            \arrow[no head, from=2-1, to=3-1]
            \arrow["\shortmid"{marking}, no head, from=2-2, to=3-2]
            \arrow["\shortmid"{marking}, no head, from=1-2, to=2-3]
            \arrow["\shortmid"{marking}, no head, from=2-3, to=3-2]
          \end{tikzcd}\]
          \captionof{figure}{Type \rom{1}}
          \label{fig:dfc_to_zpo_hasse_3}
        \end{minipage}
        \begin{minipage}{.5\textwidth}
          % https://q.uiver.app/?q=WzAsNixbMCwyLCJlIl0sWzEsMSwiZF8xIl0sWzEsMCwiYyJdLFswLDEsImQiXSxbMSwyLCJ0LmRfMSJdLFsyLDEsImRfMiJdLFswLDEsIiIsMCx7InN0eWxlIjp7ImhlYWQiOnsibmFtZSI6Im5vbmUifX19XSxbMSwyLCIiLDAseyJzdHlsZSI6eyJoZWFkIjp7Im5hbWUiOiJub25lIn19fV0sWzIsMywiIiwwLHsic3R5bGUiOnsiYm9keSI6eyJuYW1lIjoiYmFycmVkIn0sImhlYWQiOnsibmFtZSI6Im5vbmUifX19XSxbMywwLCIiLDAseyJzdHlsZSI6eyJoZWFkIjp7Im5hbWUiOiJub25lIn19fV0sWzEsNCwiIiwwLHsic3R5bGUiOnsiYm9keSI6eyJuYW1lIjoiYmFycmVkIn0sImhlYWQiOnsibmFtZSI6Im5vbmUifX19XSxbMiw1LCIiLDAseyJzdHlsZSI6eyJoZWFkIjp7Im5hbWUiOiJub25lIn19fV0sWzUsNCwiIiwxLHsic3R5bGUiOnsiaGVhZCI6eyJuYW1lIjoibm9uZSJ9fX1dXQ==
          \[\begin{tikzcd}[row sep = small]
            & c \\
            d & {d_1} & {d_2} \\
            e & {\gamma(d_1)}
            \arrow[no head, from=3-1, to=2-2]
            \arrow[no head, from=2-2, to=1-2]
            \arrow["\shortmid"{marking}, no head, from=1-2, to=2-1]
            \arrow[no head, from=2-1, to=3-1]
            \arrow["\shortmid"{marking}, no head, from=2-2, to=3-2]
            \arrow[no head, from=1-2, to=2-3]
            \arrow[no head, from=2-3, to=3-2]
          \end{tikzcd}\]
          \captionof{figure}{Type \rom{2}}
          \label{fig:dfc_to_zpo_hasse_4}
        \end{minipage}
      \end{figure}
      In the second case, we keep completing lozenges on the right until ending with a diagram of the shape of \ref{fig:dfc_to_zpo_hasse_5} below:
      \begin{figure}[H]
        \hspace{15pt}
        \begin{minipage}{.5\textwidth}
          % https://q.uiver.app/#q=WzAsMTAsWzAsMiwiZSJdLFsxLDEsImRfMSJdLFswLDEsImQiXSxbMSwyLCJcXGdhbW1hKGRfMSkiXSxbMiwxLCJkXzIiXSxbMiwyLCJcXGdhbW1hKGRfMikiXSxbMywxLCJkX3EiXSxbMywyLCJcXGdhbW1hKGRfcSkgPSBlJyJdLFs0LDEsImQiXSxbMiwwLCJjIl0sWzAsMSwiIiwwLHsic3R5bGUiOnsiaGVhZCI6eyJuYW1lIjoibm9uZSJ9fX1dLFsyLDAsIiIsMCx7InN0eWxlIjp7ImhlYWQiOnsibmFtZSI6Im5vbmUifX19XSxbMSwzLCIiLDAseyJzdHlsZSI6eyJib2R5Ijp7Im5hbWUiOiJiYXJyZWQifSwiaGVhZCI6eyJuYW1lIjoibm9uZSJ9fX1dLFs0LDMsIiIsMSx7InN0eWxlIjp7ImhlYWQiOnsibmFtZSI6Im5vbmUifX19XSxbNCw1LCIiLDEseyJzdHlsZSI6eyJib2R5Ijp7Im5hbWUiOiJiYXJyZWQifSwiaGVhZCI6eyJuYW1lIjoibm9uZSJ9fX1dLFs1LDYsIiIsMSx7InN0eWxlIjp7ImJvZHkiOnsibmFtZSI6ImRvdHRlZCJ9LCJoZWFkIjp7Im5hbWUiOiJub25lIn19fV0sWzYsNywiIiwxLHsic3R5bGUiOnsiYm9keSI6eyJuYW1lIjoiYmFycmVkIn0sImhlYWQiOnsibmFtZSI6Im5vbmUifX19XSxbOCw3LCIiLDEseyJzdHlsZSI6eyJib2R5Ijp7Im5hbWUiOiJiYXJyZWQifSwiaGVhZCI6eyJuYW1lIjoibm9uZSJ9fX1dLFs5LDIsIiIsMCx7ImN1cnZlIjoyLCJzdHlsZSI6eyJib2R5Ijp7Im5hbWUiOiJiYXJyZWQifSwiaGVhZCI6eyJuYW1lIjoibm9uZSJ9fX1dLFsxLDksIiIsMCx7ImN1cnZlIjotMSwic3R5bGUiOnsiaGVhZCI6eyJuYW1lIjoibm9uZSJ9fX1dLFs5LDQsIiIsMSx7InN0eWxlIjp7ImhlYWQiOnsibmFtZSI6Im5vbmUifX19XSxbOSw2LCIiLDEseyJjdXJ2ZSI6LTEsInN0eWxlIjp7ImhlYWQiOnsibmFtZSI6Im5vbmUifX19XSxbOSw4LCIiLDEseyJjdXJ2ZSI6LTIsInN0eWxlIjp7ImJvZHkiOnsibmFtZSI6ImJhcnJlZCJ9LCJoZWFkIjp7Im5hbWUiOiJub25lIn19fV1d
          \[\begin{tikzcd}[column sep = small]
            && c \\
            d & {d_1} & {d_2} & {d_q} & d \\
            e & {\gamma(d_1)} & {\gamma(d_2)} & {\gamma(d_q) = e'}
            \arrow[no head, from=3-1, to=2-2]
            \arrow[no head, from=2-1, to=3-1]
            \arrow["\shortmid"{marking}, no head, from=2-2, to=3-2]
            \arrow[no head, from=2-3, to=3-2]
            \arrow["\shortmid"{marking}, no head, from=2-3, to=3-3]
            \arrow[dotted, no head, from=3-3, to=2-4]
            \arrow["\shortmid"{marking}, no head, from=2-4, to=3-4]
            \arrow["\shortmid"{marking}, no head, from=2-5, to=3-4]
            \arrow["\shortmid"{marking}, curve={height=12pt}, no head, from=1-3, to=2-1]
            \arrow[curve={height=-6pt}, no head, from=2-2, to=1-3]
            \arrow[no head, from=1-3, to=2-3]
            \arrow[curve={height=-6pt}, no head, from=1-3, to=2-4]
            \arrow["\shortmid"{marking}, curve={height=-12pt}, no head, from=1-3, to=2-5]
          \end{tikzcd}\]
          \captionof{figure}{Path}
          \label{fig:dfc_to_zpo_hasse_5}
        \end{minipage}
        \begin{minipage}{.5\textwidth}
          % https://q.uiver.app/#q=WzAsNyxbMCwxLCJkX2oiXSxbMCwyLCJcXGdhbW1hKGRfaikiXSxbMSwxLCJkX3tqKzF9Il0sWzEsMiwiXFxnYW1tYShkX3tqKzF9KSJdLFsyLDEsImRfe2krMX0gPSBkX2oiXSxbMiwyLCJcXGdhbW1hKGRfe2krMX0pID0gXFxnYW1tYShkX2opIl0sWzEsMCwiYyJdLFswLDEsIiIsMCx7InN0eWxlIjp7ImJvZHkiOnsibmFtZSI6ImJhcnJlZCJ9LCJoZWFkIjp7Im5hbWUiOiJub25lIn19fV0sWzIsMSwiIiwxLHsic3R5bGUiOnsiaGVhZCI6eyJuYW1lIjoibm9uZSJ9fX1dLFsyLDMsIiIsMSx7InN0eWxlIjp7ImJvZHkiOnsibmFtZSI6ImJhcnJlZCJ9LCJoZWFkIjp7Im5hbWUiOiJub25lIn19fV0sWzMsNCwiIiwxLHsic3R5bGUiOnsiYm9keSI6eyJuYW1lIjoiZG90dGVkIn0sImhlYWQiOnsibmFtZSI6Im5vbmUifX19XSxbNCw1LCIiLDEseyJzdHlsZSI6eyJib2R5Ijp7Im5hbWUiOiJiYXJyZWQifSwiaGVhZCI6eyJuYW1lIjoibm9uZSJ9fX1dLFswLDYsIiIsMCx7InN0eWxlIjp7ImhlYWQiOnsibmFtZSI6Im5vbmUifX19XSxbNiwyLCIiLDEseyJzdHlsZSI6eyJoZWFkIjp7Im5hbWUiOiJub25lIn19fV0sWzYsNCwiIiwxLHsic3R5bGUiOnsiaGVhZCI6eyJuYW1lIjoibm9uZSJ9fX1dXQ==
          \[\begin{tikzcd}[column sep = small]
            & c \\
            {d_j} & {d_{j+1}} & {d_{i+1} = d_j} \\
            {\gamma(d_j)} & {\gamma(d_{j+1})} & {\gamma(d_{i+1}) = \gamma(d_j)}
            \arrow["\shortmid"{marking}, no head, from=2-1, to=3-1]
            \arrow[no head, from=2-2, to=3-1]
            \arrow["\shortmid"{marking}, no head, from=2-2, to=3-2]
            \arrow[dotted, no head, from=3-2, to=2-3]
            \arrow["\shortmid"{marking}, no head, from=2-3, to=3-3]
            \arrow[no head, from=2-1, to=1-2]
            \arrow[no head, from=1-2, to=2-2]
            \arrow[no head, from=1-2, to=2-3]
          \end{tikzcd}\]
          \captionof{figure}{Cycle}
          \label{fig:dfc_to_zpo_hasse_6}
        \end{minipage}
      \end{figure}
      First, we will never encounter a situation where $d_{i+1} = d_j$ for a $j\leq i$ as in \ref{fig:dfc_to_zpo_hasse_6} above, by acyclicity.
      And this construction must finish as above by finiteness of the poset $C$.
      Now, after renaming the elements, we end up with this situation:
      % https://q.uiver.app/#q=WzAsOSxbMCwyLCJlIl0sWzEsMSwiZFsxXSJdLFswLDEsImRbXSJdLFsxLDIsIlxcZ2FtbWEoZFsxXSkiXSxbMiwxLCJkWzJdIl0sWzIsMiwiXFxnYW1tYShkWzJdKSJdLFs0LDEsImRbcV97W119XSJdLFs0LDIsIlxcZ2FtbWEoZFtxX3tbXX1dKSJdLFsyLDAsImNbXSJdLFswLDEsIiIsMCx7InN0eWxlIjp7ImhlYWQiOnsibmFtZSI6Im5vbmUifX19XSxbMiwwLCIiLDAseyJzdHlsZSI6eyJoZWFkIjp7Im5hbWUiOiJub25lIn19fV0sWzEsMywiIiwwLHsic3R5bGUiOnsiYm9keSI6eyJuYW1lIjoiYmFycmVkIn0sImhlYWQiOnsibmFtZSI6Im5vbmUifX19XSxbNCwzLCIiLDEseyJzdHlsZSI6eyJoZWFkIjp7Im5hbWUiOiJub25lIn19fV0sWzQsNSwiIiwxLHsic3R5bGUiOnsiYm9keSI6eyJuYW1lIjoiYmFycmVkIn0sImhlYWQiOnsibmFtZSI6Im5vbmUifX19XSxbNSw2LCIiLDEseyJzdHlsZSI6eyJib2R5Ijp7Im5hbWUiOiJkb3R0ZWQifSwiaGVhZCI6eyJuYW1lIjoibm9uZSJ9fX1dLFs2LDcsIiIsMSx7InN0eWxlIjp7ImJvZHkiOnsibmFtZSI6ImJhcnJlZCJ9LCJoZWFkIjp7Im5hbWUiOiJub25lIn19fV0sWzgsMiwiIiwwLHsiY3VydmUiOjIsInN0eWxlIjp7ImJvZHkiOnsibmFtZSI6ImJhcnJlZCJ9LCJoZWFkIjp7Im5hbWUiOiJub25lIn19fV0sWzEsOCwiIiwwLHsiY3VydmUiOi0xLCJzdHlsZSI6eyJoZWFkIjp7Im5hbWUiOiJub25lIn19fV0sWzgsNCwiIiwxLHsic3R5bGUiOnsiaGVhZCI6eyJuYW1lIjoibm9uZSJ9fX1dLFs4LDYsIiIsMSx7ImN1cnZlIjotMiwic3R5bGUiOnsiaGVhZCI6eyJuYW1lIjoibm9uZSJ9fX1dXQ==
      
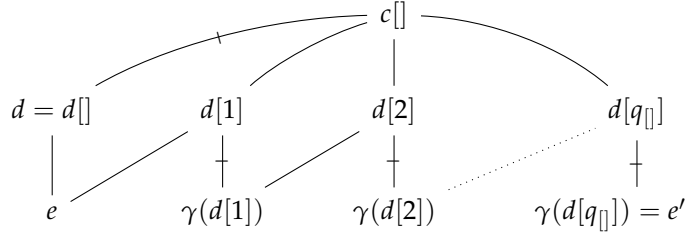
\[\begin{tikzcd}
        && {c[]} \\
        {d = d[]} & {d[1]} & {d[2]} && {d[q_{[]}]} \\
        e & {\gamma(d[1])} & {\gamma(d[2])} && {\hspace{-20pt} \gamma(d[q_{[]}]) = e'}
        \arrow[no head, from=3-1, to=2-2]
        \arrow[no head, from=2-1, to=3-1]
        \arrow["\shortmid"{marking}, no head, from=2-2, to=3-2]
        \arrow[no head, from=2-3, to=3-2]
        \arrow["\shortmid"{marking}, no head, from=2-3, to=3-3]
        \arrow[dotted, no head, from=3-3, to=2-5]
        \arrow["\shortmid"{marking}, no head, from=2-5, to=3-5]
        \arrow["\shortmid"{marking}, curve={height=12pt}, no head, from=1-3, to=2-1]
        \arrow[curve={height=-6pt}, no head, from=2-2, to=1-3]
        \arrow[no head, from=1-3, to=2-3]
        \arrow[curve={height=-12pt}, no head, from=1-3, to=2-5]
      \end{tikzcd}\]
      If each $d[i]$ is in $\Lambda_{k+1}$ then we are done.
      In the other case, for example if $d[2]\prec^{+}c[2]$, then we can once again \textit{unfold}
      % https://q.uiver.app/#q=WzAsMTQsWzUsMiwiXFxnYW1tYShkWzJdKSJdLFs2LDEsImRbMiwxXSJdLFs1LDEsImRbMl0iXSxbNiwyLCJcXGdhbW1hKGRbMiwxXSkiXSxbNywxLCJkWzIsMl0iXSxbNywyLCJcXGdhbW1hKGRbMiwyXSkiXSxbOSwxLCJkWzIscV97WzJdfV0iXSxbOSwyLCJcXGdhbW1hKGRbMixxX3tbMl19XSkiXSxbNywwLCJjWzJdIl0sWzEsMCwiY1syXSJdLFsxLDEsImRbMl0iXSxbMCwyLCJcXGdhbW1hKGRbMV0pIl0sWzEsMiwiXFxnYW1tYShkWzJdKSJdLFszLDEsIlxcbWJveHsgYXMgfSJdLFswLDEsIiIsMCx7InN0eWxlIjp7ImhlYWQiOnsibmFtZSI6Im5vbmUifX19XSxbMiwwLCIiLDAseyJzdHlsZSI6eyJoZWFkIjp7Im5hbWUiOiJub25lIn19fV0sWzEsMywiIiwwLHsic3R5bGUiOnsiYm9keSI6eyJuYW1lIjoiYmFycmVkIn0sImhlYWQiOnsibmFtZSI6Im5vbmUifX19XSxbNCwzLCIiLDEseyJzdHlsZSI6eyJoZWFkIjp7Im5hbWUiOiJub25lIn19fV0sWzQsNSwiIiwxLHsic3R5bGUiOnsiYm9keSI6eyJuYW1lIjoiYmFycmVkIn0sImhlYWQiOnsibmFtZSI6Im5vbmUifX19XSxbNSw2LCIiLDEseyJzdHlsZSI6eyJib2R5Ijp7Im5hbWUiOiJkb3R0ZWQifSwiaGVhZCI6eyJuYW1lIjoibm9uZSJ9fX1dLFs2LDcsIiIsMSx7InN0eWxlIjp7ImJvZHkiOnsibmFtZSI6ImJhcnJlZCJ9LCJoZWFkIjp7Im5hbWUiOiJub25lIn19fV0sWzgsMiwiIiwwLHsiY3VydmUiOjIsInN0eWxlIjp7ImJvZHkiOnsibmFtZSI6ImJhcnJlZCJ9LCJoZWFkIjp7Im5hbWUiOiJub25lIn19fV0sWzEsOCwiIiwwLHsiY3VydmUiOi0xLCJzdHlsZSI6eyJoZWFkIjp7Im5hbWUiOiJub25lIn19fV0sWzgsNCwiIiwxLHsic3R5bGUiOnsiaGVhZCI6eyJuYW1lIjoibm9uZSJ9fX1dLFs4LDYsIiIsMSx7ImN1cnZlIjotMiwic3R5bGUiOnsiaGVhZCI6eyJuYW1lIjoibm9uZSJ9fX1dLFs5LDEwLCIiLDEseyJzdHlsZSI6eyJib2R5Ijp7Im5hbWUiOiJiYXJyZWQifSwiaGVhZCI6eyJuYW1lIjoibm9uZSJ9fX1dLFsxMCwxMSwiIiwxLHsic3R5bGUiOnsiaGVhZCI6eyJuYW1lIjoibm9uZSJ9fX1dLFsxMCwxMiwiIiwxLHsic3R5bGUiOnsiYm9keSI6eyJuYW1lIjoiYmFycmVkIn0sImhlYWQiOnsibmFtZSI6Im5vbmUifX19XV0=
      
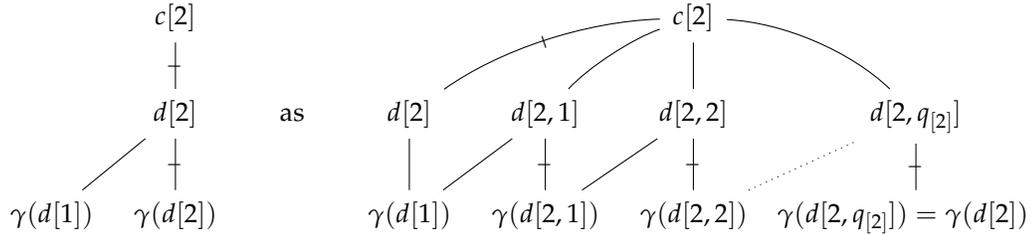
\[\begin{tikzcd}[column sep = tiny]
        & {c[2]} &&&&&& {c[2]} \\
        & {d[2]} && {\mbox{ as }} && {d[2]} & {d[2,1]} & {d[2,2]} && {d[2,q_{[2]}]} \\
        {\gamma(d[1])} & {\gamma(d[2])} &&&& {\gamma(d[1])} & {\gamma(d[2,1])} & {\gamma(d[2,2])} && {\hspace{-10pt} \gamma(d[2,q_{[2]}]) = \gamma(d[2])}
        \arrow[no head, from=3-6, to=2-7]
        \arrow[no head, from=2-6, to=3-6]
        \arrow["\shortmid"{marking}, no head, from=2-7, to=3-7]
        \arrow[no head, from=2-8, to=3-7]
        \arrow["\shortmid"{marking}, no head, from=2-8, to=3-8]
        \arrow[dotted, no head, from=3-8, to=2-10]
        \arrow["\shortmid"{marking}, no head, from=2-10, to=3-10]
        \arrow["\shortmid"{marking}, curve={height=12pt}, no head, from=1-8, to=2-6]
        \arrow[curve={height=-6pt}, no head, from=2-7, to=1-8]
        \arrow[no head, from=1-8, to=2-8]
        \arrow[curve={height=-12pt}, no head, from=1-8, to=2-10]
        \arrow["\shortmid"{marking}, no head, from=1-2, to=2-2]
        \arrow[no head, from=2-2, to=3-1]
        \arrow["\shortmid"{marking}, no head, from=2-2, to=3-2]
      \end{tikzcd}\]
      Now we may iterate this unfolding process. This will produce a tree-shaped collection of $d$'s and $c$'s, with relations
      $$d[a_1,\,\cdots,\,a_p]\prec^{+}c[a_1,\,\cdots,\,a_p] \qquad \mbox{and} \qquad c[a_1,\,\cdots,\,a_p]\succ^{-}d[a_1,\,\cdots,\,a_p,\,a]$$ for each $a\in\left[\hspace{-3pt}\left[1,\,q_{[a_1,\,\cdots,\,a_p]}\right]\hspace{-3pt}\right]$.
      Each branch $$d[]\prec^{+}c[]\succ^{-}d[a_1]\prec^{+}c[a_1]\succ^{-}d[a_1,\,a_2]\prec^{+}\cdots\succ^{-}d[a_1,\,\cdots,\,a_p]\prec^{+}\cdots$$
      must be finite, otherwise it would contradict $\mathcal{P}_{k+1}$.
      \textit{i.e.} at some point, all $d[a_1,\,\cdots,\,a_p]$ are in $\Lambda_{k+1}$.
      We thus obtain a path of the desired form, which completes the proof of \reflem{lem:DFC_path_in_Lambda}.
    \qedhere
  \end{itemize}
\end{proof}

\noindent Notice that through the proof above, we have obtained the following result (\textit{c.f.} $\mathcal{R}_k$):

\begin{lem}\label{lem:DFC_Lambda_sources_of_iterated_target}
  $\forall k\in\intInter{0}{n-1},\quad \Lambda_k \subseteq \delta(\gamma^{(k+1)}\omega)$.
\end{lem}

\begin{lem}\label{lem:facet_of_Lambda}
  If $e \prec^{\alpha} d$ for some $d \in C$ and $\alpha \in \{+,\,-\}$, then $e \prec^{\alpha} d'$ for a unique $d' \in \Lambda$.\\[5pt]
  Moreover, for $k<n$, we have $C_k \setminus \{\gamma^{(k)}\omega\} \subseteq (\delta^{-}\cup\gamma^{+})(\Lambda_{k+1})$.
\end{lem}
\begin{proof}
  The existence part of the first statement has been seen in the proof of \refprop{prop:DFC_strictness_1} above (\textit{c.f.} $\mathcal{Q}_k$).\\
  We prove the uniqueness part. It is clear if $k=n-1$, hence we assume $k\le n-2$.
  Suppose $c\prec^{\alpha}b$ and $c\prec^{\alpha}b'$, with $b,\,b'\in\Lambda_{k+1}$.
  By using \reflem{lem:DFC_Lambda_sources_of_iterated_target}, either $\gamma^{(k+2)}\omega$ is not a loop and we may construct a lozenge $(e,\,d,\,\gamma^{(k+2)}\omega,\,d')$ which, by the sign rule, implies $d = d'$.
  Or $\gamma^{(k+2)}\omega$ is a loop on $\gamma^{(k+1)}\omega$, hence $d = d' = \gamma^{(k+1)}\omega$.

  We then prove the second point. Let $\gamma^{(k)}\omega \neq e \prec^{\o} d$ for some $e,\,d \in C$.
  By (loops), there is some $c$ with $d \prec^{+} c$. Using the first point, we may assume $c \in \Lambda_{k+2}$.\\
  Since $e \neq \gamma^{(k)}\omega$, we also have $c \neq \gamma^{(k+2)}\omega$.
  Hence $k < n-2$, and \reflem{lem:DFC_Lambda_sources_of_iterated_target} yields $c \in \delta(\gamma^{(k+3)}\omega)$.
  Necessarily $c \nprec^{\o} \gamma^{(k+3)}\omega$, otherwise we would have $d = \gamma^{(k+1)}\omega$ which would contradict the hypothesis. Thus $c \prec^{-} \gamma^{(k+3)}\omega$ and we are in the situation $e \prec^{\o} d \prec^{+} c \prec^{-} \gamma^{(k+3)}\omega$.
  By successively completing lozenges from left to right, we may construct a diagram as below ($\alpha \beta \neq \o$):
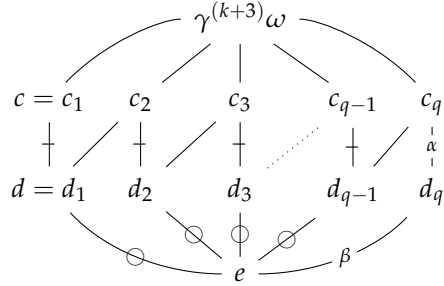
  \[\begin{tikzcd}[row sep = 15pt, column sep = tiny]
    && {\gamma^{(k+3)}\omega} \\
    {c=c_1} & {c_2} & {c_3} & {c_{q-1}} & {c_q} \\
    {d=d_1} & {d_2} & {d_3} & {d_{q-1}} & {d_q} \\
    && e
    \arrow["\Circle"{marking}, curve={height=12pt}, no head, from=3-1, to=4-3]
    \arrow[curve={height=12pt}, no head, from=1-3, to=2-1]
    \arrow["\shortmid"{marking}, no head, from=2-1, to=3-1]
    \arrow[no head, from=1-3, to=2-2]
    \arrow[no head, from=2-2, to=3-1]
    \arrow["\shortmid"{marking}, no head, from=2-2, to=3-2]
    \arrow["\Circle"{marking}, no head, from=3-2, to=4-3]
    \arrow[no head, from=1-3, to=2-3]
    \arrow[no head, from=2-3, to=3-2]
    \arrow["\shortmid"{marking}, no head, from=2-3, to=3-3]
    \arrow["\Circle"{marking}, no head, from=3-3, to=4-3]
    \arrow[no head, from=1-3, to=2-4]
    \arrow[dotted, no head, from=2-4, to=3-3]
    \arrow["\shortmid"{marking}, no head, from=2-4, to=3-4]
    \arrow["\Circle"{marking}, no head, from=3-4, to=4-3]
    \arrow[curve={height=-12pt}, no head, from=1-3, to=2-5]
    \arrow[no head, from=2-5, to=3-4]
    \arrow["\alpha"{description}, no head, from=2-5, to=3-5]
    \arrow["\beta"{description}, curve={height=-12pt}, no head, from=3-5, to=4-3]
  \end{tikzcd}\]
  Hence $e \prec^{\beta} d_q$ and we are indeed reduced to the first statement.
\end{proof}

\begin{lem}\label{lem:DFC_not_a_source_is_iterated_target}
  For all $k\le n$, $C_{k}\setminus\delta^{-}(C_{k+1}) = \{\gamma^{(k)}\omega\}$.
\end{lem}
\begin{proof} ~\\[-15pt]
  \begin{itemize}
    \item \boxed{\mbox{case }k=n}\:
      In this case, $d = \omega$ is the only $n$-dimensional cell, and it is not in $\delta^{-}(C)$.
    \item \boxed{\mbox{case }k<n}\:
      Suppose $d\in C_k\setminus\delta^{-}(C_{k+1})$ is not $\gamma^{(k)}\omega$.
      Then there is some $c\in C_{k+1}$ with $d = \gamma(c)$, and by \reflem{lem:facet_of_Lambda}, we may assume $c\in\Lambda$.\\
      Necessarily $c \neq \gamma^{(k+1)}\omega$ because we assumed $d \neq \gamma^{(k)}\omega$. Hence $k \le n-2$ and $c \prec^{-} \gamma^{(k+2)}\omega$ by \reflem{lem:DFC_Lambda_sources_of_iterated_target}.
      By remembering that $d \notin \delta^{-}(C)$, completing $d = \gamma(c) \prec c \prec^{-} \gamma^{(k+2)}\omega$ yields $d = \gamma^{(k)}\omega$.

      Conversely: if $k = n-1$, then $\gamma(\omega) \nprec^{-} \omega$, whence $\gamma(\omega) \notin \delta^{-}(C_n)$. \\
      If $k < n-1$, suppose $\gamma^{(k)}\omega \prec^{-} c$ for some $c$. By \reflem{lem:facet_of_Lambda}, we may assume $c \in \Lambda$, hence $c \in \delta(\gamma^{k+2}\omega)$ by \reflem{lem:DFC_Lambda_sources_of_iterated_target}.
      We also have $\gamma^{(k)}\omega = \gamma^2(\gamma^{(k+2)}\omega)$. If $c \prec^{\o} \gamma^{(k+2)}\omega$, then $\gamma^{(k)}\omega \in \delta^{-}(\gamma^{(k+1)}\omega)$ which is impossible.
      If $c \prec^{-} \gamma^{(k+2)}\omega$, then we are in one of the two possible situations below.
      % https://q.uiver.app/#q=WzAsOCxbMSwwLCJcXGdhbW1hXnsoaysyKX1cXG9tZWdhIl0sWzAsMSwiYyJdLFsxLDIsIlxcZ2FtbWFeeyhrKX1cXG9tZWdhIl0sWzIsMSwiXFxnYW1tYV57KGsrMSl9XFxvbWVnYSJdLFs1LDAsIlxcZ2FtbWFeeyhrKzIpfVxcb21lZ2EiXSxbNCwxLCJjIl0sWzUsMiwiXFxnYW1tYV57KGspfVxcb21lZ2EiXSxbNiwxLCJcXGdhbW1hXnsoaysxKX1cXG9tZWdhIl0sWzAsMSwiIiwxLHsic3R5bGUiOnsiaGVhZCI6eyJuYW1lIjoibm9uZSJ9fX1dLFsxLDIsIiIsMSx7InN0eWxlIjp7ImhlYWQiOnsibmFtZSI6Im5vbmUifX19XSxbMCwzLCIiLDEseyJzdHlsZSI6eyJib2R5Ijp7Im5hbWUiOiJiYXJyZWQifSwiaGVhZCI6eyJuYW1lIjoibm9uZSJ9fX1dLFszLDIsIiIsMSx7InN0eWxlIjp7ImJvZHkiOnsibmFtZSI6ImJhcnJlZCJ9LCJoZWFkIjp7Im5hbWUiOiJub25lIn19fV0sWzQsNSwiIiwxLHsic3R5bGUiOnsiaGVhZCI6eyJuYW1lIjoibm9uZSJ9fX1dLFs0LDcsIiIsMSx7InN0eWxlIjp7ImJvZHkiOnsibmFtZSI6ImJhcnJlZCJ9LCJoZWFkIjp7Im5hbWUiOiJub25lIn19fV0sWzUsNiwiIiwxLHsic3R5bGUiOnsiaGVhZCI6eyJuYW1lIjoibm9uZSJ9fX1dLFs3LDYsIiIsMSx7InN0eWxlIjp7ImhlYWQiOnsibmFtZSI6Im5vbmUifX19XV0=
      
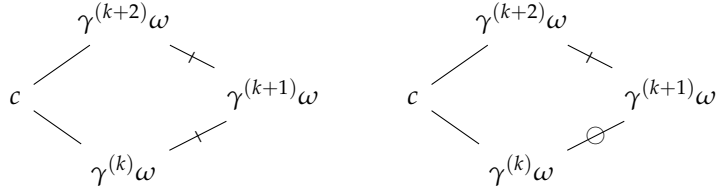
\[\begin{tikzcd}[column sep = small, row sep = small]
        & {\gamma^{(k+2)}\omega} &&&& {\gamma^{(k+2)}\omega} \\
        c && {\gamma^{(k+1)}\omega} && c && {\gamma^{(k+1)}\omega} \\
        & {\gamma^{(k)}\omega} &&&& {\gamma^{(k)}\omega}
        \arrow[no head, from=1-2, to=2-1]
        \arrow[no head, from=2-1, to=3-2]
        \arrow["\shortmid"{marking}, no head, from=1-2, to=2-3]
        \arrow["\shortmid"{marking}, no head, from=2-3, to=3-2]
        \arrow[no head, from=1-6, to=2-5]
        \arrow["\shortmid"{marking}, no head, from=1-6, to=2-7]
        \arrow[no head, from=2-5, to=3-6]
        \arrow["\Circle"{marking}, no head, from=2-7, to=3-6]
      \end{tikzcd}\]
      The left one violates the sign rule, and the right one contradicts \refprop{prop:source_confinement}.
      \qedhere
  \end{itemize}
\end{proof}

\begin{lem}\label{lem:Lambda_is_iterated_target_sources}
  For any $k<n$, $\Lambda_k = \delta(\gamma^{(k+1)}\omega)$.
\end{lem}
\begin{proof}
  One inclusion is by \reflem{lem:DFC_Lambda_sources_of_iterated_target}, so we show the other one.
  Let $d \in \delta(\gamma^{(k+1)}\omega)$, and suppose that $d \prec^{+} c$ for some $c$. By \reflem{lem:facet_of_Lambda}, we may assume $c \in \Lambda_{k+1}$, hence $c \in \delta(\gamma^{(k+2)}\omega)$ by \reflem{lem:DFC_Lambda_sources_of_iterated_target}.
  \begin{itemize}
    \item Suppose $d \prec^{-} \gamma^{(k+1)}\omega$.
      Either $\gamma^{(k+2)}\omega$ is a loop, hence $c = \gamma^{(k+1)}\omega$, so $d \prec^{-} c$ and $d \prec^{+} c$, which is impossible.
      Or $c \prec^{-} \gamma^{(k+2)}\omega$, $\gamma^{(k+1)}\omega \prec^{+} \gamma^{(k+2)}\omega$, and there is a lozenge $(d,\,c,\,\gamma^{(k+2)}\omega,\,\gamma^{(k+1)}\omega)$ violating the sign rule.
    \item Suppose $d \prec^{\o} \gamma^{(k+1)}\omega$.
      Either $\gamma^{(k+2)}\omega$ is a loop, hence $c = \gamma^{(k+1)}\omega$, so $d \prec^{-} c$ and $d \prec^{\o} c$, which is impossible.
      Or $c \prec^{-} \gamma^{(k+2)}\omega$, $\gamma^{(k+1)}\omega \prec^{+} \gamma^{(k+2)}\omega$, and there is a lozenge $(d,\,c,\,\gamma^{(k+2)}\omega,\,\gamma^{(k+1)}\omega)$ contradicting \refprop{prop:source_confinement}.
    \qedhere
  \end{itemize}
\end{proof}

\begin{prop}[Tree structure $\Delta_a$]\label{prop:DFC_tree_structure}
  For all $a \in C$, there is an associated tree structure $\Delta_a$ as follows.
  \begin{itemize}
    \item The nodes are $\Delta_a^\bullet = \delta(a) \setminus \Omega$, and the edges are $T^{\edge} = \{\gamma^2(a)\} \cup \partial\Delta_a^\bullet$.
    \item There is a relation $a \branch{b}$ when $b \prec^{-} a$ and a relation $\branch{b} a$ when $b \prec^{+} a$.
    \item The root edge is $\gamma^2(a)$.
  \end{itemize}
\end{prop}
\begin{proof}
  When $a \in \Omega$, $\Delta_a$ is either reduced to its root if $\gamma(a)\in\Omega$, or a corolla with $\delta^{-}(\gamma(a))$ as leaves if $\gamma(a) \notin \Omega$.\\
  Suppose now $a \notin \Omega$.
  Each node $b \in \Delta_a^\bullet$ is paired with exactly one target edge $\gamma(b) \prec^{+} b$. The sign rule shows the uniqueness of the target node, and also shows that $\gamma^2(a)$ has no target node.
  Let $c \neq \gamma^2(a)$ be an edge of $\Delta_a$, there is some $b \in \Delta_a^\bullet$ and $c \prec^{\alpha} b$ for some $\alpha \neq \o$. By oriented thinness, we may assume $c \prec^{-} b$.
  Now, completing lozenges from left to right as below yields the descending path from $c$ to the root $\gamma^2(a)$.
  \vspace{-15pt}
  % https://q.uiver.app/#q=WzAsMTEsWzMsMCwiYSJdLFsxLDEsImIgPSBiXzEiXSxbMSwyLCJcXGdhbW1hKGJfMSkiXSxbMiwxLCJiXzIiXSxbMiwyLCJcXGdhbW1hKGJfMikiXSxbMywxLCJiXzMiXSxbMywyLCJcXGdhbW1hKGJfMykiXSxbNCwxLCJiX3twfSJdLFs0LDIsIlxcZ2FtbWFeMihhKSJdLFs1LDEsIlxcZ2FtbWEoYSkiXSxbMCwyLCJjIl0sWzEsMiwiIiwwLHsic3R5bGUiOnsiYm9keSI6eyJuYW1lIjoiYmFycmVkIn0sImhlYWQiOnsibmFtZSI6Im5vbmUifX19XSxbMywyLCIiLDIseyJzdHlsZSI6eyJoZWFkIjp7Im5hbWUiOiJub25lIn19fV0sWzMsNCwiIiwwLHsic3R5bGUiOnsiYm9keSI6eyJuYW1lIjoiYmFycmVkIn0sImhlYWQiOnsibmFtZSI6Im5vbmUifX19XSxbNCw1LCIiLDAseyJzdHlsZSI6eyJoZWFkIjp7Im5hbWUiOiJub25lIn19fV0sWzUsNiwiIiwwLHsic3R5bGUiOnsiYm9keSI6eyJuYW1lIjoiYmFycmVkIn0sImhlYWQiOnsibmFtZSI6Im5vbmUifX19XSxbNyw2LCIiLDIseyJzdHlsZSI6eyJib2R5Ijp7Im5hbWUiOiJkb3R0ZWQifSwiaGVhZCI6eyJuYW1lIjoibm9uZSJ9fX1dLFswLDEsIiIsMCx7ImN1cnZlIjoxLCJzdHlsZSI6eyJoZWFkIjp7Im5hbWUiOiJub25lIn19fV0sWzAsMywiIiwwLHsic3R5bGUiOnsiaGVhZCI6eyJuYW1lIjoibm9uZSJ9fX1dLFswLDUsIiIsMCx7InN0eWxlIjp7ImhlYWQiOnsibmFtZSI6Im5vbmUifX19XSxbMCw3LCIiLDAseyJzdHlsZSI6eyJoZWFkIjp7Im5hbWUiOiJub25lIn19fV0sWzcsOCwiIiwwLHsic3R5bGUiOnsiYm9keSI6eyJuYW1lIjoiYmFycmVkIn0sImhlYWQiOnsibmFtZSI6Im5vbmUifX19XSxbOCw5LCIiLDAseyJzdHlsZSI6eyJib2R5Ijp7Im5hbWUiOiJiYXJyZWQifSwiaGVhZCI6eyJuYW1lIjoibm9uZSJ9fX1dLFswLDksIiIsMix7ImN1cnZlIjotMSwic3R5bGUiOnsiYm9keSI6eyJuYW1lIjoiYmFycmVkIn0sImhlYWQiOnsibmFtZSI6Im5vbmUifX19XSxbMSwxMCwiIiwyLHsic3R5bGUiOnsiaGVhZCI6eyJuYW1lIjoibm9uZSJ9fX1dXQ==
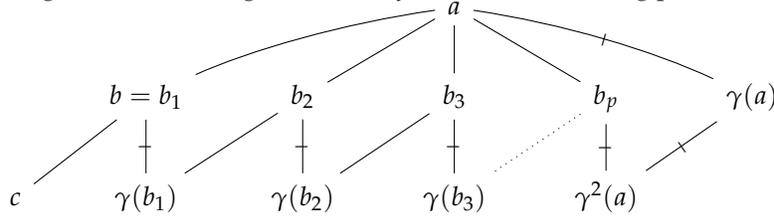
  \[\begin{tikzcd}
    &&& a \\
    & {b = b_1} & {b_2} & {b_3} & {b_{p}} & {\gamma(a)} \\
    c & {\gamma(b_1)} & {\gamma(b_2)} & {\gamma(b_3)} & {\gamma^2(a)}
    \arrow["\shortmid"{marking}, no head, from=2-2, to=3-2]
    \arrow[no head, from=2-3, to=3-2]
    \arrow["\shortmid"{marking}, no head, from=2-3, to=3-3]
    \arrow[no head, from=3-3, to=2-4]
    \arrow["\shortmid"{marking}, no head, from=2-4, to=3-4]
    \arrow[dotted, no head, from=2-5, to=3-4]
    \arrow[curve={height=6pt}, no head, from=1-4, to=2-2]
    \arrow[no head, from=1-4, to=2-3]
    \arrow[no head, from=1-4, to=2-4]
    \arrow[no head, from=1-4, to=2-5]
    \arrow["\shortmid"{marking}, no head, from=2-5, to=3-5]
    \arrow["\shortmid"{marking}, no head, from=3-5, to=2-6]
    \arrow["\shortmid"{marking}, curve={height=-6pt}, no head, from=1-4, to=2-6]
    \arrow[no head, from=2-2, to=3-1]
  \end{tikzcd}\]
\end{proof}

\begin{prop}[Tree structure $T_k$]\label{prop:DFC_tree_structure_Lambda}
  Let $C$ be a DFC and $2 \le k \le n+2$. There is a rooted tree structure $T_k$ defined as follows.
  \begin{itemize}
    \item The set of nodes is $\Lambda_{k-1}$, and the set of edges is $C_{k-2}$.
    \item There is a relation $x \branch{y}$ iff $y \prec^{-} x$, and a relation $\branch{y} x$ iff $y \prec^{+} x$.
    \item The root edge is $\gamma^{(k-2)}(\omega)$.
  \end{itemize}
\end{prop}
\begin{proof}
  Suppose $k \le n$.\\
  Because of \reflem{lem:facet_of_Lambda}, $C_{k-2} = \partial\Lambda_{k-1} \cup \{\gamma^{(k-2)}\omega\}$.
  Notice also that elements $y \in \Lambda_{k-1}$ are not loops, because there is no $x$ with $y \prec^{+} x$. Hence $\Lambda_{k-1} = \delta(\gamma^{(k)}\omega)\setminus\Omega$ according to \reflem{lem:Lambda_is_iterated_target_sources}.
  Hence $T_k$ coincides with $\Delta_{\gamma^{(k)}\omega}$ as defined in \refprop{prop:DFC_tree_structure}.

  Consider now the case $k = n+1$. There is a unique node $\omega$ with a set of leaves $T_{n+1}^\vert = \delta(\omega)$ and $\gamma(\omega)$ as root. It is a well defined rooted tree structure.

  When $k = n+2$, there is no node and a unique leaf $\omega$: $T_{n+2}$ is a unit tree.
\end{proof}

\begin{lem}\label{lem:blackdots_leaves}
  Let $2 \le k < n+2$. We have the equality $\msfbdots(T_k) = \msfleaves(T_{k+1})$.
\end{lem}
\begin{proof}
  $\msfbdots(T_k) = T_k^\bullet = \Lambda_{k-1} = C_{k-1} \setminus \gamma^{+}(C_{k}) \overset{\mbox{\reflem{lem:facet_of_Lambda}}}{=} C_{k-1} \setminus \gamma^{+}(\Lambda_{k}) = T_{k+1}^\vert = \msfleaves(T_{k+1})$.
\end{proof}

\begin{prop}[Hexagon property]
  For every hexagon as follows, with $\alpha,\,\alpha',\,\beta,\,\beta' \neq \o$.
  % https://q.uiver.app/?q=WzAsNixbMCwyLCJkIl0sWzIsMiwiZCciXSxbMCwxLCJjIl0sWzIsMSwiYyciXSxbMSwwLCJiIl0sWzEsMywiZSJdLFswLDUsIlxcYmV0YSIsMCx7InN0eWxlIjp7ImhlYWQiOnsibmFtZSI6Im5vbmUifX19XSxbMSw1LCJcXGJldGEnIiwyLHsic3R5bGUiOnsiaGVhZCI6eyJuYW1lIjoibm9uZSJ9fX1dLFsyLDAsIlxcYWxwaGEiLDAseyJzdHlsZSI6eyJoZWFkIjp7Im5hbWUiOiJub25lIn19fV0sWzMsMSwiXFxhbHBoYSciLDIseyJzdHlsZSI6eyJoZWFkIjp7Im5hbWUiOiJub25lIn19fV0sWzQsMiwiIiwxLHsic3R5bGUiOnsiaGVhZCI6eyJuYW1lIjoibm9uZSJ9fX1dLFs0LDMsIiIsMSx7InN0eWxlIjp7ImhlYWQiOnsibmFtZSI6Im5vbmUifX19XV0=
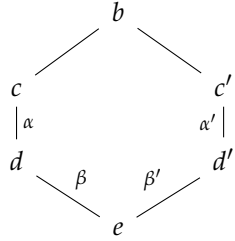
  \[\begin{tikzcd}[row sep = 12pt]
    & b \\
    c && {c'} \\
    d && {d'} \\
    & e
    \arrow["\beta", no head, from=3-1, to=4-2]
    \arrow["{\beta'}"', no head, from=3-3, to=4-2]
    \arrow["\alpha", no head, from=2-1, to=3-1]
    \arrow["{\alpha'}"', no head, from=2-3, to=3-3]
    \arrow[no head, from=1-2, to=2-1]
    \arrow[no head, from=1-2, to=2-3]
  \end{tikzcd}\]
  Either $c = c'$, or (potentially by exchanging the role of $c$ and $c'$)
  there is a simple zig-zag in $\Delta_b$ as below (with $p \ge 1$):
  
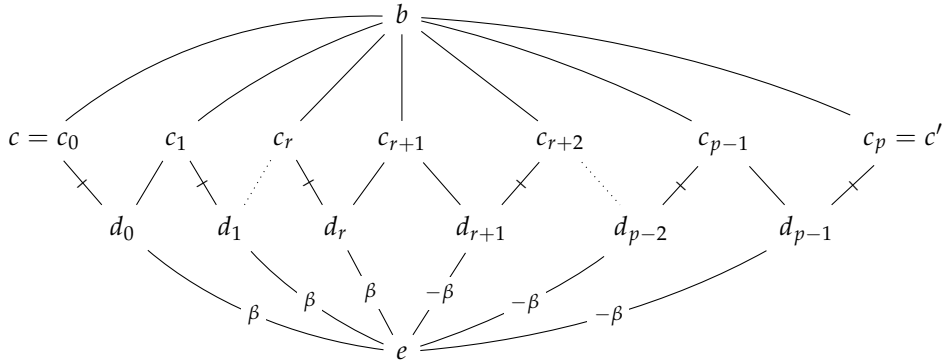
\[\begin{tikzcd}[column sep = 2pt, row sep = 15pt]
    &&&&&& b \\
    \\
    {c = c_0} && {c_1} && {c_r} && {c_{r+1}} && {c_{r+2}} && {c_{p-1}} && {c_p = c'} \\
    & {d_0} && {d_1} && {d_r} && {d_{r+1}} && {d_{p-2}} && {d_{p-1}} \\
    \\
    &&&&&& e
    \arrow["\shortmid"{marking}, no head, from=3-1, to=4-2]
    \arrow[no head, from=3-3, to=4-2]
    \arrow["\shortmid"{marking}, no head, from=3-3, to=4-4]
    \arrow["\shortmid"{marking}, no head, from=3-5, to=4-6]
    \arrow[no head, from=3-7, to=4-6]
    \arrow[no head, from=3-7, to=4-8]
    \arrow["\shortmid"{marking}, no head, from=4-8, to=3-9]
    \arrow["\shortmid"{marking}, no head, from=4-10, to=3-11]
    \arrow[no head, from=3-11, to=4-12]
    \arrow["\shortmid"{marking}, no head, from=4-12, to=3-13]
    \arrow[curve={height=-18pt}, no head, from=3-1, to=1-7]
    \arrow[curve={height=-6pt}, no head, from=3-3, to=1-7]
    \arrow[no head, from=3-5, to=1-7]
    \arrow[no head, from=3-7, to=1-7]
    \arrow[no head, from=1-7, to=3-9]
    \arrow[curve={height=-6pt}, no head, from=1-7, to=3-11]
    \arrow[curve={height=-12pt}, no head, from=1-7, to=3-13]
    \arrow["\beta"{description}, no head, from=4-6, to=6-7]
    \arrow["{-\beta}"{description}, no head, from=4-8, to=6-7]
    \arrow["{-\beta}"{description}, curve={height=-6pt}, no head, from=4-10, to=6-7]
    \arrow["{-\beta}"{description}, curve={height=-12pt}, no head, from=4-12, to=6-7]
    \arrow["\beta"{description}, curve={height=6pt}, no head, from=4-4, to=6-7]
    \arrow["\beta"{description}, curve={height=12pt}, no head, from=4-2, to=6-7]
    \arrow[dotted, no head, from=4-4, to=3-5]
    \arrow[dotted, no head, from=3-9, to=4-10]
  \end{tikzcd}\]
  \label{prop:DFC_hexagon}
\end{prop}
\begin{proof}
  It is shown exactly as the Hexagon property (Proposition $3.9$) in \cite{leclerc2023}.
\end{proof}

\begin{prop}[pencil linearity]\label{prop:DFC_pencil_linearity}
  $C$ satisfies the axiom of pencil linearity introduced in \cite{zawadowski2023positive}.
  That is:\\ $\forall k>0,\,\forall e\in C_{k-1},\,\forall \beta \in\{+,\,-\},\quad\left\{d\in C_k \mid e\prec^{\beta}d\right\}$ is linearly ordered by $<^{+}$.
\end{prop}
\begin{proof}
  By the same arguments as in \cite{leclerc2023}, \textit{mutatis mutandis}.
\end{proof}

\para As DFC, zoom-complexes carry an order for loops.
  For instance, in \refexm{exm:opetope_zoom_loops}, the subdivision of the edge $c_1$: $(a_7 < a_5 < a_4 < a_3)$ tells us in which order the cells of $\varrho$ should be composed.
  Notice that this order is not immediately part of the DFC structure defined in \refexm{exm:opetope_loops_poset}, but it may be rebuilt out of the locals order $\triangleleft$.
  To begin with, we will define in \refdefin{defin:loops_order} an order $\blacktriangleleft^1$, allowing us to compare $a_7$ to $a_3,\,a_4$ and $a_5$.
  Then in \refdefin{defin:loops_total_order}, we will introduce an order $\blacktriangleleft^2$, allowing us to compare $a_5$ to $a_3$ and $a_4$. Finally, from $\blacktriangleleft^1$ and $\blacktriangleleft^2$ we will obtain a total order $\blacktriangleleft$, allowing us to compare the loops on a given cell which are minimal for the inclusion.

\begin{defin}[Loops order]\label{defin:loops_order}
  For any chains $c \prec^{\o} b \prec^{+} a$ and $c \prec^{\o} b' \prec^{-} a$ with $a\in\Lambda$, we say that the loop $b'$ is \emph{included} in the loop $b$, which we write $b' \subset b$.
  We write $\subseteq$ for the reflexive transitive closure of this relation.
  
  Notice that for any $0 \le k \le n-2$ and $c\in C_k$, by Propositions \ref{prop:DFC_hexagon} and \ref{prop:DFC_pencil_linearity} there is a zig-zag in $T_{k+2}$ as follows:
  % https://q.uiver.app/#q=WzAsMTIsWzEsMCwiYV8xIl0sWzAsMSwiYl8xIl0sWzIsMCwiYV97cC0xfSJdLFsxLDEsImJfMiJdLFsyLDEsImJfcCJdLFszLDAsImFfcCJdLFs0LDEsImJfe3ArMX0iXSxbNCwwLCJhX3twKzF9Il0sWzUsMSwiYl97cS0xfSJdLFs1LDAsImFfe3EtMX0iXSxbNiwxLCJiX3EiXSxbMywyLCJjIl0sWzAsMSwiIiwwLHsic3R5bGUiOnsiaGVhZCI6eyJuYW1lIjoibm9uZSJ9fX1dLFswLDMsIiIsMix7InN0eWxlIjp7ImJvZHkiOnsibmFtZSI6ImJhcnJlZCJ9LCJoZWFkIjp7Im5hbWUiOiJub25lIn19fV0sWzIsMywiIiwwLHsic3R5bGUiOnsiaGVhZCI6eyJuYW1lIjoibm9uZSJ9fX1dLFsyLDQsIiIsMix7InN0eWxlIjp7ImJvZHkiOnsibmFtZSI6ImJhcnJlZCJ9LCJoZWFkIjp7Im5hbWUiOiJub25lIn19fV0sWzQsNSwiIiwyLHsic3R5bGUiOnsiaGVhZCI6eyJuYW1lIjoibm9uZSJ9fX1dLFs1LDYsIiIsMix7InN0eWxlIjp7ImhlYWQiOnsibmFtZSI6Im5vbmUifX19XSxbNiw3LCIiLDIseyJzdHlsZSI6eyJib2R5Ijp7Im5hbWUiOiJiYXJyZWQifSwiaGVhZCI6eyJuYW1lIjoibm9uZSJ9fX1dLFs3LDgsIiIsMix7InN0eWxlIjp7ImhlYWQiOnsibmFtZSI6Im5vbmUifX19XSxbOCw5LCIiLDIseyJzdHlsZSI6eyJib2R5Ijp7Im5hbWUiOiJiYXJyZWQifSwiaGVhZCI6eyJuYW1lIjoibm9uZSJ9fX1dLFs5LDEwLCIiLDIseyJzdHlsZSI6eyJoZWFkIjp7Im5hbWUiOiJub25lIn19fV0sWzEsMTEsIiIsMCx7ImN1cnZlIjoyLCJzdHlsZSI6eyJoZWFkIjp7Im5hbWUiOiJub25lIn19fV0sWzMsMTEsIiIsMCx7ImN1cnZlIjoxLCJzdHlsZSI6eyJoZWFkIjp7Im5hbWUiOiJub25lIn19fV0sWzQsMTEsIiIsMCx7InN0eWxlIjp7ImhlYWQiOnsibmFtZSI6Im5vbmUifX19XSxbNiwxMSwiIiwwLHsic3R5bGUiOnsiYm9keSI6eyJuYW1lIjoiYmFycmVkIn0sImhlYWQiOnsibmFtZSI6Im5vbmUifX19XSxbOCwxMSwiIiwwLHsiY3VydmUiOi0xLCJzdHlsZSI6eyJib2R5Ijp7Im5hbWUiOiJiYXJyZWQifSwiaGVhZCI6eyJuYW1lIjoibm9uZSJ9fX1dLFsxMCwxMSwiIiwyLHsiY3VydmUiOi0yLCJzdHlsZSI6eyJib2R5Ijp7Im5hbWUiOiJiYXJyZWQifSwiaGVhZCI6eyJuYW1lIjoibm9uZSJ9fX1dXQ==
  
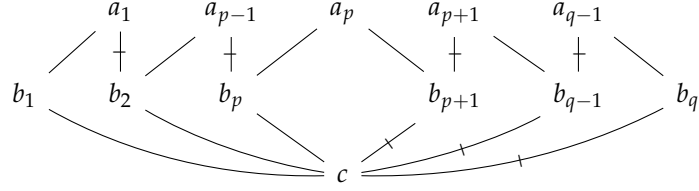
\[\begin{tikzcd}[sep=scriptsize]
    & {a_1} & {a_{p-1}} & {a_p} & {a_{p+1}} & {a_{q-1}} \\
    {b_1} & {b_2} & {b_p} && {b_{p+1}} & {b_{q-1}} & {b_q} \\
    &&& c
    \arrow[no head, from=1-2, to=2-1]
    \arrow["\shortmid"{marking}, no head, from=1-2, to=2-2]
    \arrow[no head, from=1-3, to=2-2]
    \arrow["\shortmid"{marking}, no head, from=1-3, to=2-3]
    \arrow[no head, from=2-3, to=1-4]
    \arrow[no head, from=1-4, to=2-5]
    \arrow["\shortmid"{marking}, no head, from=2-5, to=1-5]
    \arrow[no head, from=1-5, to=2-6]
    \arrow["\shortmid"{marking}, no head, from=2-6, to=1-6]
    \arrow[no head, from=1-6, to=2-7]
    \arrow[curve={height=12pt}, no head, from=2-1, to=3-4]
    \arrow[curve={height=6pt}, no head, from=2-2, to=3-4]
    \arrow[no head, from=2-3, to=3-4]
    \arrow["\shortmid"{marking}, no head, from=2-5, to=3-4]
    \arrow["\shortmid"{marking}, curve={height=-6pt}, no head, from=2-6, to=3-4]
    \arrow["\shortmid"{marking}, curve={height=-12pt}, no head, from=2-7, to=3-4]
  \end{tikzcd}\]
  which lists every chain of the form $c \prec^{\beta} b \prec^{\alpha} a \in \Lambda$ for some $a,\,b$ and $\alpha,\,\beta \in \{+,\,-\}$.
  Note that the left or the right part of the diagram above may be empty, and even the zig-zag itself may be empty.
  For any two loops $c \prec^{\o} b \prec^{-} a_i$ and $c \prec^{\o} b' \prec^{-} a_j$ with $i < j$ and any two loops $b_1 \subseteq b$,\, $b_1' \subseteq b'$, we let $b_1 \blacktriangleleft^1 b_1'$.
\end{defin}

\begin{lem}\label{lem:loop_path}
  Let $c \prec^{\o} b \prec^{+} a \in \Lambda_k$ for some $k$.
  Then we may find a (unique) path as follows in $T_{k+1}$:
  % https://q.uiver.app/#q=WzAsOSxbMCwwLCJhID0gYV8wIl0sWzAsMSwiYiA9IGJfMCJdLFsxLDAsImFfMSJdLFsxLDEsImJfMSJdLFsyLDAsImFfe3AtMSB9Il0sWzIsMSwiYl97cC0xfSJdLFszLDAsImFfcCJdLFszLDEsImJfcCJdLFszLDIsImMiXSxbMCwxLCIiLDAseyJzdHlsZSI6eyJib2R5Ijp7Im5hbWUiOiJiYXJyZWQifSwiaGVhZCI6eyJuYW1lIjoibm9uZSJ9fX1dLFsxLDIsIiIsMCx7InN0eWxlIjp7ImhlYWQiOnsibmFtZSI6Im5vbmUifX19XSxbMiwzLCIiLDAseyJzdHlsZSI6eyJib2R5Ijp7Im5hbWUiOiJiYXJyZWQifSwiaGVhZCI6eyJuYW1lIjoibm9uZSJ9fX1dLFszLDQsIiIsMCx7InN0eWxlIjp7ImJvZHkiOnsibmFtZSI6ImRvdHRlZCJ9LCJoZWFkIjp7Im5hbWUiOiJub25lIn19fV0sWzQsNSwiIiwwLHsic3R5bGUiOnsiYm9keSI6eyJuYW1lIjoiYmFycmVkIn0sImhlYWQiOnsibmFtZSI6Im5vbmUifX19XSxbNiw1LCIiLDIseyJzdHlsZSI6eyJoZWFkIjp7Im5hbWUiOiJub25lIn19fV0sWzYsNywiXFxhbHBoYSIsMCx7InN0eWxlIjp7ImhlYWQiOnsibmFtZSI6Im5vbmUifX19XSxbNyw4LCJcXHZhcmVwc2lsb24iLDAseyJzdHlsZSI6eyJoZWFkIjp7Im5hbWUiOiJub25lIn19fV0sWzUsOCwiIiwwLHsic3R5bGUiOnsiaGVhZCI6eyJuYW1lIjoibm9uZSJ9fX1dLFszLDgsIiIsMCx7ImN1cnZlIjoxLCJzdHlsZSI6eyJoZWFkIjp7Im5hbWUiOiJub25lIn19fV0sWzEsOCwiIiwwLHsiY3VydmUiOjIsInN0eWxlIjp7ImhlYWQiOnsibmFtZSI6Im5vbmUifX19XV0=
  
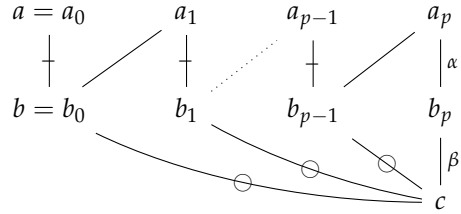
\[\begin{tikzcd}
    {a = a_0} & {a_1} & {a_{p-1 }} & {a_p} \\
    {b = b_0} & {b_1} & {b_{p-1}} & {b_p} \\
    &&& c
    \arrow["\shortmid"{marking}, no head, from=1-1, to=2-1]
    \arrow[no head, from=2-1, to=1-2]
    \arrow["\shortmid"{marking}, no head, from=1-2, to=2-2]
    \arrow[dotted, no head, from=2-2, to=1-3]
    \arrow["\shortmid"{marking}, no head, from=1-3, to=2-3]
    \arrow[no head, from=1-4, to=2-3]
    \arrow["\alpha", no head, from=1-4, to=2-4]
    \arrow["\beta", no head, from=2-4, to=3-4]
    \arrow["\Circle"{marking}, no head, from=2-3, to=3-4]
    \arrow["\Circle"{marking}, curve={height=6pt}, no head, from=2-2, to=3-4]
    \arrow["\Circle"{marking}, curve={height=12pt}, no head, from=2-1, to=3-4]
  \end{tikzcd}\]
  where $p \ge 0,\,\alpha\neq\o$ and ($\beta \neq \o$ or $b_p = \gamma^{(k-1)}\omega$).
\end{lem}
\begin{proof}
  We proceed inductively as follows:
  First, if $b \notin \delta^{-}(C)$, then $b = b_0 = \gamma^{(k-1)}\omega$ by \reflem{lem:DFC_not_a_source_is_iterated_target} and we are done taking $p=0$.
  If $b \in \delta^{-}(C)$, then $b \prec^{-} a_1$ for some unique $a_1 \in \Lambda$ by \reflem{lem:facet_of_Lambda}.
  Then according to \refprop{prop:loop_lozenge}, 
  \begin{itemize}
    \item either $c \prec^{\o} b_1 \prec^{+} a_1$, and we may continue the process, taking $c \prec^{\o} b_1 \prec^{+} a_1$ instead of $c \prec^{\o} b \prec^{-} a$.
    \item or there is a completion $c \prec^{\beta} b_1 \prec^{\alpha} a_1$ for some $\alpha,\,\beta \neq \o$ and we are done taking $p=1$.
  \end{itemize}
  The process ends by \refprop{prop:DFC_strictness_1}.
\end{proof}

\begin{defin}[Loops total order]\label{defin:loops_total_order}
  Fix some $c \in C$, we define an order $\blacktriangleleft^2$ for loops of the form $c \prec^{\o} b \prec^{+} a$ for some $a \in \Lambda \cap N$.
  % That is, we consider a relation $\blacktriangleleft$ defined as containing the restriction of $\triangleleft$ to the loops of the form stated above, and the following new relations:

  Let $c \prec^{\o} b \prec^{+} a$ and $c \prec^{\o} b' \prec^{+} a'$ for some $b \neq b'$ and $a,\, a' \in \Lambda \cap N$. Suppose they are not comparable by $\blacktriangleleft^1$.
  Then by \reflem{lem:loop_path}, we may find some $a'' = a_p = a_q' \in \Lambda$ with the following diagram in $C$, with every $a_i$ and $a_i'$ in $\Lambda$. Note that it is uniquely determined from $a,\,a' \in \Lambda \cap N$.
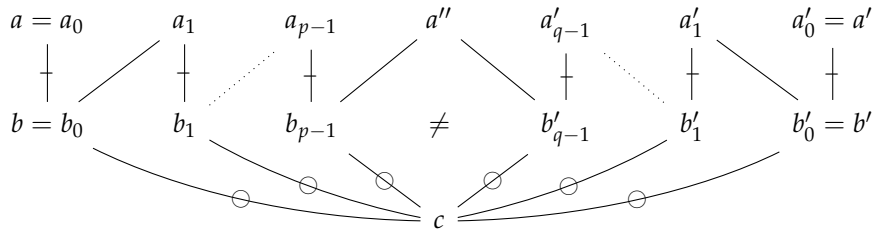
  \[\begin{tikzcd}
    {a = a_0} & {a_1} & {a_{p-1 }} & {a''} & {a_{q-1}'} & {a_1'} & {a_0'=a'} \\
    {b = b_0} & {b_1} & {b_{p-1}} & \neq & {b_{q-1}'} & {b_1'} & {b_0'=b'} \\
    &&& c
    \arrow["\shortmid"{marking}, no head, from=1-1, to=2-1]
    \arrow[no head, from=2-1, to=1-2]
    \arrow["\shortmid"{marking}, no head, from=1-2, to=2-2]
    \arrow[dotted, no head, from=2-2, to=1-3]
    \arrow["\shortmid"{marking}, no head, from=1-3, to=2-3]
    \arrow[no head, from=1-4, to=2-3]
    \arrow["\Circle"{marking}, no head, from=2-3, to=3-4]
    \arrow["\Circle"{marking}, curve={height=6pt}, no head, from=2-2, to=3-4]
    \arrow["\Circle"{marking}, curve={height=12pt}, no head, from=2-1, to=3-4]
    \arrow[no head, from=1-4, to=2-5]
    \arrow["\Circle"{marking}, no head, from=2-5, to=3-4]
    \arrow["\shortmid"{marking}, no head, from=2-5, to=1-5]
    \arrow[dotted, no head, from=1-5, to=2-6]
    \arrow["\shortmid"{marking}, no head, from=2-6, to=1-6]
    \arrow[no head, from=1-6, to=2-7]
    \arrow["\shortmid"{marking}, no head, from=1-7, to=2-7]
    \arrow["\Circle"{marking}, curve={height=6pt}, no head, from=3-4, to=2-6]
    \arrow["\Circle"{marking}, curve={height=12pt}, no head, from=3-4, to=2-7]
  \end{tikzcd}\]
  % and $b_{p-1}$ may be taken to be different from $b_{q-1}'$ because $a_0$ and $a_0'$ have no sources, 
  % neither $b$ nor $b'$ may be included in the other one.
  Then we set $b \blacktriangleleft^2 b'$ if $b_{p-1} \triangleleft b_{q-1}'$.
  In this sense, $\blacktriangleleft^2$ is inherited from the already existing local order $\triangleleft$ (c.f. \refdefin{defin:MOP}), and it is an extension of the latter.
\end{defin}

\begin{lem}\label{lem:loops_total_order_is_a_total_order}
  Fix some $c \in C$ and let $V(c)$ denotes the set of loops $b$ of the form $c \prec^{\o} b \prec^{+} a$ for some $a \in \Lambda \cap N$. Then the relation $\blacktriangleleft \,:=\, \left(\blacktriangleleft^1 \cup \blacktriangleleft^2\right)\vert_{V(c)}$ is a total order.
\end{lem}
\begin{proof}
  $\blacktriangleleft^1$ and $\blacktriangleleft^2$ are irreflexive by definition, hence so is $\blacktriangleleft$.
  
  By definition of $\blacktriangleleft^2$, having $b \blacktriangleleft^1 b'$ and $b' \blacktriangleleft^2 b$ at the same time is impossible, hence in order to show that $\blacktriangleleft$ is antisymmetric, it suffices to show that so are $\blacktriangleleft^1$ and $\blacktriangleleft^2$.
  Since the diagram depicted in \refdefin{defin:loops_total_order} is uniquely determined by $a$ and $a'$ in $\Lambda \cap N$, we have either $b \blacktriangleleft^2 b'$ or $b' \blacktriangleleft^2 b$ depending on whether $b_{p-1} \triangleleft b_{q-1}'$ or $b_{q-1}' \triangleleft b_{p-1}$ but nether both. Hence $\blacktriangleleft^2$ is antisymmetric.
  By definition, $\blacktriangleleft^1$ is also antisymmetric, hence so is $\blacktriangleleft$.

  We now check the transitivity.
  First, suppose $b \blacktriangleleft^2 b' \blacktriangleleft^1 b''$. By definition, $b$ and $b'$ are incomparable by $\blacktriangleleft^1$. Hence $b \blacktriangleleft^1 b''$, whence $b \blacktriangleleft b''$.
  The symmetric case is shown similarly.\\
  Suppose now that $b \blacktriangleleft^1 b'$ and $b' \blacktriangleleft^1 b''$. And suppose they are related by the following two diagrams ($a_p = \hat{a} = a_q'$ and $a_r' = \hat{\hat{a}} = a_s''$):
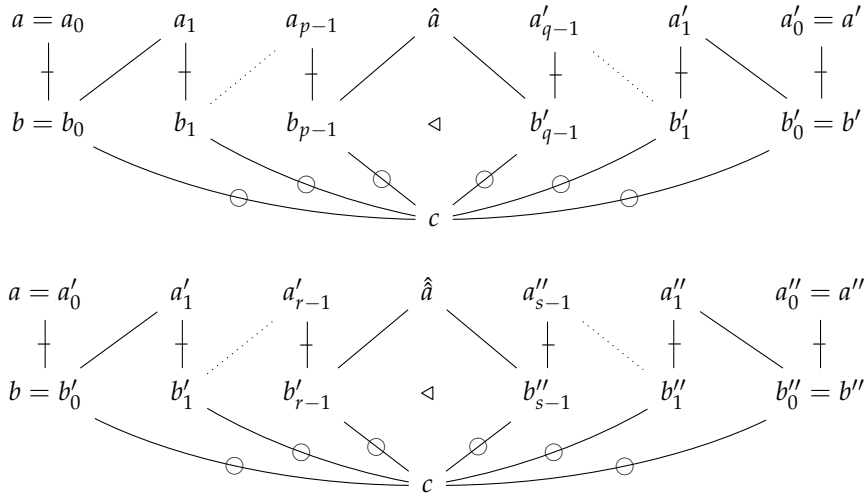
    \[\begin{tikzcd}
      {a = a_0} & {a_1} & {a_{p-1}} & {\hat{a}} & {a_{q-1}'} & {a_1'} & {a_0'=a'} \\
      {b = b_0} & {b_1} & {b_{p-1}} & {\triangleleft} & {b_{q-1}'} & {b_1'} & {b_0'=b'} \\
      &&& c
      \arrow["\shortmid"{marking}, no head, from=1-1, to=2-1]
      \arrow[no head, from=2-1, to=1-2]
      \arrow["\shortmid"{marking}, no head, from=1-2, to=2-2]
      \arrow[dotted, no head, from=2-2, to=1-3]
      \arrow["\shortmid"{marking}, no head, from=1-3, to=2-3]
      \arrow[no head, from=1-4, to=2-3]
      \arrow["\Circle"{marking}, no head, from=2-3, to=3-4]
      \arrow["\Circle"{marking}, curve={height=6pt}, no head, from=2-2, to=3-4]
      \arrow["\Circle"{marking}, curve={height=12pt}, no head, from=2-1, to=3-4]
      \arrow[no head, from=1-4, to=2-5]
      \arrow["\Circle"{marking}, no head, from=2-5, to=3-4]
      \arrow["\shortmid"{marking}, no head, from=2-5, to=1-5]
      \arrow[dotted, no head, from=1-5, to=2-6]
      \arrow["\shortmid"{marking}, no head, from=2-6, to=1-6]
      \arrow[no head, from=1-6, to=2-7]
      \arrow["\shortmid"{marking}, no head, from=1-7, to=2-7]
      \arrow["\Circle"{marking}, curve={height=6pt}, no head, from=3-4, to=2-6]
      \arrow["\Circle"{marking}, curve={height=12pt}, no head, from=3-4, to=2-7]
    \end{tikzcd}\]
  \[\begin{tikzcd}
    {a = a_0'} & {a_1'} & {a_{r-1}'} & {\hat{\hat{a}}} & {a_{s-1}''} & {a_1''} & {a_0''=a''} \\
    {b = b_0'} & {b_1'} & {b_{r-1}'} & {\triangleleft} & {b_{s-1}''} & {b_1''} & {b_0''=b''} \\
    &&& c
    \arrow["\shortmid"{marking}, no head, from=1-1, to=2-1]
    \arrow[no head, from=2-1, to=1-2]
    \arrow["\shortmid"{marking}, no head, from=1-2, to=2-2]
    \arrow[dotted, no head, from=2-2, to=1-3]
    \arrow["\shortmid"{marking}, no head, from=1-3, to=2-3]
    \arrow[no head, from=1-4, to=2-3]
    \arrow["\Circle"{marking}, no head, from=2-3, to=3-4]
    \arrow["\Circle"{marking}, curve={height=6pt}, no head, from=2-2, to=3-4]
    \arrow["\Circle"{marking}, curve={height=12pt}, no head, from=2-1, to=3-4]
    \arrow[no head, from=1-4, to=2-5]
    \arrow["\Circle"{marking}, no head, from=2-5, to=3-4]
    \arrow["\shortmid"{marking}, no head, from=2-5, to=1-5]
    \arrow[dotted, no head, from=1-5, to=2-6]
    \arrow["\shortmid"{marking}, no head, from=2-6, to=1-6]
    \arrow[no head, from=1-6, to=2-7]
    \arrow["\shortmid"{marking}, no head, from=1-7, to=2-7]
    \arrow["\Circle"{marking}, curve={height=6pt}, no head, from=3-4, to=2-6]
    \arrow["\Circle"{marking}, curve={height=12pt}, no head, from=3-4, to=2-7]
  \end{tikzcd}\]
  Suppose that $r \ge q$ (the case $r \le q$ is symmetric). Then we have the following diagram, whence $b \blacktriangleleft b''$:
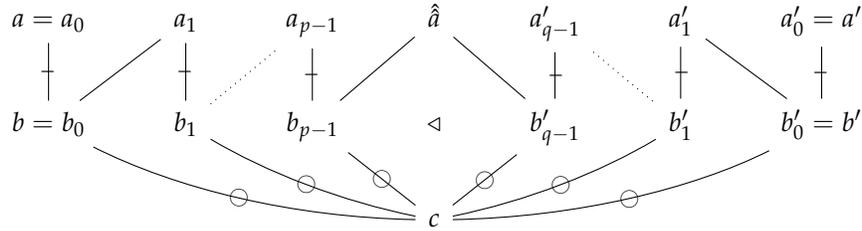
  \[\begin{tikzcd}
    {a = a_0} & {a_1} & {a_{p-1}} & {\hat{\hat{a}}} & {a_{q-1}'} & {a_1'} & {a_0'=a'} \\
    {b = b_0} & {b_1} & {b_{p-1}} & {\triangleleft} & {b_{q-1}'} & {b_1'} & {b_0'=b'} \\
    &&& c
    \arrow["\shortmid"{marking}, no head, from=1-1, to=2-1]
    \arrow[no head, from=2-1, to=1-2]
    \arrow["\shortmid"{marking}, no head, from=1-2, to=2-2]
    \arrow[dotted, no head, from=2-2, to=1-3]
    \arrow["\shortmid"{marking}, no head, from=1-3, to=2-3]
    \arrow[no head, from=1-4, to=2-3]
    \arrow["\Circle"{marking}, no head, from=2-3, to=3-4]
    \arrow["\Circle"{marking}, curve={height=6pt}, no head, from=2-2, to=3-4]
    \arrow["\Circle"{marking}, curve={height=12pt}, no head, from=2-1, to=3-4]
    \arrow[no head, from=1-4, to=2-5]
    \arrow["\Circle"{marking}, no head, from=2-5, to=3-4]
    \arrow["\shortmid"{marking}, no head, from=2-5, to=1-5]
    \arrow[dotted, no head, from=1-5, to=2-6]
    \arrow["\shortmid"{marking}, no head, from=2-6, to=1-6]
    \arrow[no head, from=1-6, to=2-7]
    \arrow["\shortmid"{marking}, no head, from=1-7, to=2-7]
    \arrow["\Circle"{marking}, curve={height=6pt}, no head, from=3-4, to=2-6]
    \arrow["\Circle"{marking}, curve={height=12pt}, no head, from=3-4, to=2-7]
  \end{tikzcd}\]

  Note that $\blacktriangleleft^2$ allows us to compare precisely those elements which are not comparable using $\blacktriangleleft^1$, hence $\blacktriangleleft$ is a total order.
\end{proof}

\begin{lem}[$\gamma(N_{\ge2}) \subseteq \Omega$]\label{lem:nulldot_target_is_a_loop}
  We have the inclusion $\gamma(N_{\ge 2}) \subseteq \Omega$.
\end{lem}
\begin{proof}
  let $y = \gamma(x)$ for some $x \in N_{\ge2}$.
  Then $y \prec^{+} x$ because $\delta(x) = \emptyset$. Suppose $\gamma(y) \prec^{+} y$, then completing the chain $\gamma(y) \prec^{+} y \prec^{+} x$ yields a chain $z \prec^{+} y' \prec^{-} x$, which is impossible because $\delta(x) = \emptyset$.
  Hence $\gamma(y) \prec^{\o} y$, whence $y \in \Omega$.
\end{proof}

\begin{defin}[Nulldot total order]\label{defin:nulldots_order}
  According to \reflem{lem:facet_of_Lambda}, each $y \in \Omega$ admits exactly one $x \in \Lambda$ with $y \prec^{+} x$.
  Hence, by \reflem{lem:nulldot_target_is_a_loop}, the order $\blacktriangleleft$ introduced in \refdefin{defin:loops_total_order} induces a (strict) total order on $W(y) = \{w \in \Lambda \cap N,\, \gamma^2(w) = y\}$.
  We still denote it $\blacktriangleleft$ (there is no confusion because $N \cap \Omega = \emptyset$).\\
  When $b \blacktriangleleft b'$ we say that $b$ \emph{is below} $b'$ or that $b'$ \emph{is above} $b$.
\end{defin}

\begin{defin}[$W_k,\,T_k'$]
  For any $2 \le k < n+2$, we define a subdivision $W_k$ of $T_k$, by $W_k(y) = \{w \in \Lambda \cap N,\, \gamma^2(w) = y\}$.
  As order, we take $\blacktriangleleft$ as introduced in \refdefin{defin:nulldots_order}.
\end{defin}

\begin{lem}[$\msfwdots(T_k') = \msfndots(T_{k+1})$]\label{lem:whitedots_nulldots}
  For $2 \le k < n+2$, we have $\msfwdots(T_k') = \msfndots(T_{k+1})$.
\end{lem}
\begin{proof}
  We have $\msfwdots(T_k') = \bigcup_{z \in C_{k-2}}{W_k(z)} = N \cap \Lambda_k = \msfndots(T_{k+1})$.
\end{proof}

\begin{lem}[Kernel rule]\label{lem:kernel_rule}
  Let $2 \le k < n+2$ and $\sqsubseteq$ the preorder on $T_{k+1}^\bullet \sqcup T_{k+1}^\vert$ given by $x \sqsubseteq y$ iff there is a descending path from $x$ to $y$ in $T_{k+1}$.
  For any $x \in T_{k+1}^\bullet$, $\mathcal{O}_x = \{y \in \msfdots(T_k') \mid y \sqsubseteq x\}$ induces a connected full subgraph of $T_{k}'$.
\end{lem}
\begin{proof}
  Since $T_{n+2}^\bullet = \emptyset$, we have the result for $k = n+1$.

  If $k = n$ or $k<n$ and $\gamma^{(k+1)}\omega$ is a loop, $T_{k+1}$ is a corolla with $a = \gamma^{(k)}\omega$ as only node, whence $\mathcal{O}_{a} = \msfdots(T_k')$ inducing a connected full subgraph of $T_k'$. We suppose now $k<n$ and $\gamma^{(k+1)}\omega\notin\Omega$.

  Let $a \in \Lambda_{k} = T_{k+1}^\bullet$, we prove that $\mathcal{O}_a$ induces a connected full subgraph in $T_k'$ by induction, assuming the result for all nodes $a'$ with a descending path from $a'$ to $a$ in $T_{k+1}$.

  To begin with, we handle the case where $\gamma(a) \in \Omega$. In this case let $c = \gamma^2(a)$. First, we observe that there is no leaf $x \sqsubseteq a$, indeed it would implies having a path $a = a_1 \branch{b = b_1} \cdots \branch{b_{p-1}} a_{p-1} \branch{b_p}$ in $T_{k+1}$, hence a diagram as below by successively applying \refprop{prop:source_confinement}:
  \begin{figure}[H]
    % https://q.uiver.app/#q=WzAsMTEsWzAsMiwiYiA9IGJfMSJdLFsxLDIsImJfMiJdLFsyLDIsImJfe3AtMX0iXSxbMywyLCJiX3twfSJdLFswLDEsImE9YV8xIl0sWzEsMSwiYV8yIl0sWzIsMSwiYV97cC0xfSJdLFszLDEsImFfcCJdLFsyLDAsIlxcZ2FtbWFeeyhrKzIpfVxcb21lZ2EiXSxbNCwxLCJhX3twKzF9Il0sWzIsMywiYyJdLFs0LDAsIiIsMCx7InN0eWxlIjp7ImhlYWQiOnsibmFtZSI6Im5vbmUifX19XSxbNSwxLCIiLDAseyJzdHlsZSI6eyJoZWFkIjp7Im5hbWUiOiJub25lIn19fV0sWzYsMiwiIiwwLHsic3R5bGUiOnsiaGVhZCI6eyJuYW1lIjoibm9uZSJ9fX1dLFs3LDMsIiIsMCx7InN0eWxlIjp7ImhlYWQiOnsibmFtZSI6Im5vbmUifX19XSxbMCwxMCwiIiwwLHsiY3VydmUiOjIsInN0eWxlIjp7ImhlYWQiOnsibmFtZSI6Im5vbmUifX19XSxbMSwxMCwiIiwxLHsic3R5bGUiOnsiaGVhZCI6eyJuYW1lIjoibm9uZSJ9fX1dLFsyLDEwLCIiLDAseyJzdHlsZSI6eyJoZWFkIjp7Im5hbWUiOiJub25lIn19fV0sWzMsMTAsIiIsMCx7InN0eWxlIjp7ImhlYWQiOnsibmFtZSI6Im5vbmUifX19XSxbOCw0LCIiLDAseyJjdXJ2ZSI6Miwic3R5bGUiOnsiaGVhZCI6eyJuYW1lIjoibm9uZSJ9fX1dLFs4LDUsIiIsMSx7InN0eWxlIjp7ImhlYWQiOnsibmFtZSI6Im5vbmUifX19XSxbOCw2LCIiLDAseyJzdHlsZSI6eyJoZWFkIjp7Im5hbWUiOiJub25lIn19fV0sWzgsNywiIiwwLHsic3R5bGUiOnsiaGVhZCI6eyJuYW1lIjoibm9uZSJ9fX1dLFs4LDksIiIsMCx7ImN1cnZlIjotMiwic3R5bGUiOnsiYm9keSI6eyJuYW1lIjoiYmFycmVkIn0sImhlYWQiOnsibmFtZSI6Im5vbmUifX19XSxbNSwwLCIiLDAseyJzdHlsZSI6eyJib2R5Ijp7Im5hbWUiOiJiYXJyZWQifSwiaGVhZCI6eyJuYW1lIjoibm9uZSJ9fX1dLFs2LDEsIiIsMCx7InN0eWxlIjp7ImJvZHkiOnsibmFtZSI6ImJhcnJlZCJ9LCJoZWFkIjp7Im5hbWUiOiJub25lIn19fV0sWzcsMiwiIiwxLHsic3R5bGUiOnsiYm9keSI6eyJuYW1lIjoiYmFycmVkIn0sImhlYWQiOnsibmFtZSI6Im5vbmUifX19XSxbOSwzLCIiLDEseyJzdHlsZSI6eyJoZWFkIjp7Im5hbWUiOiJub25lIn19fV1d
    \[\begin{tikzcd}[row sep = 15pt, column sep = small]
      && {\gamma^{(k+1)}\omega} \\
      {a=a_1} & {a_2} & {a_{p-1}} & {a_p} & {a_{p+1}} \\
      {b = b_1} & {b_2} & {b_{p-1}} & {b_{p}} \\
      && {c}
      \arrow[no head, from=2-1, to=3-1]
      \arrow[no head, from=2-2, to=3-2]
      \arrow[no head, from=2-3, to=3-3]
      \arrow[no head, from=2-4, to=3-4]
      \arrow["\Circle"{marking}, curve={height=12pt}, no head, from=3-1, to=4-3]
      \arrow["\Circle"{marking}, no head, from=3-2, to=4-3]
      \arrow["\Circle"{marking}, no head, from=3-3, to=4-3]
      \arrow["\Circle"{marking}, no head, from=3-4, to=4-3]
      \arrow[curve={height=12pt}, no head, from=1-3, to=2-1]
      \arrow[no head, from=1-3, to=2-2]
      \arrow[no head, from=1-3, to=2-3]
      \arrow[no head, from=1-3, to=2-4]
      \arrow["\shortmid"{marking}, curve={height=-12pt}, no head, from=1-3, to=2-5]
      \arrow["\shortmid"{marking}, no head, from=2-2, to=3-1]
      \arrow[dotted, "\shortmid"{marking}, no head, from=2-3, to=3-2]
      \arrow["\shortmid"{marking}, no head, from=2-4, to=3-3]
      \arrow[no head, from=2-5, to=3-4]
    \end{tikzcd}\]
    \captionof{figure}{A leaf over $a$}\label{fig:loops_over_a}
  \end{figure}
  But this implies $b_p \in \Lambda$, contradicting the loops property.
  By a very similar argument, we see that each nulldot $a' \sqsubseteq a$ satisfies $\gamma^2(a') = c$, whence $a'$ lying on the edge $c$ of $T_k'$.
  Moreover, they consist in all the whitedots contained in $a$, which must appear consecutively on the edge $c$, whence the connectivity.
  
  We suppose now that $\gamma(a) \notin \Omega$. For any $b \prec^{-} a$, we let $\mathcal{O}_b := \{x \in \msfdots(T_k) \mid x \sqsubseteq b\}$. By induction hypothesis, $\mathcal{O}_b$ induces a connected full subgraph of $T_k'$ for any $b \prec^{-} a$, and we have $\mathcal{O}_a = \bigcup_{b\prec^{-}a}\mathcal{O}_b$.
  We prove the following properties.
  \begin{enumerate}
    \item Let $b \branch{c} b'$ in $\Delta_a$, then $\mathcal{O}_b$ and $\mathcal{O}_{b'}$ are linked by an edge.
    \item Let $b \in \delta(a) \cap \Omega$, then $\mathcal{O}_b \subseteq \msfwdots(T_k')$.
    \item Let $b \branch{c}$ in $T_k'$. If $b \in \mathcal{O}_a$ then $\mathcal{O}_a \cap W_k(c)$ is a lower set of $W_k(c)$.
    \item Let $\branch{c} c$ in $T_k'$. If $b \in \mathcal{O}_a$ then $\mathcal{O}_a \cap W_k(c)$ is an upper set of $W_k(c)$. 
  \end{enumerate}
  By using the tree structure of $\Delta_a$, 1. and 2. imply that $\mathcal{O}_a \cap \msfbdots(T_k')$ induces a connected full subgraph of $T_k$. Then 3. and 4. imply that $\mathcal{O}_a$ induces a connected full subgraph of $T_k'$.\\
  We are left to prove those few properties.
  \begin{enumerate}
    \item Suppose $b \branch{c} b'$ in $\Delta_a$. Then there is a lozenge $(c \prec^{-} b \prec^{-} a,\, c \prec^{+} b' \prec^{-} a)$ in $C$. By successively completing lozenges, we may form a diagram as below in $C$:
      \[\begin{tikzcd}[row sep = scriptsize, column sep = tiny]
        &&& {\gamma^{(k+1)}\omega} \\
        {\gamma^{(k)}\omega} & {a_0} & {a_{p-1}} & a & {a_q'} & {a_0'} & {\gamma^{(k)}\omega} \\
        {b_0} & {b_2} & {b = b_p} && {b' = b_q'} & {b_1'} & {b_0'} \\
        &&& c
        \arrow[no head, from=2-2, to=3-1]
        \arrow["\shortmid"{marking}, no head, from=2-2, to=3-2]
        \arrow[dotted, no head, from=2-3, to=3-2]
        \arrow["\shortmid"{marking}, no head, from=2-3, to=3-3]
        \arrow[no head, from=3-3, to=2-4]
        \arrow[no head, from=2-4, to=3-5]
        \arrow["\shortmid"{marking}, no head, from=3-5, to=2-5]
        \arrow[dotted, no head, from=2-5, to=3-6]
        \arrow["\shortmid"{marking}, no head, from=3-6, to=2-6]
        \arrow[no head, from=2-6, to=3-7]
        \arrow[curve={height=12pt}, no head, from=3-1, to=4-4]
        \arrow[curve={height=6pt}, no head, from=3-2, to=4-4]
        \arrow[no head, from=3-3, to=4-4]
        \arrow["\shortmid"{marking}, no head, from=3-5, to=4-4]
        \arrow["\shortmid"{marking}, curve={height=-6pt}, no head, from=3-6, to=4-4]
        \arrow["\shortmid"{marking}, curve={height=-12pt}, no head, from=3-7, to=4-4]
        \arrow[no head, from=1-4, to=2-4]
        \arrow[no head, from=1-4, to=2-3]
        \arrow[curve={height=6pt}, no head, from=1-4, to=2-2]
        \arrow[no head, from=1-4, to=2-5]
        \arrow[curve={height=-6pt}, no head, from=1-4, to=2-6]
        \arrow["\shortmid"{marking}, curve={height=-12pt}, no head, from=1-4, to=2-7]
        \arrow["\shortmid"{marking}, curve={height=12pt}, no head, from=1-4, to=2-1]
        \arrow[no head, from=2-1, to=3-1]
        \arrow[no head, from=2-7, to=3-7]
      \end{tikzcd}\]
      hence finding a relation $b_0 \branch{c} b_0'$ in $T_k$.
      Since $b_0 \in \mathcal{O}_b$ and $b_0' \in \mathcal{O}_{b'}$, it yields the result.
  
    \item Let $b \in \delta(a) \cap \Omega$. We may prove that there is no leaf $x \sqsubseteq b$ by the same argument we used when dealing with the case $\gamma(a) \in \Omega$.

    \item Let $b \branch{c}$ in $T_k'$ with $b \in \mathcal{O}_a$.
      Notice that $\gamma^{(k-1)}\omega \notin \Omega$, because it would imply that $b$ is a loop by \refprop{prop:source_confinement}, which is not the case since $b \in \Lambda$.
      Hence, for $c \prec^{\o} b' \prec^{+} a' \in \Lambda$, using repeatedly oriented thinness and \refprop{prop:loop_lozenge} yields the following diagram in $C$, where $\alpha,\,\beta \neq \o$.
      % https://q.uiver.app/#q=WzAsMTAsWzMsMCwiXFxnYW1tYV57KGsrMSl9XFxvbWVnYSJdLFszLDEsImEnIl0sWzMsMiwiYiciXSxbMywzLCJjIl0sWzIsMSwiXFxidWxsZXQiXSxbMiwyLCJcXGJ1bGxldCJdLFsxLDEsIlxcYnVsbGV0Il0sWzEsMiwiXFxidWxsZXQiXSxbMCwyLCJiJyciXSxbMCwxLCJhJyciXSxbMiwzLCIiLDAseyJzdHlsZSI6eyJoZWFkIjp7Im5hbWUiOiJub25lIn19fV0sWzEsMiwiIiwwLHsic3R5bGUiOnsiYm9keSI6eyJuYW1lIjoiYmFycmVkIn0sImhlYWQiOnsibmFtZSI6Im5vbmUifX19XSxbMCwxLCIiLDAseyJzdHlsZSI6eyJoZWFkIjp7Im5hbWUiOiJub25lIn19fV0sWzQsMiwiIiwwLHsic3R5bGUiOnsiaGVhZCI6eyJuYW1lIjoibm9uZSJ9fX1dLFs0LDUsIiIsMix7InN0eWxlIjp7ImJvZHkiOnsibmFtZSI6ImJhcnJlZCJ9LCJoZWFkIjp7Im5hbWUiOiJub25lIn19fV0sWzYsNSwiIiwwLHsic3R5bGUiOnsiYm9keSI6eyJuYW1lIjoiZG90dGVkIn0sImhlYWQiOnsibmFtZSI6Im5vbmUifX19XSxbNiw3LCIiLDIseyJzdHlsZSI6eyJib2R5Ijp7Im5hbWUiOiJiYXJyZWQifSwiaGVhZCI6eyJuYW1lIjoibm9uZSJ9fX1dLFs1LDMsIiIsMix7InN0eWxlIjp7ImhlYWQiOnsibmFtZSI6Im5vbmUifX19XSxbNywzLCIiLDIseyJjdXJ2ZSI6MSwic3R5bGUiOnsiaGVhZCI6eyJuYW1lIjoibm9uZSJ9fX1dLFs4LDMsIlxcYmV0YSIsMix7ImN1cnZlIjoyLCJzdHlsZSI6eyJoZWFkIjp7Im5hbWUiOiJub25lIn19fV0sWzksNywiIiwyLHsic3R5bGUiOnsiaGVhZCI6eyJuYW1lIjoibm9uZSJ9fX1dLFs5LDgsIlxcYWxwaGEiLDIseyJzdHlsZSI6eyJoZWFkIjp7Im5hbWUiOiJub25lIn19fV0sWzAsOSwiIiwxLHsiY3VydmUiOjIsInN0eWxlIjp7ImhlYWQiOnsibmFtZSI6Im5vbmUifX19XSxbMCw2LCIiLDEseyJjdXJ2ZSI6MSwic3R5bGUiOnsiaGVhZCI6eyJuYW1lIjoibm9uZSJ9fX1dLFswLDQsIiIsMSx7InN0eWxlIjp7ImhlYWQiOnsibmFtZSI6Im5vbmUifX19XV0=
      \[\begin{tikzcd}[sep=scriptsize]
        &&& {\gamma^{(k+1)}\omega} \\
        {a''} & \bullet & \bullet & {a'} \\
        {b''} & \bullet & \bullet & {b'} \\
        &&& c
        \arrow["\Circle"{marking}, no head, from=3-4, to=4-4]
        \arrow["\shortmid"{marking}, no head, from=2-4, to=3-4]
        \arrow[no head, from=1-4, to=2-4]
        \arrow[no head, from=2-3, to=3-4]
        \arrow["\shortmid"{marking}, no head, from=2-3, to=3-3]
        \arrow[dotted, no head, from=2-2, to=3-3]
        \arrow["\shortmid"{marking}, no head, from=2-2, to=3-2]
        \arrow["\Circle"{marking}, no head, from=3-3, to=4-4]
        \arrow["\Circle"{marking}, curve={height=6pt}, no head, from=3-2, to=4-4]
        \arrow["\beta"', curve={height=12pt}, no head, from=3-1, to=4-4]
        \arrow[no head, from=2-1, to=3-2]
        \arrow["\alpha"', no head, from=2-1, to=3-1]
        \arrow[curve={height=12pt}, no head, from=1-4, to=2-1]
        \arrow[curve={height=6pt}, no head, from=1-4, to=2-2]
        \arrow[no head, from=1-4, to=2-3]
      \end{tikzcd}\]
      Then, for some $a_1$ such that $b \prec^{-} a_1 \in \Lambda$, applying \refprop{prop:DFC_hexagon} to the pair $$(c \prec^{-} b \prec^{-} a_1 \prec^{-} \gamma^{(k+1)}\omega,\,\,\, c \prec^{\beta} b'' \prec^{\alpha} a'' \prec^{-} \gamma^{(k+1)}\omega)$$ yields a diagram as follows in $C$, with $p < q \le r$.
      \[\begin{tikzcd}[sep=scriptsize]
        &&&& {\gamma^{(k+1)}\omega} \\
        \\
        {\gamma^{(k)}\omega} & {a_1} & {a_p} & {a_{p+1}} & {a_{p+2}} & {a_q} & {a_{q+1}} & {a_r = a'} \\
        {b = b_0} & {b_1} & {b_p} && {b_{p+1}} & {b_{q-1}} & {b_q} & {b_r = b'} \\
        \\
        &&&& c
        \arrow["\shortmid"{marking}, no head, from=4-5, to=6-5]
        \arrow[curve={height=6pt}, no head, from=1-5, to=3-4]
        \arrow[curve={height=8pt}, no head, from=4-3, to=6-5]
        \arrow[curve={height=14pt}, no head, from=4-2, to=6-5]
        \arrow[curve={height=20pt}, no head, from=4-1, to=6-5]
        \arrow["\shortmid"{marking}, no head, from=3-1, to=4-1]
        \arrow["\shortmid"{marking}, curve={height=20pt}, no head, from=1-5, to=3-1]
        \arrow[curve={height=14pt}, no head, from=1-5, to=3-2]
        \arrow[curve={height=8pt}, no head, from=1-5, to=3-3]
        \arrow[no head, from=4-1, to=3-2]
        \arrow["\shortmid"{marking}, no head, from=3-2, to=4-2]
        \arrow[dotted, no head, from=3-3, to=4-2]
        \arrow["\shortmid"{marking}, no head, from=3-3, to=4-3]
        \arrow[no head, from=4-3, to=3-4]
        \arrow[no head, from=3-4, to=4-5]
        \arrow["\shortmid"{marking}, no head, from=4-5, to=3-5]
        \arrow[dotted, no head, from=3-5, to=4-6]
        \arrow["\shortmid"{marking}, no head, from=3-6, to=4-6]
        \arrow[no head, from=1-5, to=3-5]
        \arrow[curve={height=-6pt}, no head, from=1-5, to=3-6]
        \arrow["\shortmid"{marking}, curve={height=-6pt}, no head, from=4-6, to=6-5]
        \arrow[no head, from=3-6, to=4-7]
        \arrow["\shortmid"{marking}, no head, from=3-7, to=4-7]
        \arrow[dotted, no head, from=3-7, to=4-8]
        \arrow["\shortmid"{marking}, no head, from=3-8, to=4-8]
        \arrow[curve={height=-20pt}, no head, from=1-5, to=3-8]
        \arrow[curve={height=-14pt}, no head, from=1-5, to=3-7]
        \arrow["\Circle"{marking}, curve={height=-8pt}, no head, from=4-7, to=6-5]
        \arrow["\Circle"{marking}, curve={height=-14pt}, no head, from=4-8, to=6-5]
      \end{tikzcd}\]
      Then $a_{p+1} \sqsubseteq a$.
      If $c \prec^{\o} b'' \prec^{+} a'' \in \Lambda$ for some $b'' \blacktriangleleft b'$ in $W_k(C)$, then similarly we may find a zig-zag
      
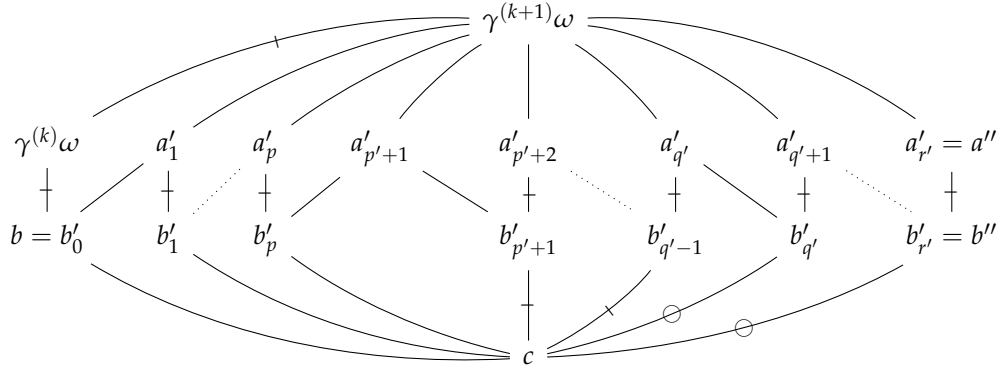
\[\begin{tikzcd}[sep=scriptsize]
        &&&& {\gamma^{(k+1)}\omega} \\
        \\
        {\gamma^{(k)}\omega} & {a_1'} & {a_p'} & {a'_{p'+1}} & {a'_{p'+2}} & {a'_{q'}} & {a'_{q'+1}} & {a'_{r'} = a''} \\
        {b = b'_0} & {b'_1} & {b'_p} && {b'_{p'+1}} & {b'_{q'-1}} & {b'_{q'}} & {b'_{r'} = b''} \\
        \\
        &&&& c
        \arrow["\shortmid"{marking}, no head, from=4-5, to=6-5]
        \arrow[curve={height=6pt}, no head, from=1-5, to=3-4]
        \arrow[curve={height=8pt}, no head, from=4-3, to=6-5]
        \arrow[curve={height=14pt}, no head, from=4-2, to=6-5]
        \arrow[curve={height=20pt}, no head, from=4-1, to=6-5]
        \arrow["\shortmid"{marking}, no head, from=3-1, to=4-1]
        \arrow["\shortmid"{marking}, curve={height=20pt}, no head, from=1-5, to=3-1]
        \arrow[curve={height=14pt}, no head, from=1-5, to=3-2]
        \arrow[curve={height=8pt}, no head, from=1-5, to=3-3]
        \arrow[no head, from=4-1, to=3-2]
        \arrow["\shortmid"{marking}, no head, from=3-2, to=4-2]
        \arrow[dotted, no head, from=3-3, to=4-2]
        \arrow["\shortmid"{marking}, no head, from=3-3, to=4-3]
        \arrow[no head, from=4-3, to=3-4]
        \arrow[no head, from=3-4, to=4-5]
        \arrow["\shortmid"{marking}, no head, from=4-5, to=3-5]
        \arrow[dotted, no head, from=3-5, to=4-6]
        \arrow["\shortmid"{marking}, no head, from=3-6, to=4-6]
        \arrow[no head, from=1-5, to=3-5]
        \arrow[curve={height=-6pt}, no head, from=1-5, to=3-6]
        \arrow["\shortmid"{marking}, curve={height=-6pt}, no head, from=4-6, to=6-5]
        \arrow[no head, from=3-6, to=4-7]
        \arrow["\shortmid"{marking}, no head, from=3-7, to=4-7]
        \arrow[dotted, no head, from=3-7, to=4-8]
        \arrow["\shortmid"{marking}, no head, from=3-8, to=4-8]
        \arrow[curve={height=-20pt}, no head, from=1-5, to=3-8]
        \arrow[curve={height=-14pt}, no head, from=1-5, to=3-7]
        \arrow["\Circle"{marking}, curve={height=-8pt}, no head, from=4-7, to=6-5]
        \arrow["\Circle"{marking}, curve={height=-14pt}, no head, from=4-8, to=6-5]
      \end{tikzcd}\]
      Then by oriented thinness and the uniqueness of lozenge completions, both zig-zags must coincide up to $a_{\mathsf{min}(q,\,q')}$, hence at least up to $a_{\mathsf{min}(p+1,\,p'+1)}$.
      However, if $p' \le p$, then $a'_{p'+1} = a_{p'+1} \sqsubseteq a_{p+1} \sqsubseteq a$. Whence $a'' \sqsubseteq a$, that is $a'' \in \mathcal{O}_a$, which shows that $\mathcal{O}_a \cap W_k(c)$ is a lower set in $W_k(c)$.

      \item The case $c \branch{b}$ is symmetric to the previous one.
    \qedhere
  \end{enumerate}
\end{proof}

\begin{lem}[$T_2$ is linear]\label{lem:tree_2}
  The tree $T_2$ is linear.
\end{lem}
\begin{proof}
  Let $x\in\Lambda_1 = T_2^\bullet$. By \reflem{lem:nulldot_target_is_a_loop}, $x$ is not a nulldot because $\ast = \gamma^2(x) \nprec^{\o} \gamma(x)$.
  Hence $x$ admits at least one source. By oriented thinness this source is unique, hence $T_2$ is linear.
\end{proof}

\begin{defin}[$K_0 : T_0 \to T_1,\, K_1 : T_1 \to T_2$]\label{defin:augmentation_map}
  Since $T_2$ is a linear tree (\reflem{lem:tree_2}), there is a unique exact zoom complex (up to renamings)
  $$T_0 \overset{K_0}{\to} T_1 \overset{K_1}{\to} T_2$$
  We define the constellations $K_0$ and $K_1$ as such.
\end{defin}

\begin{defin}\label{defin:DFC_to_Zoom_object}
  For $i>0$, Let $K_i$ denote the exact constellation $T_i \to T_{i+1}$. It is well defined because of Lemmas \ref{lem:blackdots_leaves},\,\ref{lem:whitedots_nulldots} and \ref{lem:kernel_rule}.
  Then $$T_0 \overset{K_0}{\longrightarrow} T_1 \overset{K_1}{\longrightarrow} \cdots \overset{K_{n-1}}{\longrightarrow} T_n$$ is an opetope of dimension $n$. It is denoted $Z(C)$.
\end{defin}

\begin{defin}\label{defin:DFC_to_Zoom_morphism}
  Let $f : C \to D$ be an isomorphism of DFC. There is an isomorphism $Z(f) : Z(C) \to Z(D)$ of opetopes, as follows.
  Let $Z(D) = \left(U_0 \overset{L_0}{\longrightarrow} U_1 \overset{L_1}{\longrightarrow} \cdots \overset{L_{n-1}}{\longrightarrow} U_n\right)$ and $(U_k')_{0\le k \le n}$ the associated subdivided trees. Then $Z(f)_k : T_k' \to U_k'$ is defined as $f$ on blackdots, whitedots and leaves.
\end{defin}
\begin{proof}
  We check that $Z(f)$ defined as above is an isomorphism of opetopes.\\
  Let $k \in \intInter{1}{n+2}$.
  First, $f$ sends bijectively $\Lambda_{k-1} \subseteq C$ to $\Lambda_{k-1} \subseteq D$, hence bijectively $T_k^\bullet$ to $S_k^\bullet$.
  It also sends bijectively $C_{k-2}$ to $D_{k-2}$, hence bijectively $T_k^{\edge}$ to $U_k^{\edge}$.\\
  Since $(y \prec^{\alpha} x) \leftrightarrow (f(y) \prec^{\alpha} f(x))$ for any $\alpha$, $Z(f)_k$ preserves relations $x \branch{y}$ or $y \branch{x}$ from $T_k$ to $U_k$.
  Hence $f$ induces an isomorphism of rooted trees $T_k \simeq U_k$.

  We then check that it extends as an isomorphism $T_k' \simeq U_k'$. For $y \in T_k^{\edge}$, 
  $$\begin{array}{rcl}
    W_k(y) & = & \{w \in \Lambda_{k} \cap N \mid \gamma^2(w) = y\} \\
    ~ & = & \{w \in C_k \mid (\delta(w) = \emptyset) \wedge (\nexists v,\, w \prec^{+} v) \wedge \gamma^2(w) = y \} \\
    ~ & \overset{f}{\simeq} & \{w \in D_k \mid (\delta(w) = \emptyset) \wedge (\nexists v,\, w \prec^{+} v) \wedge \gamma^2(w) = f(y)\} \\
    ~ & =: & X_k(f(y))
  \end{array}$$ is the set of whitedots on $f(y)$ in $U_k'$. Since $f$ is an isomorphism, it sends diagrams of the form given in \refdefin{defin:loops_total_order} to diagrams of the same form in $D$. It also sends sourceless cells to sourceless cells, and by assumption it preserves the order $\blacktriangleleft$ of \refdefin{defin:MOP}.
  Hence it preserves $\blacktriangleleft$ as given by \refdefin{defin:loops_total_order} and \refdefin{defin:nulldots_order}.
  Hence $Z(f)_k$ is a well-defined isomorphism of subdivided trees.

  Since the structural maps in $K_i$ and $L_i \,(0<i<n+2)$ are identities, $Z(f)$ is compatible with them.
  Hence $Z(f)$ is an isomorphism of opetopes.
\end{proof}

\begin{thm}[$Z : \DFC^\core \to \KJBM$]\label{thm:DFC_to_Ope}
  There is a functor $Z : \DFC^\core \to \KJBM$ given by \refdefin{defin:DFC_to_Zoom_object} on objects and by \refdefin{defin:DFC_to_Zoom_morphism} on isomorphisms.
\end{thm}
\begin{proof}
  The functoriality is clear from \refdefin{defin:DFC_to_Zoom_morphism}.
\end{proof}

\section{From opetopes to dendritic face complexes}
\label{sec:Zoom_to_DFC}

\para In this section, our aim is to define a functor $P : \KJBM \to \DFC$.\\
In order to do so, it will be enlightening to extends canonically an opetope $Y = \left(S_0 \to S_1 \to \cdots \to S_n\right)$ of dimension $n$ to an opetope $\left(S_0 \to S_1 \to \cdots \to S_{n+2}\right)$.
Indeed, in \refsec{sec:DFC_to_Zoom}, any $n$-dimensional DFC $C$ gave rise to a sequence of trees $T_k$ for $0 \le k \le n+2$, the two last of them being retrieved by this canonical extension.
Formally, it goes as follows.

\begin{defin}\label{defin:zoom_extension}
  Consider an opetope $Y = \left(S_0 \to S_1 \to \cdots \to S_n\right)$ of dimension $n$.
  We let
  \begin{itemize}
    \item $S_{n+1}$ be the corolla whose node is $\{\top\}$ and set of leaves is $S_n^\bullet$,
    \item $S_{n+2}$ be a unit tree whose root edge is $\top$,
    \item $V_n$ (resp. $V_{n+1}$) be the empty subdivision on $S_{n}$ (resp. $S_{n+1}$),
    \item $S_k' = (S_k,\,V_k)$ the subdivided tree associated to $S_k$ for any $k<n+2$.
  \end{itemize}
  Note that there are two exact constellations $S_n \to S_{n+1} \to S_{n+2}$.
\end{defin}

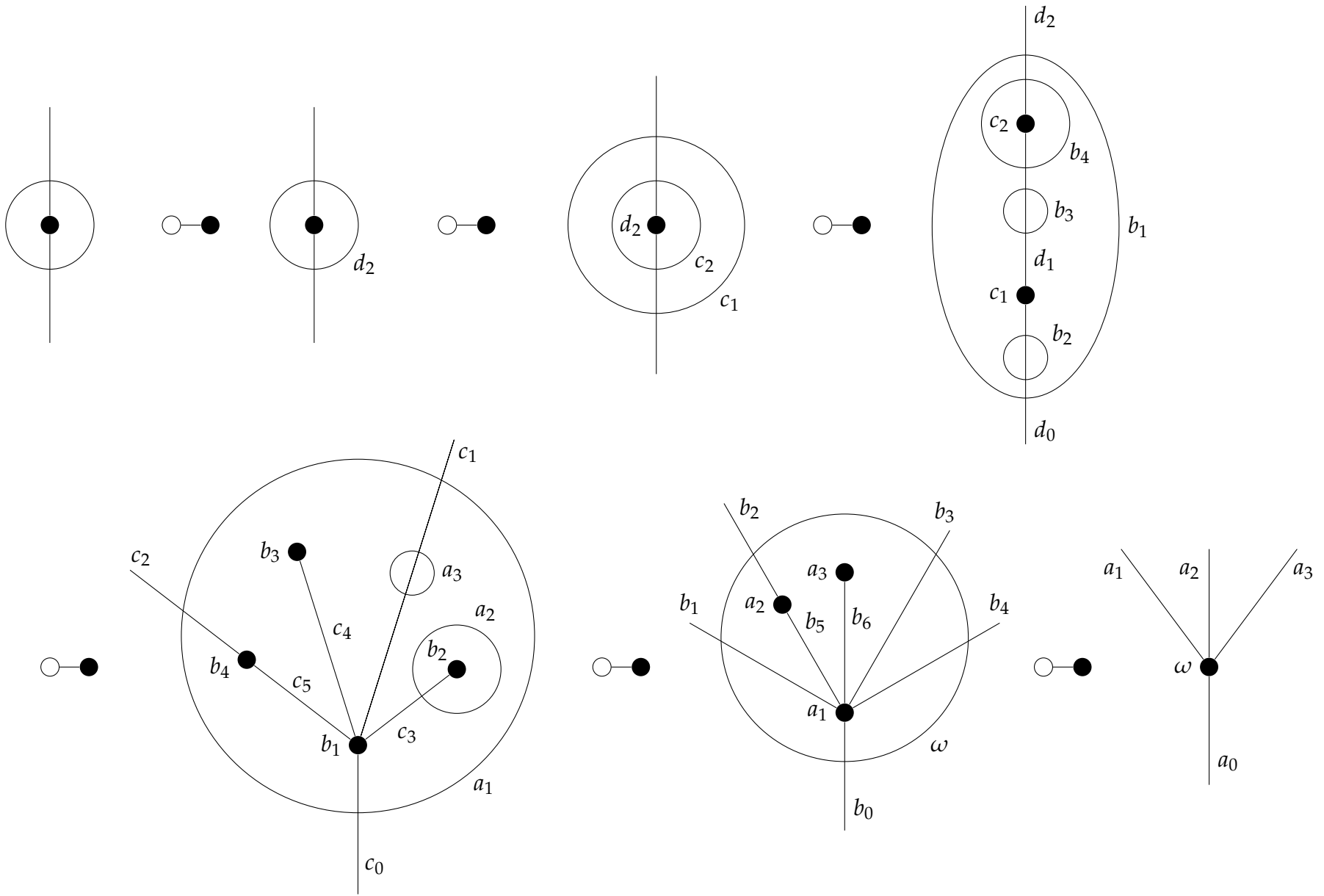
\begin{sidewaysfigure}[p]
  % \documentclass[landscape]{article}
% \usepackage[utf8]{inputenc}
% \usepackage{graphics, graphicx, tikz, geometry}
% \usepackage{amsfonts}
% \usepackage{amssymb,amsmath,amsthm,stmaryrd,mathrsfs,wasysym}
% %
% %%% Set the fonts
% \usepackage{mathpazo}
% \usepackage[scaled=0.95]{helvet}
% \usepackage{courier}
% \usepackage{adjustbox}
% \geometry{hmargin=0.5cm,vmargin=1cm}
% \usetikzlibrary{arrows, arrows.meta}
% \linespread{1.05} % Palatino looks better with this
% \definecolor{mypurple}{rgb}{0.7,0,0.7}
% \usetikzlibrary{calc, trees}
% %
% \begin{document}
% %
  \adjustbox{scale=1, center}{

    \begin{tikzpicture}
      [ blackdot/.style = {circle, fill}
      , whitedot/.style = {circle, draw}
      ]

      % constellation 0
      \node (X0) at (0,0) {
        \begin{tikzpicture}
          [ level distance = 80px
          , bdot/.style = {circle, fill}
          , boundary/.style = {draw = none}
          , scale = 0.8
          ] 
          \node at (0,0) {} [grow' = up, boundary]
            child {
              node (mid_0) [bdot] {}
                child{
                  node {} [boundary]
                }
            };
          \draw (mid_0) circle (1);
        \end{tikzpicture}
      };

      % zoom
      \node at (2.2,0) [whitedot] {}
        child[level distance = 20px, grow = right]{
          node [blackdot] {}
        };

      % constellation 1
      \node (X1) at (5,0) {
        \begin{tikzpicture}
          [ level distance = 80px
          , bdot/.style = {circle, fill}
          , boundary/.style = {draw = none}
          , scale = 0.8
          ] 
          \node {} [grow' = up, boundary]
            child {
              node (mid_1) [bdot] {}
                child{
                  node {} [boundary]
                }
            };
          \draw (mid_1) circle (1)
            node [xshift = 26, yshift = -20] {\large$d_2$};
        \end{tikzpicture}
      };

      % zoom
      \node at (7.2,0) [whitedot] {}
        child[level distance = 20px, grow = right]{
          node [blackdot] {}
        };

      % constellation 2
      \node (X2) at (11,0) {
        \begin{tikzpicture}
          [ level distance = 80px
          , bdot/.style = {circle, fill}
          , boundary/.style = {draw = none}
          , scale = 0.8
          ]
          \node {} [grow' = up, boundary, level distance = 100px]
            child {
              node (d2) [bdot, label=left:{\large$d_2$\hspace{-1px}~}] {}
                child{
                  node {} [boundary]
                }
            };
          \draw (d2) circle (1)
            node [xshift = 25, yshift = -20] {\large$c_2$};
          \draw (d2) circle (2)
            node [xshift = 38, yshift = -40] {\large$c_1$};
        \end{tikzpicture}
      };

      % zoom
      \node at (14,0) [whitedot] {}
        child[level distance = 20px, grow = right]{
          node [blackdot] {}
        };

      % constellation 3
      \node (X3) at (18,0) {
        \begin{tikzpicture}
          [ level distance = 80px
          , bdot/.style = {circle, fill}
          , boundary/.style = {draw = none}
          , scale = 0.8
          ] 
          \node (root_3) {} [grow' = up, boundary, level distance = 100px]
            child {
              node (c1) [bdot, label=left:{\large$c_1$}] {}
                child[level distance = 110px]{
                  node (c2) [bdot, label=left:{\large$c_2$}] {} [level distance = 80px]
                    child{
                      node {} [boundary]
                      edge from parent node [right = 1, pos = 0.9] {\large$d_2$}
                    }
                }
              edge from parent node [right = 1, pos = 0.1] {\large$d_0$}
            };
          \node (mid_3) at ($(c1)!0.4!(c2)$) {};
          \draw (mid_3) ellipse (60px and 110px) 
            node [xshift = 58] {\large$b_1$};
          \draw (mid_3) [yshift = 10] circle (1/2)
            node [xshift = 20] {\large$b_3$};
          \node at (mid_3) [yshift = -17, xshift = 10] {\large$d_1$};
          \draw (c2) circle (1)
            node [xshift = 28, yshift = -15] {\large$b_4$};
          \draw ($(root_3)!0.6!(c1)$) circle (1/2)
            node [xshift = 19, yshift = 12] {\large$b_2$};
        \end{tikzpicture}
      };

      % zoom
      \node at (0,-8) [whitedot] {}
        child[level distance = 20px, grow = right]{
          node [blackdot] {}
        };

      % constellation 4
      \node (X4) at (5,-8) {
        \begin{tikzpicture}
          [ level distance = 80px
          , bdot/.style = {circle, fill}
          , boundary/.style = {draw = none}
          , scale = 0.8
          ]
          \pgfmathsetmacro{\theta}{35}
          \node (root_3) {} [grow' = up, boundary, level distance = 100px]
            child {
              node (b1) [bdot, label=left:{\large$b_1$}] {} [clockwise from = 90+1.5*\theta, sibling angle = \theta]
                child[level distance = 90px]{
                  node (b4) [bdot, label={left, yshift = -4}:{\large$b_4$}] {} [level distance = 100px]
                    child{
                      node [boundary] {}
                      edge from parent node [pos = 0.9, above = 3] {\large$c_2$}
                    }
                  edge from parent node [above = 2] {\large$c_5$}
                }
                child[level distance = 130px]{
                  node (b3) [bdot, label=left:{\large$b_3$}] {}
                  edge from parent node [right = 2, pos = 0.6] {\large$c_4$}
                }
                child[level distance = 210px]{
                  node [boundary] {}
                  edge from parent node [pos = 0.55] (a3) {}
                  edge from parent node [pos = 0.95, right = 1] {\large$c_1$}
                }
                child[level distance = 80px]{
                  node (b2) [bdot] {}
                  edge from parent node [below right, pos = 0.3] {\large$c_3$}
                }
              edge from parent node [right, pos = 0.2] {\large$c_0$}
            };
          \draw ($1.7*(b1)$) circle (4)
            node [right = 65, below = 70] {\large$a_1$};
          \draw (a3) circle (1/2)
            node [right = 12] {\large$a_3$};
          \draw (b2) circle (1)
            node [above = 29, right = 6] {\large$a_2$};
          \node at (b2) [xshift = -10, yshift = 10] {\large$b_2$};
        \end{tikzpicture}
      };

      % zoom
      \node at (10,-8) [whitedot] {}
        child[level distance = 20px, grow = right]{
          node [blackdot] {}
        };

      % constellation 5
      \node (X5) at (14.4,-8) {
        \begin{tikzpicture}
          [ level distance = 80px
          , bdot/.style = {circle, fill}
          , boundary/.style = {draw = none}
          , scale = 0.8
          ]
          \pgfmathsetmacro{\theta}{30}
          \node (root_3) {} [grow' = up, boundary, level distance = 80px]
            child {
              node (a1) [bdot, label={left}:{\large$a_1$}] {} [clockwise from = 90+2*\theta, sibling angle = \theta]
                child [level distance = 120px] {
                  node [boundary] {}
                  edge from parent node [above = 1, pos = 1] {\large$b_1$}
                }
                child{
                  node (a2) [bdot, label={left}:{\large$a_2$}] {}
                    child [clockwise from = 90+\theta] {
                      node [boundary] {}
                      edge from parent node [right = 3, pos = 0.95] {\large$b_2$}
                    }
                  edge from parent node [right = 3, pos = 0.9] {\large$b_5$}
                }
                child [level distance = 90px] {
                  node (a3) [bdot, label={left}:{\large$a_3$}] {}
                  edge from parent node [right, pos = 0.7] {\large$b_6$}
                }
                child [level distance = 140px] {
                  node [boundary] {}
                  edge from parent node [above = 5, pos = 0.95] {\large$b_3$}
                }
                child [level distance = 120px] {
                  node [boundary] {}
                  edge from parent node [above = 1, pos = 1] {\large$b_4$}
                }
              edge from parent node [right = 1, pos = 0.2] {\large$b_0$}
            };
          \draw ($1.6*(a1)$) circle (2.8)
            node [right = 49, below = 49] {\large$\omega$};
        \end{tikzpicture}
      };

      % zoom
      \node at (18,-8) [whitedot] {}
        child[level distance = 20px, grow = right]{
          node [blackdot] {}
        };

      % constellation 6
      \node (X6) at (21,-8) {
        \begin{tikzpicture}
          [ level distance = 80px
          , bdot/.style = {circle, fill}
          , boundary/.style = {draw = none}
          , scale = 0.8
          ]
          \node (root_3) {} [grow' = up, boundary, level distance = 80px]
            child {
              node (omega) [bdot, label={left}:{\large$\omega$}] {} [sibling distance = 60]
                child{
                  node [boundary] {}
                  edge from parent node [left = 3, pos = 0.8] {\large$a_1$}
                }
                child{
                  node [boundary] {}
                  edge from parent node [left = 1, pos = 0.8] {\large$a_2$}
                }
                child{
                  node [boundary] {}
                  edge from parent node [right = 3, pos = 0.8] {\large$a_3$}
                }
              edge from parent node [right = 1, pos = 0.2] {\large$a_0$}
            };
        \end{tikzpicture}
      };

      \end{tikzpicture}

  }
% %
% \end{document}
  \captionof{figure}{Extended zoom complex}
  \label{fig:four_ope_zoom_ext}
\end{sidewaysfigure}

\begin{exm}
  The $4$-opetope introduced in \refexm{exm:four_opetope_zoom} may be extended as the zoom-complex depicted in \ref{fig:four_ope_zoom_ext}.
\end{exm}

\para In the following, we fix such an $n$-opetope $Y = \left(S_0 \to S_1 \to \cdots \to S_n\right)$ extended according to \refdefin{defin:zoom_extension}.
  We also write $\rho(S_1)$ as $\bot$.

\begin{defin}[$S_x$]\label{defin:nesting}
  Let $0 \le k \le n+1$, for $x \in S_{k+1}^\bullet \sqcup S_{k+1}^{\edge}$, we let $\mathcal{O}_{k,\,x} := \{y \in \msfdots(S_k') \mid y \sqsubseteq x\}$, where $y \sqsubseteq x$ means that there is a descending path from $y$ to $x$ in $S_{k+1}$.
  By the kernel rule, $\mathcal{O}_{k,\,x}$ induces a connected full subgraph of $S_k'$, that is a subdivided subtree of $S_k'$. We then let $S_{k,\,x}$ denote the underlying tree structure, and $S_{k,\,x}' = (S_{k,\,x},\,V_{k,\,x})$ the associated subdivided tree.
\end{defin}

\begin{defin}[$P(Y)$]\label{defin:ope_poset}
  We let $P(Y)$ be the following many-to-one poset.
  \begin{itemize}
    \item $P(Y)_{-1} = \{\bot\}$ and $P(Y)_k = S_{k+2}^{\edge}$ when $0 \le k \le n$.
    \item Let $\delta$ and $\gamma$ be defined by:
      \begin{itemize}
        \item[(i)] $\forall k\ge0,\, \forall x \in P(Y)_{k},\,\, \gamma(x) = \rho(S_{k+1,\,x})$.
        \item[(ii)] $\forall k>0,\, \forall x \in P(Y)_{k},\,\, \delta(x) = S_{k+1,\,x}^\vert$.
        \item[(iii)] $\forall x \in P(Y)_0,\,\, \delta(x) = \emptyset$.
      \end{itemize}
    \item When $z \prec^{\o} y \prec^{-} x \in P(Y)_k$ and $z \prec^{\o} y' \prec^{-} x$ for $y \neq y'$, then both $y,\,y'$ are distinct leaves of $S_{k+1,\,x}$, hence $\mathcal{O}_{k,\,y}$ and $\mathcal{O}_{k,\,y'}$ are disjoint sets of successive whitedots on the same edge $z$ of $S_k$.
      Hence they may be compared by the order on the subdivision of $z$.
      If $\mathcal{O}_{k,\,y} < \mathcal{O}_{k,\,y'}$, we let $y \blacktriangleleft y'$.
  \end{itemize}
\end{defin}

\begin{rmk}\label{rmk:ope_poset_loops}
  Notice that $\gamma(x) \prec^{\o} x \in P(Y)_k$ happens iff $k>0$ and $\mathcal{O}_{k+1,\,x} \subseteq \msfndots(S_{k+2})$,\\ that is iff $S_{k+1,\,x}$ is a unit tree. In particular, there is some $w \in S_{k+2}^\bullet$ with $\branch{x} w$.
\end{rmk}

\begin{rmk}\label{rmk:nodes_cofaces}
  When $x\in S_{k+1}^\bullet$ for some $k \le (n-1)$, then $x \in S_{k+2}^\vert$, hence $S_{k+1,\,x}$ is the corolla whose only node is $x$.\\
  Hence for $k > 0$, $\delta(x) = \{y \mid x \branch{y} \mbox{ in } S_{k+1}\}$ and for $k \ge 0$, $\gamma(x)$ is the target edge of $x$ in $S_{k+1}$.\\
  In particular, $\delta(x) \cap \gamma(x) = \emptyset$.
\end{rmk}

\begin{lem}\label{lem:ope_poset_mto}
  \refdefin{defin:ope_poset} indeed defines a many-to-one poset.
\end{lem}
\begin{proof}
  First, we see that $\gamma,\,\delta : P_{k+1} \to \mathcal{P}(P_k)$.
  By definition, we always have $\delta(x) \cup \gamma(x) \subseteq P(Y)_{\dim(x)-1}$ when $\dim(x)\ge1$.
  When $\dim(x) = 0$, $\delta(x) = \emptyset \subseteq P(Y)_{-1}$, and
  the only leaf of the linear tree $S_2$ must be in $\mathcal{O}_{1,\,x}$, which is the unique node of $S_1$.
  Hence $\gamma(x) = \rho(S_{1,\,x}) = \rho(S_1) = \bot$.\\[5pt]
  We then check the properties of many-to-one posets.
  \begin{itemize}
    \item When $\dim(x) \ge 0$, $\gamma(x)$ is given by (i). It consists of a unique element.
    \item Suppose $y \in \gamma(x) \cap \delta(x)$ for some $x$, then by \refrmk{rmk:ope_poset_loops}, $S_{k+1,\,x}$ is a unit tree whose only edge is $y$.
      Hence $\gamma(x) = \{\rho(S_{k+1,\,x})\} = \{y\} = S_{k+1,\,x}^\vert = \delta(x)$.
    \item By definition of $P(Y)$, $P(Y)_{-1} = \{\bot\}$. And by rule (iii), for all $x \in P(Y)_0,\, \bot \prec^{+} x$.
    \qedhere
  \end{itemize}
\end{proof}

\begin{prop}[Greatest element]\label{prop:greatest_element}
  $P(Y)$ admits $\top$ as a greatest element.
\end{prop}
\begin{proof}
  $\top$ is maximal since $P(Y)_{\ge n} = P(Y)_n = \{\top\}$.
  Let $y \in P(Y)_{k}$ for some $k<n$.
  \begin{itemize}
    \item If $S_{k+2}^{\edge} = \{y\}$, since $S_{k+3}$ cannot be empty, it must admit at least one edge $x$. 
      Then $S_{k+2,\,x}$ is a subtree of $S_{k+2}$, i.e. it is $S_{k+2}$, and $y \prec^{\o} x$.
    \item If $\branch{y} x$ for some $x$ in $S_{k+2}$, then $y = \gamma(x)$ by \refrmk{rmk:nodes_cofaces}.
    \item If $x \branch{y}$ for some $x$ in $S_{k+2}$, then $y \in \delta(x)$ by \refrmk{rmk:nodes_cofaces}.
  \end{itemize}
  Hence, for any $k<n$ and $y \in P(Y)_{k}$, there is an element $x \in P(Y)_{k+1}$ with $y \prec x$. It shows that $\top$ is a greatest element.
\end{proof}

\begin{prop}[Loops]\label{prop:loops}
  When $z \prec^{\o} y$ in $P(Y)$, there is some $x \in P(Y)$ with $y \prec^{+} x$.
\end{prop}
\begin{proof}
  Let $k>0$ and $z \prec^{\o} y \in P(Y)_{k}$. By \refrmk{rmk:ope_poset_loops}, $\branch{y} x$ for some $x$ in $S_{k+2}$, hence $y = \gamma(x)$. Since $x\in S_{k+2}^\bullet$, we have $\delta(x) \cap \gamma(x) = \emptyset$ by \refrmk{rmk:nodes_cofaces}, we have $y \nprec^{\o} x$.
\end{proof}

\begin{lem}[Source trees]\label{lem:source_trees}
  Let $x \in P(Y)_k$ for some $2 \le k \le n$.
  There is a rooted tree $\Sigma_x$ as follows.
  \begin{itemize}
    \item $\Sigma_x^\bullet = \{ y \in \delta(x) \mid S_{k,\,y} \mbox{ is not a unit tree} \}$. And $\Sigma_x^{\edge} = \{\gamma^2(x)\} \cup \partial\Sigma_x^\bullet$.
    \item $y \branch{z}$ in $\Sigma_x$ iff $z \in S_{k,\,y}^\vert$\,\,\quad and \quad $\branch{z} y$ in $\Sigma_x$ iff $z = \rho(S_{k,\,y})$.
    \item The root edge is $\rho(\Sigma_x) = \gamma^2(x)$.
  \end{itemize}
\end{lem}
\begin{proof}
  We check that it yields a well-defined rooted tree.
  \begin{itemize}
    \item Let $y \in \delta(x)$. Then its target edge is uniquely defined as $\rho(S_{k,\,y})$.
    \item Let $z \in \Sigma_x^{\edge} \setminus \gamma^2(x)$. Since the $\mathcal{S}_{k,\,y}^\bullet$ ($y\in\Sigma_x^\bullet$) form a partition of $\mathcal{S}_{k,\,\gamma(x)}^\bullet$ and $z \neq \gamma^2(x) = \rho(S_{k,\,\gamma(x)})$,\\
      $z$ admits a target node $y_0$ in $S_{k,\,\gamma(x)}$, and we have $y_0 \in \mathcal{S}_{k,\,y}^\bullet$ for a unique $y \in \Sigma_x^\bullet$.
      Hence $z \in S_{k,\,y}^\vert$ and $y \branch{z}$ in $\Sigma_x$.
    \item Since $\gamma^2(x)$ is the root edge of $S_{k,\,\gamma(x)}$, it may not be in $S_{k,\,y}^\vert$ for some $y \in \Sigma_x^\bullet$, as it would imply $\gamma^2(x)$ having a target node inside $\mathcal{S}_{k,\,\gamma(x)}^\bullet$.
      Hence, $\gamma^2(x)$ is targetless as an edge of $\Sigma_x$.
    \item Since the $\mathcal{S}_{k,\,y}^\bullet$ $(y \in \Sigma_x^\bullet)$ are non-empty and disjoint, we have $S_{k,\,y}^\vert \cap S_{k,\,y'}^\vert = \emptyset$ for $y \neq y'$ in $\Sigma_x^\bullet$. Whence the uniqueness for the target node of an edge $z \in \Sigma_x^{\edge}$.
    \item When there is a relation $\branched{z} y \branched{z'}$ in $\Sigma_x$, we have $z' \in S_{k,\,y}^\vert$ and $z' = \rho(S_{k,\,y})$.
      Since $S_{k,\,y}$ contains at least one node, we then find a path $\branched{z = z_1} y_1 \branched{z_2} \cdots \branched{z_p} y_p \branched{z_{p+1} = y'}$ in $S_k$ for some $p \ge 1$.
      Hence the length of descending paths in $\Sigma_x$ is bounded, which (using the second point) proves that there exists a descending path from every edge $z \in \Sigma_x^{\edge}$ to $\gamma^2(x)$. The uniqueness follows from the fourth point.
    \qedhere
  \end{itemize}
\end{proof}

\begin{prop}[Acyclicity]\label{prop:acyclicity}
  $P(Y)$ satisfies the property (acyclicity).
\end{prop}
\begin{proof}
  It is a direct consequence of the tree structure $\Sigma_x$ as given by \reflem{lem:source_trees}.
\end{proof}

\begin{lem}[Leaf lozenge]\label{lem:leaf_lozenge}
  Let $x \in P(Y)_{\ge2}$ where $x$ is not a loop. \\ There is a lozenge $(z \prec^{-} y \prec^{-} x,\, z \prec^{-} \gamma(x) \prec^{+} x)$ in $P(Y)$ iff $y \branch{z}$ is a leaf in $\Sigma_x$.
\end{lem}
\begin{proof}
  It follows from the observation $\Sigma_x^\vert = S_{k,\,\gamma(x)}^\vert$, which follows from the fact that the $\mathcal{S}_{k,\,y}^\bullet$ $(y \in \Sigma_x^\bullet)$ form a partition of $\mathcal{S}_{k,\,\gamma(x)}^\bullet$.
\end{proof}

\begin{lem}[Root lozenge]\label{lem:root_lozenge}
  Let $x \in P(Y)_{\ge2}$ where $x$ is not a loop. \\ There is a lozenge $(\gamma(y) \prec^{+} y \prec^{-} x,\, \gamma^2(x) \prec^{+} \gamma(x) \prec^{+} x)$ in $P(Y)$ iff $\branch{\gamma^2(x)} y$ in $\Sigma_x$.
\end{lem}
\begin{proof}
  It follows from the observation $\rho(\Sigma_x) = \rho(S_{k,\,y})$ for $y$ the root node of $\Sigma_x$, i.e. with $\branch{\rho(\Sigma_x)} y$ in $\Sigma_x$.
\end{proof}

\begin{lem}[Sign rule]\label{lem:sign_rule}
  Let $x \in P(Y)_{\ge2}$ where $x$ is not a loop. \\ There is no lozenge of the form $(z \prec^{\alpha} y \prec^{-} x,\, z \prec^{\alpha} y' \prec^{-} x)$ in $P(Y)$  for $y \neq y'$ and $\alpha \in \{+,\,-\}$.
\end{lem}
\begin{proof}
  By observing that the sets $\mathcal{S}_{k,\,y}^\bullet$ $(y \in \Sigma_x^\bullet)$ are disjoint and non-empty.
\end{proof}

\begin{lem}[loop lozenge]\label{lem:loop_lozenge}
  For $(z \prec^{\o} y \prec^{-} x \in P(Y)_{\ge2})$, either there is some $y'$ and $\alpha \neq \o$ such that $(z \prec^{\alpha} y' \prec^{-} x)$, or there is a chain $(z \prec^{\o} \gamma(x) \prec^{+} x)$.
\end{lem}
\begin{proof}
  Suppose that $\gamma(x)$ is a loop, then $S_{k,\,\gamma(x)}$ is an empty tree whose edge is $\gamma^2(x)$. Since $y \in S_{k+1,\,x}^\vert$, we have $\mathcal{O}_{k,\,y} \subseteq \mathcal{O}_{k,\,\gamma(x)}$, hence $S_{k,\,y}$ must also be a unit tree, with $z = \gamma^2(x)$ as only edge. Whence the chain $z \prec^{\o} \gamma(x) \prec^{+} x$.\\
  Suppose now that $\gamma(x)$ is not a loop. Then $S_{k,\,\gamma(x)}^\bullet \neq \emptyset$. 
  \begin{itemize}
    \item Either $z$ admits a source node $y'$ in $\Sigma_x$, whence $z \prec^{+} y' \prec^{-} x$ by definition of $\Sigma_x$.
    \item Or $z$ admits a target node $y'$, whence $z \prec^{-} y' \prec^{-} x$ by definition of $\Sigma_x$.
    \qedhere
  \end{itemize}
\end{proof}

\begin{prop}[Oriented thinness]\label{prop:oriented_thinness}
  $P(Y)$ satisfies the property (oriented thinness).
\end{prop}
\begin{proof}
  For chains of the form $(z \prec^{\beta} y \prec^{\alpha} x)$ with $x \in P(Y)_{\ge{2}}$, this is shown by Lemmas \ref{lem:leaf_lozenge}, \ref{lem:root_lozenge}, \ref{lem:sign_rule}, \ref{prop:loop_lozenge} and using the tree structure of $\Sigma_x$.\\
  Consider $x \in P(Y)_1 = S_3^{\edge}$ which is not a loop, then $S_{2,\,x}$, as a (non-trivial) subtree of the linear tree $S_2$, must contain exactly one leaf and one root (distinct from each other).
  Hence $x$ admits exactly one source $y \prec^{-} x$ and one target $y' \prec^{+} x$, whence the lozenge completion for the chains of the form $\bot \prec y \prec x$.
\end{proof}

\begin{thm}\label{thm:ope_poset_is_a_dfc}
  The many-to-one poset $P(Y)$ is a dendritic face complex.
\end{thm}
\begin{proof}
  By Propositions \ref{prop:greatest_element}, \ref{prop:loops}, \ref{prop:acyclicity} and \ref{prop:oriented_thinness}.
\end{proof}

\begin{defin}\label{defin:ope_poset_morphism}
  Let $f : Y \to Z$ be an isomorphism of opetopes. There is an isomorphism $P(f) : P(Y) \to P(Z)$ defined as $f$ on every element.
\end{defin}
\begin{proof}
  First, denote $Y = (S_i)_{i\le n+2}$ and $Z = (T_i)_{i\le n+2}$ (by extending the respective opetopes according to \ref{defin:zoom_extension}, and extending $f$ accordingly).
  Then for any $k \ge 0$, $f$ sends $x \in S_{k+2}^{\edge} = P(Y)_k$ to $f(x) \in T_{k+2}^{\edge} = P(Z)_k$.
  Since $f$ is an isomorphism of trees, we have $\mathcal{O}_x \overset{f}\simeq \mathcal{O}_{f(x)}$ for any $x \in S_{k+2}^{\edge}$. Hence $\gamma(x) \overset{f}{\simeq} \gamma(f(x))$ and $\delta(x) \overset{f}{\simeq} \delta(f(x))$. Hence $P(f)$ preserves also $\prec^{\alpha}$ for $\alpha \in \{+,\,-,\,\o\}$.
\end{proof}

\begin{thm}[$\KJBM \to \DFC^\core$]
  There is a functor $P : \KJBM \to \DFC^\core$ given by \refdefin{defin:ope_poset} on objects and \refdefin{defin:ope_poset_morphism} on isomorphisms.  
\end{thm}
\begin{proof}
  The functoriality is clear from \refdefin{defin:ope_poset_morphism}
\end{proof}

\section{An equivalence of categories}
\label{sec:cat_equiv}

\para In this section, we show that $Z$ and $P$ form an equivalence of categories.

\begin{defin}[$\theta$]\label{defin:PZ_trans}
  Let $C$ be a DFC of dimension $n \ge 0$, then by definition $P(Z(C))$ has same dimension as $C$.
  We let $\theta_C : C \to (P \circ Z)(C)$ be the following natural isomorphism of DFC:
  \begin{itemize}
    \item $\theta_C$ sends the bottom element to the bottom element.
    \item For $x \in C_k \,(0 \le k < n)$, $\theta_C(x) := x \in P(Z(C))_k = S_{k+2}^{\edge} = C_k$.
  \end{itemize}
\end{defin}
\begin{prop}\label{prop:PZ_trans_well_defined}
  The natural isomorphism $\theta$ from \refdefin{defin:PZ_trans} is well defined.
\end{prop}
\begin{proof}
  That $\theta_C$ is natural in $\DFC^\core$ is clear.
  We check that it is a well-defined morphism of DFC.
  \begin{itemize}
    \item First, we see $\forall x \in C,\, \delta(x) = \delta(\theta_C(x))$.\\
      Let $x \in C_k$, we have $\delta(\theta_C(x)) = S_{k+1,\,\theta_C(x)}^\vert$.
      \begin{itemize}
        \item When $k = n$, we have $x = \omega$ and $\theta_C(\omega) = \top$. Hence $\delta(x) = C_{n-1} \setminus \{\gamma(\omega)\}$ and $\delta(\theta_C(x)) = S_{n+1,\,\top} = \Lambda_n = \delta(\omega) \setminus \gamma(\top)$.
        \item When $0 < k < n$, we have $x = \theta_C(x)$, and using \reflem{lem:facet_of_Lambda}, $S_{k+1,\,x}^\bullet = \{x' \sqsubseteq x \mid x' \in \Lambda_k\}$, where $x' \sqsubseteq x$ means that there is a descending path from $x'$ to $x$ in $S_{k+2}$.
        
          Let $y \in \delta(x)$. Then, repeatedly applying oriented thinness and \reflem{lem:facet_of_Lambda} yields the following pattern in $C$, where $w_1,\,\cdots,\,w_p,\,x' \in \Lambda$:
          % https://q.uiver.app/#q=WzAsOCxbMCwxLCJ4Il0sWzEsMCwid18xIl0sWzIsMSwiXFxidWxsZXQiXSxbMywwLCJ3XzIiXSxbNCwxLCJcXGJ1bGxldCJdLFs1LDAsIndfcCJdLFs2LDEsIngnIl0sWzMsMiwieSJdLFsxLDAsIiIsMCx7InN0eWxlIjp7ImJvZHkiOnsibmFtZSI6ImJhcnJlZCJ9LCJoZWFkIjp7Im5hbWUiOiJub25lIn19fV0sWzEsMiwiIiwyLHsic3R5bGUiOnsiaGVhZCI6eyJuYW1lIjoibm9uZSJ9fX1dLFsyLDMsIiIsMix7InN0eWxlIjp7ImJvZHkiOnsibmFtZSI6ImJhcnJlZCJ9LCJoZWFkIjp7Im5hbWUiOiJub25lIn19fV0sWzMsNCwiIiwyLHsic3R5bGUiOnsiYm9keSI6eyJuYW1lIjoiZG90dGVkIn0sImhlYWQiOnsibmFtZSI6Im5vbmUifX19XSxbNCw1LCIiLDIseyJzdHlsZSI6eyJib2R5Ijp7Im5hbWUiOiJiYXJyZWQifSwiaGVhZCI6eyJuYW1lIjoibm9uZSJ9fX1dLFs1LDYsIiIsMix7InN0eWxlIjp7ImhlYWQiOnsibmFtZSI6Im5vbmUifX19XSxbMCw3LCIiLDAseyJjdXJ2ZSI6Miwic3R5bGUiOnsiaGVhZCI6eyJuYW1lIjoibm9uZSJ9fX1dLFsyLDcsIiIsMCx7InN0eWxlIjp7ImhlYWQiOnsibmFtZSI6Im5vbmUifX19XSxbNCw3LCIiLDAseyJzdHlsZSI6eyJoZWFkIjp7Im5hbWUiOiJub25lIn19fV0sWzYsNywiIiwyLHsiY3VydmUiOi0yLCJzdHlsZSI6eyJoZWFkIjp7Im5hbWUiOiJub25lIn19fV1d
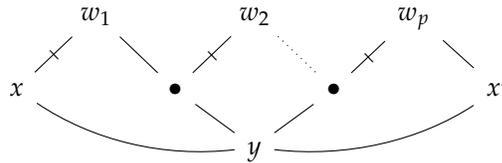
          \[\begin{tikzcd}[sep = small]
            & {w_1} && {w_2} && {w_p} \\
            x && \bullet && \bullet && {x'} \\
            &&& y
            \arrow["\shortmid"{marking}, no head, from=1-2, to=2-1]
            \arrow[no head, from=1-2, to=2-3]
            \arrow["\shortmid"{marking}, no head, from=2-3, to=1-4]
            \arrow[dotted, no head, from=1-4, to=2-5]
            \arrow["\shortmid"{marking}, no head, from=2-5, to=1-6]
            \arrow[no head, from=1-6, to=2-7]
            \arrow[curve={height=12pt}, no head, from=2-1, to=3-4]
            \arrow[no head, from=2-3, to=3-4]
            \arrow[no head, from=2-5, to=3-4]
            \arrow[curve={height=-12pt}, no head, from=2-7, to=3-4]
          \end{tikzcd}\]
          Hence $x' \branch{y}$ in $S_{k+1,\,x}$.
          Suppose now having a relation $\branch{y} x''$ in $S_{k+1,\,x}$, then there is a path $x \prec^{+} w_1' \succ^{-} \cdots \succ^{-} \bullet \prec^{+} w_q' \succ^{-} x'' \in \Lambda$.
          And using successively oriented thinness as above yields the relation $y \prec^{-} x''$, which is absurd. Whence $y \in S_{k+1,\,x}^\vert$.
          
          Conversely, let $y \in S_{k+1,\,x}^\vert$. Then $y \prec^{-} x'$ for some $x' \in \Lambda_k,\, x' \sqsubseteq x$. Hence there is a path $x = x_0 \prec^{+} w_1 \succ^{-} \cdots \succ^{-} x_{p-1} \prec^{+} w_p \succ^{-} x' = x_p$, and $y \prec^{-} x_i$ for all $i$'s. Indeed suppose there is some $i$ s.t. $y \nprec^{-} x_i$ and take a minimal such index $i_0$.
          Then the completion of the chain $y \prec^{-} x_{i_0+1} \prec^{-} w_{i_0+1}$ must yield $y \prec^{+} x_1' \prec^{-} w_{i_0+1}$ and by repeatedly aplying oriented thinness and \reflem{lem:facet_of_Lambda} as above we may find some $x'' \sqsubseteq x,\, x'' \in \Lambda_{k}$ with $y \prec^{+} x''$.
          Hence there is a configuration $x' \branch{y} x''$ in $S_{k+1,\,x}$, which is absurd.
          Whence $y \prec^{-} x_0 = x$, i.e. $y \in \delta(x)$.
      \end{itemize}

    \item With a similar reasoning, we show that for $x \in C$, $\gamma(\theta_C(x)) = \rho(S_{k+1,\,\theta_C(x)})$, whence $\gamma(x) = \gamma(\theta_C(x))$.
  \end{itemize}
  It proves that $\theta_C$ preserves $\delta$ and $\gamma$, hence it preserves $\prec^{+},\,\prec^{-}$ and $\prec^{\o}$.
    
  We finally check that $\theta_C$ preserves the order $\blacktriangleleft$. Let $z \prec^{\o} y \prec^{-} x \in \Lambda_k$ and $z \prec^{\o} y' \prec^{-} x$ with $y \blacktriangleleft y'$ in $C$. Then for any loop $\bar{y}$ (resp. $\bar{y}'$) included in or equal to $y$ (resp. $y'$), we have $\bar{y} \blacktriangleleft \bar{y'}$ according to \refdefin{defin:loops_total_order}.
  Hence $\mathcal{O}_{k,\,y} \blacktriangleleft \mathcal{O}_{k,\,y'}$, whence $y \blacktriangleleft y'$ in $(P\circ Z)(C)$. Hence $\theta_C$ is an isomorphism of DFC.
\end{proof}

\begin{prop}[$\tau$]\label{prop:ZP_trans}
  Let $Y = (S_i)_{0\le i\le n}$ be an $n$-opetope for some $n\ge0$. Then by definition $(Z \circ P)(Y)$ is an $n$-opetope and there is an isomorphism $\tau_Y : Y \to (Z \circ P)(Y)$, natural in $Y$, as follows.
  \begin{itemize}
    \item Denote $(Z \circ P)(Y)$ as $(T_i)_{0 \le i \le n}$ (as given by \refprop{prop:DFC_tree_structure_Lambda}).
    \item Define $S_{n+1},\, S_{n+2}$ according to \ref{defin:zoom_extension}. And consider $T_{n+1},\, T_{n+2}$ as in \refprop{prop:DFC_tree_structure_Lambda}.
    \item Define $\tau_{Y,\,n+2} : S_{n+2} \to T_{n+2}$ in the only possible way (both are unit trees).
    \item For $2 \le i \le n+1$, let $\tau_{Y,\,k} : S_k \to T_k$ be defined as the identity on all elements.
    \item Extend $\tau_C$ in degree $0$ and $1$ in the only possible way.
  \end{itemize} 
\end{prop}
\begin{proof}
  The naturality is clear from the definition. We check that $\tau_C$ is well defined.
  \begin{itemize}
    \item First, we need to show that $S_i^\bullet = \Lambda_{i-1} = T_i^\bullet$ for $2 \le i \le n$.
      \begin{itemize}
        \item Let $y \in S_i^\bullet$ and suppose $y \prec^{+} x$ in $P(Y)$. Hence $y = \rho(S_{i+1,\,x})$. Since $y \nprec^{\o} x$, $S_{i+1,\,x}$ contains at least one node, and $\branch{y} x'$ where $x'$ is the root node of $S_{i+1,\,x}$.
          Hence $y$ is not a leaf in $S_i$.
        \item Let $y \in T_i^\bullet = \Lambda_i$ and suppose $y$ is not a leaf in $S_{i+2}$. Then $\branch{y} x$ for some $x \in S_{i+2}^\bullet$ and we have  $\{x\} = S_{i+1,\,x}^\bullet$. Hence $y = \rho(S_{i+1,\,x})$, i.e. $y \prec^{+} x$, which is absurd.
      \end{itemize}
    \item We have $T_i^{\edge} = P(Y)_{i-2} = S_i^{\edge}$ for $2 \le i \le n$.
    \item Let $1 \le i \le n$, We now check that $\tau_{C,\,i}$ preserves the relations $x \branch{y}$ and $\branch{y} x$.
      \begin{itemize}
        \item Suppose $x \branch{y}$ in $S_i$, then $y \prec^{-} x$ because $S_{i+1,\,x}^\bullet = \{x\}$.
          Since $x \in S_i^\bullet$, $x = \Lambda_i = T_i^\bullet$ by the first point, hence by definition of $T_i$, we have $x \branch{y}$ in $T_i$.
        \item The case $\branch{y} x$ is similar.
        \qedhere
      \end{itemize}
  \end{itemize}
\end{proof}

\begin{thm}[$\DFC^\core \simeq \KJBM$]
  There is an equivalence of categories
  % https://q.uiver.app/#q=WzAsMixbMCwwLCJcXERGQyJdLFszLDAsIlxcT3BlIl0sWzAsMSwiWiIsMSx7ImN1cnZlIjotNH1dLFsxLDAsIlAiLDEseyJjdXJ2ZSI6LTR9XV0=
  \[\begin{tikzcd}
    \DFC^\core &&& \KJBM
    \arrow["Z"{description}, curve={height=-24pt}, from=1-1, to=1-4]
    \arrow["P"{description}, curve={height=-24pt}, from=1-4, to=1-1]
  \end{tikzcd}\]
\end{thm}
\begin{proof}
  This is by \refdefin{defin:PZ_trans} and \refprop{prop:ZP_trans}.
\end{proof}

\section{Conclusion and related works}

In this article, we have presented a new way to encode opetopes, and shown our formalism to be equivalent to the one of zoom complexes presented in \cite{Kock2010} through an equivalence of categories
$\DFC^\core \simeq \KJBM$.
As a consequence, and because of the other equivalence obtained in \cite{leclerc2023}, it gives a embeding of positives opetopes as described by \textsc{Zawadowski} in the category $\KJBM$ defined in \refdefin{defin:op_cat}.

A further work about DFC would be to figure out whether the whole formalism of DFC for non-necessarily positive opetopes is equivalent to the formalism of ordered face structures introduced by \textsc{Zawadowski} in \cite{zawadowski2008}.
The present paper has already done a step toward this goal, and it would not be surprising that those two formalisms are indeed equivalent.

\section{Acknowledgements}

I am grateful to Marek \textsc{Zawadowski} for our discussions about opetopes, and for offering me ideas for writing this article.
I also express a special thank you to the supervisor of my master's thesis, Pierre Louis \textsc{Curien} who has accompanied me for an entire year of study of higher categories and opetopes.

  % - Bibliography
  \bibliographystyle{plain}
  \bibliography{main.bib}

\end{document}